\documentclass[notitlepage,leqno,10pt]{article}

\usepackage{latexsym}
\usepackage{amsmath}
\usepackage{amsfonts}
\usepackage{amssymb}
\usepackage{graphicx}
\renewcommand{\a }{\alpha }
\renewcommand{\d}{\delta }
\newcommand{\D }{\Delta }

\newcommand{\e }{\varepsilon }
\newcommand{\g }{\gamma}

\newcommand{\G }{\Gamma }
\renewcommand{\l }{\lambda }
\renewcommand{\L }{\Lambda }

\newcommand{\n }{\nabla }
\newcommand{\var }{\varphi }

\newcommand{\s }{\sigma }

\renewcommand{\t }{\tau }
\renewcommand{\th }{\theta }
\renewcommand{\o }{\omega }

\renewcommand{\O }{\Omega }

\newcommand{\ov}{\overline}
\newcommand{\intbar}{\mathop{\int\makebox(-13.5,0){\rule[4pt]{.7em}{0.3pt}}%
\kern-6pt}\nolimits}

\newcommand{\be}{\begin{equation}}
\newcommand{\ee}{\end{equation}}
\newenvironment{pf}{\noindent{\sc Proof}.\enspace}{\rule{2mm}{2mm}\medskip}

\newtheorem{corollary}{Corollary}[section]

\newcommand{\R}{\mathbb{R}}
\newcommand{\C}{\mathbb{C}}
\renewcommand{\H}{\mathbb{H}}

\newcommand{\Z}{\mathbb{Z}}

\newcommand{\N}{\mathbb{N}}
\newcommand{\pa}{\partial}

\newcommand{\ou}{\overline{1}}

\def\func#1{\mathop{\rm #1}\nolimits}%

\author{}

\date{}

\title{A positive mass theorem in three dimensional \\ Cauchy-Riemann geometry}
\author{Jih-Hsin Cheng$^{(a)}$, Andrea Malchiodi$^{(b)}$, and Paul Yang$^{(c)}$}

\begin{document}

\newtheorem{lem}{Lemma}[section]
\newtheorem{pro}[lem]{Proposition}
\newtheorem{thm}[lem]{Theorem}
\newtheorem{rem}[lem]{Remark}
\newtheorem{cor}[lem]{Corollary}
\newtheorem{df}[lem]{Definition}

\maketitle

\begin{center}

{\small }

\end{center}

\footnotetext[1]{E-mail address: cheng@math.sinica.edu.tw, A.Malchiodi@warwick.ac.uk, malchiod@sissa.it, yang@math.princeton.edu}

\vspace{-1cm}

\begin{center}

{\small $^{(a)}$ Institute of Mathematics, Academia Sinica and NCTS, Taipei office

6F, Astronomy-Mathematics Building
No. 1, Sec. 4

Roosevelt Road,
Taipei 10617, TAIWAN

\

$^{(b)}$ University of Warwick, Mathematics Institute - Zeeman Building 
Coventry CV4 7AL \\ and  SISSA, Via Bonomea 265, 34136 Trieste, ITALY

\

$^{(c)}$Princeton University, Department of Mathematics

Fine Hall, Washington Road, 
Princeton NJ 08544-1000 USA

}

\end{center}

\

\

\noindent {\sc abstract}. We define an ADM-like mass, called p-mass, for an
asymptotically flat pseudohermitian manifold. The p-mass for the blow-up of a
compact pseudohermitian manifold (with no boundary) is identified with the first
nontrivial coefficient
in the expansion of the Green function for the CR Laplacian. We deduce an
integral formula for the p-mass, and we reduce its positivity to a solution
of Kohn's equation. We prove that the p-mass is non-negative for (blow-ups of) compact
3-manifolds of positive Tanaka-Webster class and with non-negative CR Paneitz operator.
Under these assumptions, we also characterize the zero mass case as the
standard three dimensional CR sphere. We then show the existence of (non-embeddable) CR
3-manifolds having nonpositive Paneitz  operator or negative p-mass through
a second variation formula. Finally, we apply our main result to find solutions
of the CR Yamabe problem with minimal energy. \

\begin{center}

\bigskip


\bigskip

\centerline{\bf AMS subject classification: 32V20, 53C17, 35J75, 35J20, 32V30}

\end{center}

\

\begin{center}
{\Large Table of Contents}\bigskip
\end{center}

\ \ \ \ \ \ \ \ \ \ \ \ \ {\large 1. Introduction and statement of the
results}

\ \ \ \ \ \ \ \ \ \ \ \ \ {\large 2. Asymptotically flat pseudohermitian manifolds and definition of the  p-mass}

{\large \ \ \ \ \ \ \ \ \ \ \ 3. Proof of Theorem \ref{t:pm}}

{\large \ \ \ \ \ \ \ \ \ \ \ 4. Some examples}

{\large \ \ \ \ \ \ \ \ \ \ \ 5. Proof of Theorem \ref{t:y}}

%{\large \ \ \ \ \ \ \ \ \ \ \ 7. The Yamabe quotient in the locally spherical case}

{\large \ \ \ \ \ \ \ \ \ \ \ 6. Appendix: useful facts in pseudohermitian geometry}

\bigskip

\section{Introduction and statement of the
results}\label{s:in}

Around 1960 the three physicists R.Arnowitt, S.Deser, and C.Misner developed the
Hamiltonian formulation for Einstein's general theory of relativity in a
series of papers, see e.g. \cite{ADM1}, \cite{ADM2}. The total energy was called
the \emph{ADM-mass} later in the literature. The positive mass conjecture asserts
the non-negativity of such a quantity assuming non-negativity of the local
energy density, and classifies the situation for the
zero mass case. The conjecture was proved by R.Schoen and S.T.Yau using minimal
surface theory (\cite{SY1}, \cite{SY2}, \cite{SY3}, \cite{SY4}) and by Witten
using the Dirac equation for spinors (\cite{Wi}; a rigorous proof is given in \cite{PT})
in late 70's and early 80's. Since there is a strong analogy
between CR geometry and conformal geometry, one may wonder about what happens in
CR geometry correspondingly. This is the main goal of this paper.

We recall some well known definitions and facts, referring to \cite{LP}
for more details. A Riemannian manifold $N$ (of dimension $n$) with a
smooth metric $g$ is called {\em asymptotically flat} of order $\t > 0$ if $N$
decomposes as $N = N_0 \cup N_\infty$, where $N_0$ is compact and $N_\infty$ is
diffeomorphic to $\R^n \setminus B_R(0)$ for some $R > 0$, and if $g$
on $N_\infty$ satisfies
$g_{ij} = \d_{ij} + O(|x|^{-\t})$, $|\pa g_{ij}| = O(|x|^{-\t-1})$,
$|\pa^2 g_{ij}| = O(|x|^{-\t-2})$ as $|x| \to + \infty$, in some system of coordinates.
Such coordinates are called {\em asymptotic coordinates}. One then considers
the Einstein-Hilbert action
$$
  \mathcal{A}(g) := - \int_N S_g dV_g,
$$
where $S_g$ is the scalar curvature of $(N,g)$. Given a smooth variation
$g_s$ of the metric with $\frac{d}{ds}|_{s=0} g_s = v$, $v = (v_{ij})$,
one has
$$
  \frac{d}{ds}|_{s=0} (S_{g_s} dV_{g_s}) = - \left( v^{ij} G_{ij} + \nabla^*
  \xi \right) dV_{g},
$$
where $G_{ij} = R_{ij} - \frac 12 S_g g_{ij}$
is the Einstein tensor and where $\xi = (v_{jk,}^{\; \; \; \; \; k} - v^k_{k,j}) dx^j$.
Then, if one defines the {\em mass} $m(g)$ to be
\begin{equation}\label{eq:massriem}
    m(g) = \lim_{R \to + \infty} \o_n^{-1} \oint_{\Sigma_R}
    \mu \lrcorner dV_g \qquad \qquad \mu=(\partial_{i}g_{ij}-
   \partial_{j}g_{ii})\partial_{j}, \quad \Sigma_R = \{ |x| = R \}.
\end{equation}
varying the metric it turns out that
\begin{equation}\label{eq:varmass}
    \frac{d}{ds} (\mathcal{A}(g_s) + m(g_s)) = \int_N v^{ij} G_{ij} dV_g.
\end{equation}
The positive mass theorem was then used by Schoen in \cite{Sc} to prove the
Yamabe conjecture (finding conformal metrics with constant scalar curvature) in
the cases left open in 1976 by T.Aubin (\cite{aub}), namely
in dimension 3, 4, 5 or when the metric is locally conformally flat. As shown in
\cite{trud}, to solve the problem (in the {\em positive case}, the most difficult
one) it is sufficient to show that the Yamabe quotient of the manifold is smaller
than the one of $\R^n$. This was done both in \cite{aub} and in \cite{Sc} using
suitable test functions. In particular Schoen's argument consists in gluing
a {\em standard bubble} near any $p \in M$ to $\mathcal{G}_p$, the Green's function
of the conformal Laplacian with pole at $p$. It turns out that the regular part
of $\mathcal{G}_p$ coincides with the mass of $(M \setminus \{p\},
\mathcal{G}_p^{\frac{4}{n-2}} g)$, the {\em blow-up} of $M$ at $p$.
The positive mass theorem  played a key role in
other problems in conformal geometry, for example in the compactness of the
Yamabe equation or the Yamabe flow, see \cite{BR1}, \cite{BR2}, \cite{KMS}.

\

\noindent The goal of this paper is to introduce a notion of asymptotically flat
pseudohermitian manifold, define a {\em pseudohermitian mass} and to prove its positivity
under suitable conditions.
We give then applications concerning existence of minimizers for the Tanaka-Webster quotient,
which is useful to study the CR analogue of the Yamabe quotient.

We consider a compact three dimensional pseudohermitian manifold $(M,J,\th)$ (with no boundary)
of {\em positive Tanaka-Webster
class}. This means that the first eigenvalue of the {\em conformal sublaplacian}
$$
  L_b := - 4 \D_b + R,
$$
is strictly positive. Here $\D_b$ stands for the sublaplacian of $M$ and $R$
for the Tanaka-Webster curvature, see Subsection \ref{ss:notprel} for their
definition. The conformal sublaplacian has the following covariance property
under a conformal change of contact form
$$
  \hat{L}_b (\var) = u^{-\frac{Q+2}{Q-2}} L_b(u \var); \qquad
  \quad \hat{\th} = u^2 \th,
$$
where $Q = 4$ is the {\em homogeneous dimension} of the manifold. The conformal
sublaplacian rules the change of the Tanaka-Webster curvature under
the above conformal deformation, through the following formula
$$
  - 4 \D_b u + R u = \hat{R} u^{\frac{Q+2}{Q-2}},
$$
where $\hat{R}$ is the Tanaka-Webster curvature corresponding to the pseudohermitian structure
$(J, \hat{\th})$. The positivity of the Tanaka-Webster class is  equivalent to the condition
\begin{equation}\label{eq:Y(J)}
    \mathcal{Y}(J) := \inf_{\hat{\th}} \frac{\int_M R_{J,\hat{\th}}
  \hat{\th} \wedge  d \hat{\th}}{\left( \int_M \hat{\th}
\wedge d \hat{\th} \right)^{\frac 12}} > 0,
\end{equation}
where $\hat{\th}$ is any contact form which annihilates $\xi$ (the contact bundle).
Under the assumption $ \mathcal{Y}(J) > 0$ we have that $L_b$ is
invertible so for any $p \in M$ there exists
a Green's function $G_p$ for which
$$
  \left( - 4 \D_b + R \right) G_p = 16 \d_p.
$$
One can show that in CR normal coordinates $(z,t)$ (see Subsection \ref{ss:CRcoord})
$G_p$ admits the following expansion
$$
  G_p = \frac{1}{2\pi} \rho^{-2} + A + O(\rho),
$$
where $A$ is some real constant and where we have set $\rho^4(z,t) =
|z|^4 + t^2$, $z \in \C, t \in \R$. Having in mind the Riemannian construction
for the blow-up of a compact manifold, we consider in Subsection \ref{ss:bu}
the new pseudohermitian manifold with a blow-up of contact form
\begin{equation}\label{eq:bbuu}
    N = (M \setminus \{p\}, J, \th = G_p^2 \hat{\th}),
\end{equation}
where $\hat{\th}$ is chosen so that near $p$ it has the behavior as described in Proposition
\ref{p:CRcoord}. With an {\em inversion of coordinates} (described in Subsection \ref{ss:inv})
we then obtain a  pseudohermitian manifold which has asymptotically  the
geometry of the Heisenberg group. Starting from this model, in Subsection \ref{ss:asyint}
we give a definition of asymptotically flat pseudohermitian manifold and we
introduce its {\em pseudohermitian mass} (p-mass) by the formula
$$
  m(J, \th) := i \oint_{\infty} \o^1_1 \wedge \th :=
  \lim_{\L \to + \infty} i \oint_{S_\L} \o^1_1 \wedge \th,
$$
where we have set $S_\L = \left\{ \rho = \L \right\}$, $\rho^4 = |z|^4 + t^2$, and
where $\o^1_1$ stands for the connection form of the structure, see Subsection
\ref{ss:notprel}. The above quantity is indeed a natural candidate, since it satisfies a property
analogous to \eqref{eq:varmass}, see Remark \ref{r:mass}, and moreover it
coincides with the zero-th
order term in the expansion of the Green's function for $L_b$, see Lemma \ref{l:m=A}.

In Subsection \ref{ss:asyint} we prove an integral formula for the p-mass, in the
spirit of \cite{Wi}. To state this formula we need to introduce another conformally
covariant operator, the CR Paneitz operator
$$
  P \varphi :=4(\varphi {_{\bar{1}}}^{\bar{1}}{_{1}}+iA_{11}\varphi
^{1})^{1} 
$$
(where we are following the notation of Subsection \ref{ss:notprel}). The operator
$P$ satisfies the covariance property (\cite{Hi})
\begin{equation}\label{eq:transfP}
      P_{(J, \hat{\th})} = u^{-4} P_{(J,\th)}; \qquad \quad \hat{\th} = u^2 \th.
\end{equation}
In Proposition \ref{p:genfor} we prove then the following integral formula,
which holds for an asymptotically flat pseudohermitian manifold $N$
\begin{equation}\label{eq:massintintr}
    \frac{2}{3} m(J,\th) = - \int_N |\Box_b \beta|^2 \th \wedge d \th +
    2 \int_N |\beta_{, \ou \ou}|^2 \th \wedge d \th  + 2 \int_N R |\beta_{, \ou}|^2
    \th \wedge d \th + \frac 12 \int_N \ov{\beta} P \beta \, \th \wedge d \th.
\end{equation}
Here $\beta : N \to \C$ is a function satisfying
$$
  \beta = \ov{z} + \beta_{-1} + O(\rho^{-2+\e}) \quad \hbox{ near }
     \infty; \qquad  \qquad \Box_b \beta = O(\rho^{-4}),
$$
with $\Box_b = - 2 \beta_{, \ou 1}$ and with $\beta_{-1}$ a suitable
function with homogeneity $-1$ in $\rho$ (satisfying condition
\eqref{eq:assasybeta} below).

\

\noindent To explain the link between \cite{Wi} and the $\Box_b$ operator, 
let us consider a general spin$^{c}$ structure on a contact bundle $\xi $
over an asymptotically flat pseudohermitian manifold $M$ of dimension $2n+1.$
Let $W$ denote the spinor bundle with a spin$^{c}$ connection $\nabla $
compatible with the pseudohermitian connection $\nabla ^{p.h.}.$ Let $e_{1},$
$...,$ $e_{2n}$ denote an orthonormal basis of $\xi $ with respect to the
Levi metric. Denote the contact Dirac operator $D_{\xi }$ by%
\[
D_{\xi }\psi =\sum_{\alpha =1}^{2n}\Gamma (e_{\alpha })\nabla _{e_{\alpha
}}\psi 
\]

\noindent (summation convention) for a section $\psi $ of $W,$ where $\Gamma 
$ denotes the Clifford multiplication. Let $T$ denote the Reeb vector field
associated to the contact form $\theta .$ Let $D_{\xi }^{\ast },$ $\nabla
^{\ast }$ denote the adjoint operator of $D_{\xi },$ $\nabla ,$ resp.. We
then have the following formula%
\begin{eqnarray}
D_{\xi }^{\ast }D_{\xi }\psi &=&\sum_{\alpha =1}^{2n}\nabla _{e_{\alpha
}}^{\ast }\nabla _{e_{\alpha }}\psi -2\sum_{\alpha =1}^{n}\Gamma (e_{\alpha
}e_{n+\alpha })\nabla _{T}\psi  \label{c1} \\
&&+\sum_{a<b}\Gamma (e_{a})\Gamma (e_{b})R^{\nabla }(e_{a},e_{b})\psi 
\nonumber
\end{eqnarray}

\noindent where $R^{\nabla }(e_{a},e_{b})$, the curvature operator, is
defined by%
\[
R^{\nabla }(e_{a},e_{b}):=\nabla _{e_{a}}\nabla _{e_{b}}-\nabla
_{e_{b}}\nabla _{e_{a}}-\nabla _{\lbrack e_{a},e_{b}]}.
\]

\noindent Decompose $R_{ab}^{\nabla }$ := $R^{\nabla }(e_{a},e_{b})$ as a
sum of the trace free part $\mathring{R}_{ab}^{\nabla }$ and the trace part $%
2^{-n}tr_{W}R_{ab}^{\nabla }.$ A standard deduction shows that%
\[
\sum_{a<b}\Gamma (e_{a})\Gamma (e_{b})\mathring{R}_{ab}^{\nabla }\psi =\frac{%
1}{4}R\psi
\]

\noindent where, again, $R$ denotes the Tanaka-Webster scalar curvature. On the
other hand, $2^{-n}tr_{W}R_{ab}^{\nabla }$ $=$ $F_{A}(e_{a},e_{b})$ in which
the curvature 2-form $F_{A}$ :$=$ $dA$ and $2A$ is the connection form of an
associated line bundle $L_{\Gamma }$ ($\det W$ $=$ $L_{\Gamma }^{\otimes
2^{n-1}}).$ Therefore we can reduce (\ref{c1}) to 
\begin{eqnarray}
D_{\xi }^{\ast }D_{\xi }\psi &=&\sum_{\alpha =1}^{2n}\nabla _{e_{\alpha
}}^{\ast }\nabla _{e_{\alpha }}\psi -2\sum_{\alpha =1}^{n}\Gamma (e_{\alpha
}e_{n+\alpha })\nabla _{T}\psi  \label{c2} \\
&&+\frac{1}{4}R\psi +\rho (F_{A})\psi  \nonumber
\end{eqnarray}

\noindent where $\rho (F_{A})$ $=$ $\sum_{a<b}\Gamma (e_{a})\Gamma
(e_{b})F_{A}(e_{a},e_{b}).$

To deal with the $T$-derivative term in (\ref{c2}), we consider the
canonical spin$^{c}$ structure with $W$ $=$ $\Lambda ^{0,\ast },$ the bundle
of all $(0,q)$ forms. In particular, we take $\psi $ $=$ $\bar{\partial}%
_{b}u $ $=$ $u_{,\bar{\beta}}\theta ^{\bar{\beta}}$ (summation convention
throughout the remaining part)$,$ a $(0,1)$ form with components being
derivatives of a complex function $u.$ Then we have%
\begin{equation*}
D_{\xi }\psi =(\bar{\partial}_{b}+\bar{\partial}_{b}^{\ast })\circ \bar{%
\partial}_{b}u = \bar{\partial}_{b}^{\ast }\circ \bar{\partial}_{b}u=\square _{b}u
\end{equation*}

\noindent where $\square _{b}$ $:=$ $\bar{\partial}_{b}\circ \bar{\partial}%
_{b}^{\ast }$ $+$ $\bar{\partial}_{b}^{\ast }\circ \bar{\partial}_{b}$ is
Kohn's Laplacian. Note that $\bar{\partial}_{b}^{\ast }u$ $=$ $0.$ To solve
the contact Dirac equation $D_{\xi }\psi $ $=$ $0$ for $\psi $ with suitably
asymptotic behaviour at the infinity is reduced to solving $\square _{b}u$ $=$
$0$ for $u$ with corresponding behaviour at the infinity. Associated to the $T
$-derivative term is a term involving so-called CR Paneitz operator $P$
after integration. Observe that for $n$ $\geq $ $2$ (i.e. in real dimension $2n+1$ $\geq $
$5)$, $P$ is non-negative (for closed $M$ and open $M$ with suitably decaying 
test functions) (see \cite{CC}, \cite{GL}). In dimension 3, the embeddability of
the underlying CR structure is necessary for its non-negativity, see Remark 
\ref{r:necppos}. 

On the other hand, in dimension 3, the trace curvature term $\rho
(F_{A})\psi $ is absorbed in the scalar curvature term. So by further
assuming $R$ $\geq $ $0,$ we can have the non-negativity of the p-mass
(which we pick up from the boundary terms).
In our main theorem we give some general conditions which ensure
the non-negativity of the p-mass for blow-ups of compact manifolds, 
characterizing also the zero case as
(CR equivalent to) the standard three dimensional CR sphere.

\begin{thm}\label{t:pm} Let M be a smooth, strictly pseudoconvex three
dimensional compact CR manifold.  Suppose $\mathcal{Y}(J) > 0$,
and that the CR Paneitz operator is non-negative. Let $p \in M$ and let
$\theta$ be a blow-up of contact form as in \eqref{eq:bbuu}.
Then
\begin{description}
  \item[(a)] $m(J,\th) \geq 0$;
  \item[(b)] if $m(J,\th) = 0$, M is CR equivalent (or, together with $\hat{\theta}$, isomorphic as pseudohermitian
  manifold) to $S^3$, endowed with its standard CR structure (and its standard contact form).
\end{description}
\end{thm}

\noindent We prove this theorem in Section \ref{s:pfthm}. The assumptions we give
here are conformally invariant, and are needed to ensure the positivity of the
right-hand side in \eqref{eq:massintintr}. By the result in \cite{CHCY1}, the
conditions on $\mathcal{Y}(J)$ and $P$ imply the embeddability of $M$: we use this
property to find a solution of $\Box_b \beta = 0$ with the above asymptotics (and hence to make
the first term in the right-hand side of \eqref{eq:massintintr} vanish): we first
find an approximate solution through the expansion of $\Box_b \ov{z}$ at infinity,
and then through the analysis of the Szeg\"{o} projection of this quantity, see Subsection
\ref{ss:exbeta}. To obtain the full solvability of $\Box_b \beta = 0$ we then
employ a mapping theorem in weighted spaces from \cite{HY}.
The positivity of the CR Paneitz operator is used instead to control the last term in the
right-hand side of \eqref{eq:massintintr}, showing that it is the sum of a non-negative
term and a (negative) multiple of $m(J,\th)$ which can be reabsorbed
into the left-hand side, see Subsection \ref{ss:P>0}. More comments on the relation
between the embeddability and the non-negativity of $P$ for manifolds of
positive Tanaka-Webster class are given in Remark \ref{r:necppos}. As a matter of fact, non-negativity of the CR Paneitz operator is preserved under embedded analytic deformations.

\begin{corollary}
Let M be a smooth, strictly pseudoconvex three dimensional compact CR manifold. Suppose M is an embedded, small enough, analytic deformation of the standard CR three sphere. Let $p \in M$ and let
$\theta$ be a blow-up of contact form as in \eqref{eq:bbuu}.
Then the same conclusions of Theorem \ref{t:pm} hold. 
\end{corollary}

\

\noindent In Section \ref{s:ex} we construct some examples of structures
using the deformation formulas in Subsection \ref{ss:variations}. First, using second
variation formulas, in Subsection \ref{ss:pneg} we consider perturbations of the spherical
structure for which $P$ fails to be non-negative, see Proposition
\ref{p:P<0} (and also \cite{CHCY1}). Then,
in Subsection  \ref{ss:12mass} we derive the first and second variations
of the mass near the standard sphere. We also construct examples of manifolds
with positive Tanaka-Webster class and negative
mass (when the blow-up is done at suitable points), see Proposition \ref{p:example}.
This is in striking contrast with respect to the Riemannian case, where all
perturbations of the sphere give rise to blown-up manifolds with positive mass
(except for metrics conformally equivalent to the spherical one).  We also describe 
an example of CR structure on $S^2 \times S^1$ with non-negative Paneitz operator and 
non-vanishing torsion, 
obtained as quotient of $\H^1 \setminus \{0\}$.

\

\noindent Our next main goal is to apply Theorem \ref{t:pm} to the study of
the CR Yamabe problem, namely finding conformal changes of contact form
in order to obtain constant Tanaka-Webster curvature. As for the classical Yamabe
problem, the cases $\mathcal{Y}(J) \leq 0$ are more directly treatable (see \cite{CHHAM}), while
the case $\mathcal{Y}(J) > 0$ is the most difficult one. Calling $\mathcal{Y}_0$
the quotient for the standard CR three sphere, by a result in \cite{JL0} one
always has
\begin{equation}\label{eq:ineqJ}
    \mathcal{Y}(J) \leq \mathcal{Y}_0,
\end{equation}
and if strict inequality holds, then the problem is solvable. The strict inequality
is needed to ensure compactness of the minimizing sequences in \eqref{eq:Y(J)}.
This condition was verified in \cite{JL} for (real) dimension greater or
equal to five, and for non-spherical structures, in the spirit
of \cite{aub} through some expansions involving the local geometry.

The positivity of the mass is instead a more global property, and it enters when 
$G_{p}$ has the following expansion near $p:$%
\begin{equation}
G_{p}=c_{n}\rho ^{-2n}+A+O(\rho ),  \label{c3}
\end{equation}
where $\rho $ is the Heisenberg distance in CR normal
coordinates. It turns out  that the term $A$ is a multiple of the mass defined for the
blow-up $M.$ We observe that (\ref{c3}) holds for $n$ $=$ $1$ (dimension 3
case) and for $N$ being spherical of all dimensions.

For such manifolds of dimension greater or equal to 5 (some extra technical condition in dimension
5) with positive CR Yamabe or Tanaka-Webster class, one can prove a
positive mass theorem for $A$ (and hence find solutions of the CR Yamabe
problem with minimal energy) through another approach (\cite{CCY}).

Our next result gives the strict Webster-Sobolev inequality in the three
dimensional case (the only one left), if $M$ is not CR equivalent to the standard
CR three-sphere, under the same assumptions as in the previous theorem.

\begin{thm}\label{t:y}
Suppose we are under the assumptions of Theorem \ref{t:pm}. Then either $M$
is the standard CR three-sphere or if M is not CR equivalent to the standard CR
three-sphere one has
$\mathcal{Y}(J) < \mathcal{Y}_0$. In both cases,
the Tanaka-Webster quotient admits a smooth minimizer.
\end{thm}

\noindent The CR Yamabe problem for the case of three-dimensional
CR manifolds and for  spherical CR manifolds was solved in
\cite{GA} and \cite{GJ} respectively (we also refer to \cite{CHHAM} and \cite{CCY}). While the proof in
these papers relies on topological arguments and may not provide 
energy extremals, in the spirit of \cite{babr},
our argument is based on direct minimization and gives an extra variational
characterization on the solutions (see more comments in Remark \ref{r:nonmin}).
To prove strict inequality we follow Schoen's argument in \cite{Sc}, finding test
functions which resemble a CR bubble at a small scale, and the Green's function
$G_p$ at a larger one. The construction is performed in Section \ref{s:y}.
More in general, the analysis of the Yamabe problem in the CR case has been
so far less precise than the Riemannian case: for example a basic difficulty is
the lack of a moving plane method, which is useful in general to derive a
priori estimates and to classify entire solutions.

\

\noindent In the appendix we collect some useful tools: first the variations of
some geometric quantities under a deformation of the CR structure, and
then the CR normal coordinates introduced in \cite{JL}. The latter are needed here
to derive asymptotic estimates on the Green's function near its singularity and
on the geometric quantities on asymptotically flat manifolds.

\subsection{Notation and preliminaries}\label{ss:notprel}

We collect here some useful material: throughout this paper, we will mostly
use notations taken from \cite{L1}. Consider a three
dimensional CR manifold endowed with a contact structure $\xi$ and a CR structure
$J : \xi \to \xi$ such that
$J^2 = - 1$. We assume that there exists a global choice of contact
form $\th$ which annihilates $\xi$ and for which $\th \wedge d \th$
is always nonzero. We define the {\em Reeb vector field} as the unique
vector field $T$ for which
$$
  \th(T) \equiv 1; \qquad \qquad T \lrcorner \; d \th = 0.
$$
Given $J$ as above, we have a local choice of a vector field $Z_1$
such that
\begin{equation}\label{eq:J0}
    J Z_1 = i Z_1; \qquad  J Z_{\ou} = - i Z_{\ou} \qquad \quad
  \hbox{ where } \quad Z_{\ou} = \overline{(Z_1)}.
\end{equation}
We also define $(\th, \th^1, \th^{\ou})$ as the dual triple to $(T,
Z_1, Z_{\ou})$, so that
$$
  d \th = i h_{1 \ou} \th^1 \wedge \th^{\ou} \qquad \hbox{ for some }
  h_{1 \ou} > 0   \quad (\hbox{possibly replacing } \theta \hbox{ by } -\theta). 
$$
In the following we will always assume that $h_{1 \ou} \equiv 1$.

The connection 1-form $\o^1_1$ and the torsion are uniquely
determined by the equations
$$
   \left\{
     \begin{array}{ll}
       d \th^1 = \th^1 \wedge \o^1_1 + A^1_{\ou} \th \wedge \th^{\ou};
         & \\
       \o^1_1 + \o^{\ou}_{\ou} = 0. &
     \end{array}
   \right.
$$
The {\em Tanaka-Webster curvature} is then defined by the formula
$$
  d \o^1_1 = R \th^1 \wedge \th^{\ou} \; (\hbox{mod } \th).
$$
In the (3-dimensional) Heisenberg group $\mathbb{H}^1$ standard choices for the dual
forms are
$$
  \left\{
    \begin{array}{ll}
      \stackrel{\circ}{\theta} = dt +  i z d \ov{z} - i \ov{z} d z; &  \\
      \stackrel{\circ}{\th^1} = \sqrt{2} d z; &  \\
      \stackrel{\circ}{\th^{\ou}} = \sqrt{2} d \ov{z}, &
    \end{array}
  \right.
$$
while for the vector fields we can take
\begin{equation}\label{eq:standvf}
    \left\{
    \begin{array}{ll}
      \stackrel{\circ}{T} = \frac{\pa}{\pa t}; & \\
      \stackrel{\circ}{Z}_1 = \frac{1}{\sqrt{2}} \left( \frac{\pa}{\pa z} + i \ov{z}
       \frac{\pa}{\pa t} \right); & \\
      \stackrel{\circ}{Z}_{\ou} = \frac{1}{\sqrt{2}} \left( \frac{\pa}{\pa \ov{z}} - i z
      \frac{\pa}{\pa t} \right). &
    \end{array}
  \right.
\end{equation}

The standard contact structure $\xi_0$ on $\mathbb{H}^1$ is spanned by real and imaginary
parts of $\stackrel{\circ}{Z}_1$ at each point. We define the standard CR structure
$J_{0}:\xi_{0}\to \xi_{0}$ by $J_{0}\stackrel{\circ}{Z}_1 = i\stackrel{\circ}{Z}_1$ on the complexification
of $\xi_0$. $(\mathbb{H}^1 ,J_{0},\stackrel{\circ}{\theta})$ is a pseudohermitian 3-manifold with
$\o^1_1 = A^1_{\ou} = R = 0$.

We recall the rules for the covariant differentiation with respect
to $\o^1_1$. For a complex function $f$ we set
$$
  f_1 = f,_{1} \equiv Z_1 f; \qquad \quad f,_{1 \ou} = Z_{\ou} Z_1 f
  - \o^1_1(Z_{\ou}) Z_1 f; \qquad \quad f,_0 = T f, \dots.
$$
We also define the operators
$$
  \D_b f = f,_1^{\; \; \; 1} + f,_{\ou}^{\; \; \; \ou} = f_{, 1\ou} + f_{, \ou 1};
  \qquad \qquad \Box_b f = - \D_b f + i T f,
$$
where we have used $h^{1 \ou} = h_{1 \ou} = 1$ to raise or lower the indices.
We also recall the commutation relations
\begin{equation}\label{eq:comm}
     \left\{
     \begin{array}{ll}
       c,_{1 \ou} - c,_{\ou 1} = i c,_0 + k c R;  & \\
       c,_{01} - c,_{10} = c,_{\ou} A_{11} - k c A_{11},_{\ou}; &  \\
       c,_{0\ou} - c,_{\ou 0} = c,_1 A_{\ou \ou} + k
       c A_{\ou \ou},_1,  &
     \end{array}
   \right.
\end{equation}
 (see Lemma 2.3 in \cite{L2}; here we are considering the $n$ $%
=$ $1$ case) where $c$ is a tensor with $1$ or $\bar{1}$ as subindices, $k$
is the number of $1$-subindices of $c$ minus the number of $\ov{1}$-subindices
of $c$ and where, we recall, we are assuming
that $h_{1\overline{1}}=1$ (so $A_{\bar{1}\bar{1}}$ $=$ $A_{\bar{1}}^{1}$
and $A_{11}$ is the complex conjugate of $A_{\bar{1}\bar{1}}).$

The CR Paneitz operator $P$ is defined by%
\begin{equation}
P\varphi :=4(\varphi {_{\bar{1}}}^{\bar{1}}{_{1}}+iA_{11}\varphi
^{1})^{1}.  \label{Pan}
\end{equation}

\noindent Let $\tilde{P}_3 \varphi $ $:=$ $\varphi {_{\bar{1}}}^{\bar{1}}{_{1}}%
+iA_{11}\varphi ^{1}.$ The CR pluriharmonic functions are characterized by
$\tilde{P}_{3}\varphi $ $=$ $0$ (\cite{L2}) while $P$ is identified with the
compatibility operator for solving a certain degenerate Laplace equation
(see \cite{GL}). The operator $P$ is also a CR analogue of
the Paneitz operator in conformal geometry (see \cite{Hi} for the relation
to a CR analogue of the $Q$-curvature and the $\log $-term coefficient in
the Szeg\"{o} kernel expansion). On a compact pseudohermitian 3-dimensional
manifold $(M,$ $J,$ $\theta ),$ we call $P$ {\em non-negative} if
\begin{equation}
\int_{M}\varphi P\varphi \theta \wedge d\theta \geq 0  \label{NG}
\end{equation}

\noindent for all real ($C^{\infty })$ smooth functions $\varphi.$
When we want to
emphasize the dependence of $\theta ,$ we write $P_{\theta }$ instead of $P.$
$P$ is pseudohermitian covariant in the sense that
\begin{equation}
P_{\hat{\theta}}\varphi =e^{4f} P_{\theta }\varphi  \label{PT}
\end{equation}
for the contact form change $\theta = e^{2f} \hat{\theta}$ (\cite{Hi}). The CR Paneitz
operator enters in the assumptions of the following embeddability
theorem.

\begin{thm} (\cite{CHCY1}) Let $M$ be a compact three-dimensional CR manifold.
If $P$ is non-negative and $R > 0$, then every eigenvalue $\l \neq 0$ of $\Box_b$
is greater or equal to $\min_M R$. In particular the range of $\Box_b$ is closed.
If $P$ is non-negative and $\mathcal{Y}(J) > 0$, then $M$ can be embedded into
$\C^N$, for some integer $N$.
\end{thm}

\begin{rem}\label{r:necppos} For a large class of pseudohermitian structures close to the
spherical one (see Theorem 1.9 in \cite{CHCY1}) the non-negativity of $P$ implies the
embeddability of the structure. A partial converse is shown in \cite{CHCY2}, 
where non-negativity of the CR Paneitz operator is shown to be preserved 
under embedded analytic deformations. The non-negativity of $P$ is also used in
\cite{CM} to obtain some hessian bounds
on complex functions in terms of integral bounds on the sublaplacian.

Examples with P non-negative include torsion-free (compact pseudohermitian) 3-manifolds (see Theorem 5.3 in
\cite{CCC}) and non torsion-free 3-manifolds like $S^{2} \times
S^{1}$, see Subsection \ref{ss:s2s1}
\end{rem}

\

\noindent In the systems of coordinates we will use below (both near zero
or near infinity), we will set
\begin{equation}\label{eq:rho}
    \rho^4 =  |z|^4 + t^2.
\end{equation}
In the Heisenberg group $\H^1$ the function $\rho$ has homogeneity $1$ with
respect to the natural dilation $(z,t) \mapsto (\l z, \l^2 t)$, $\l > 0$.

For a compact CR manifold $M$ and for $q > 1$ we define the Folland-Stein space
$\mathfrak{S}^{1,q}$ to be the completion of the (complex-valued) $C^\infty$
functions on $M$ with respect to the norm
$$
  \left( \int_M (u_{, 1} \ov{u}_{, \ou} + u_{, \ou} \ov{u}_{, 1})^{\frac q2}
  \th \wedge d \th \right)^{\frac 1q} + \left( \int_M |u|^q \th \wedge d \th
  \right)^{\frac 1q}.
$$
We define the spaces $\mathfrak{S}^{k,q}$ ($k \in \N$) in a similar way,
taking all possible combinations of $1$ and $\ou$ derivatives up to
order $k$.

For $k \in \Z$ we  denote by $\tilde{O}(\rho^k)$ a function $f(z,\ov{z},t)$ for
which $|f| \leq C
\rho^k$ for some $C > 0$; we  use instead the symbol $\tilde{O}'(\rho^k)$ for a
function $f(z,\ov{z},t)$ such that
$$
  |f| \leq C \rho^k, \qquad |\pa_z f| \leq C \rho^{k-1} \left| \pa_z \rho \right|,
  \qquad |\pa_{\ov{z}} f| \leq C \rho^{k-1} \left| \pa_{\ov{z}} \rho \right|,
  \qquad |\pa_t f| \leq C \rho^{k-2} \left| \pa_t \rho \right|.
$$
One can define similarly the symbols $\tilde{O}''(\rho^k)$, $\tilde{O}'''(\rho^k)$, etc.
We will use the symbol $O(\rho^k)$ for a function which is of the form
$\tilde{O}^{(j)}(\rho^k)$ for every integer $j$, or for $j$ large enough for our purposes.

Large positive constants are always denoted by $C$, and the value of
$C$ is allowed to vary from formula to formula and also within the
same line. When we want to stress the dependence of the constants on
some parameter (or parameters), we add subscripts to $C$, as $C_\d$,
etc.. Also constants with this kind of subscripts are allowed to
vary.

\

\noindent {\bf Acknowledgements} J.-H. C. (P.Y., resp.) are grateful to SISSA and
Princeton University (Academia Sinica in Taiwan, resp.) for the kind hospitality.
He also would like to thank Hung-Lin Chiu and Jack Lee for many discussions on
this topic at an early stage of this research project. J.-H. C. would like to thank
Chin-Yu Hsiao for teaching him the distributional argument used in the proof of
\eqref{eq:beta-1ou}.
A.M. has been supported by the project FIRB-IDEAS "Analysis and Beyond" and
by the Giorgio and Elena Petronio Fellowship while visiting IAS in Princeton in the
Fall Semester 2008-2009. He also would like to thank Academia Sinica in Taiwan
and Princeton University for the kind hospitality.

\section{Asymptotically flat pseudohermitian manifolds and definition of the  p-mass}\label{s:afman}

One of our main goals is to study compact manifolds with positive Tanaka-Webster invariant.
For these, the conformal sublaplacian is positive definite, and therefore by classical
arguments it admits
a Green's function. By Proposition \ref{p:exG} (we refer to Subsection \ref{ss:green}
for details), the Green's function $G_p$ with pole at $p \in M$, namely the solution of
$$
- 4 \D_b G_p +  R G_p = 16 \d_p,
$$
in CR normal coordinates $(z,t)$ centered at $p$ (see Subsection \ref{ss:CRcoord}) has the form
$$
  G_p = \frac{1}{2 \pi \rho^2} + w,
$$
where $w$ is a function of class $C^{1,\a}$ for any $\a \in (0,1)$ and where, following
our notation, $\rho^4 = |z|^4 + t^2$.

In fact, a more refined regularity property on $w$ can be proved, see Appendix 1 in Section 10 of
\cite{HY}: this fact will be needed to apply the main theorem there.

\subsection{Blow-up through the Green's function}\label{ss:bu}

Let $(M,J)$ be a CR three-manifold, and let $p \in M$.
By Proposition \ref{p:CRcoord} (which relies on results by D.Jerison and
J.Lee from \cite{JL}), we can find a contact form $\hat{\th}$ near $p$ and
local coordinates $(z, t)$ such that
$$
  \left\{
    \begin{array}{ll}
      \hat{\th} = \left( 1 + O(\rho^4) \right) \theta_0
     + O(\rho^5) d z + O(\rho^5) d \ov{z}; &  \\
      \hat{\th}^1 = \sqrt{2} \left( 1 + O(\rho^4) \right) d z
     + O(\rho^4) d \ov{z} + O(\rho^3) \theta_0,  &
    \end{array}
  \right.
$$
where we have set
$$
      \theta_0 = d t +  i z d \ov{z} - i \ov{z} d z.
$$
If $\hat{\th}$ is as above and if $\mathcal{Y}(J) > 0$ we consider the following
pseudohermitian manifold
\begin{equation}\label{eq:bubu}
    \left( M \setminus \{p\}, J, \th \equiv G_p^2 \hat{\th} \right),
\end{equation}
where we made a conformal change as in Subsection \ref{sss:conf}, taking $f = \log G_p$.
Then by \eqref{eq:th1f} we also have
\begin{equation}\label{eq:mucca}
  \th^1 = G_p \left( \hat{\th}^1 + 2 i \left(
  \log G_p \right)_{,\ou} \hat{\th} \right).
\end{equation}
By the expression of $G_p$ (see Proposition \ref{p:exG}), setting $A = w(0)$ we have that
\begin{equation}\label{eq:GpA}
      G_p = \frac{1}{2 \pi \rho^2} + A + O(\rho) \qquad \quad \hbox{ near } p.
\end{equation}
Hence, from Proposition \ref{p:CRcoord} we have that
\begin{eqnarray}\label{eq:th} \nonumber
  \th & = & G_p^2 \hat{\th} = \left( \frac{1}{2\pi} \rho^{-2} + A + O(\rho) \right)^2
  \left[ \left( 1 + O(\rho^4) \right) \theta_0 + O(\rho^5) dz
  + O(\rho^5) d\ov{z} \right] \\
   & = & \left( \frac{1}{(2\pi)^2} \rho^{-4} + 2 \frac{1}{2\pi}
 A \rho^{-2} + O(\rho^{-1}) \right)
  \left[ \left( 1 + O(\rho^4) \right) \theta_0 + O(\rho^5) dz
  + O(\rho^5) d\ov{z} \right] \\
   & = & \left( \frac{1}{(2\pi)^2} \rho^{-4} + 2 \frac{1}{2\pi}
A \rho^{-2} + O(\rho^{-1}) \right) \theta_0
  + O(\rho) dz + O(\rho) d\ov{z}. \nonumber
\end{eqnarray}
Recall also the second equation in \eqref{eq:mucca}, which we
rewrite as
$$
  \th^1 = G_p \left[ \hat{\th}^1 + 2 i \frac{(G_p)_{,\ou}}{G_p}
  \hat{\theta} \right] = G_p  \hat{\th}^1 + 2 i (G_p)_{,\ou} \hat{\theta}.
$$
Using the fact that
$$
  (\rho^{-2})_{,\ou} = \frac{1}{\sqrt{2}} \frac{i z v}{\rho^6} + O(\rho),
$$
where $v = t + i |z|^2$, we get
\begin{equation}\label{eq:Gpou}
    (G_p)_{, \ou} = \frac{1}{2\pi} \frac{1}{\sqrt{2}} \frac{i z v}{\rho^6} + O(1),
\end{equation}
from which we deduce
\begin{eqnarray}\label{eq:th1}\nonumber
  \th^1 & = & G_p  \hat{\th}^1 + 2 i (G_p)_{,\ou} \hat{\theta} \\ \nonumber
   & = & \left( \frac{1}{2\pi} \rho^{-2} + A + O(\rho) \right) \left(\sqrt{2}
   \left( 1 + O(\rho^4) \right) dz + O(\rho^4) d\ov{z} + O(\rho^3) \theta_0 \right)
   \\ & + &  2 i \left( \frac{1}{2\pi} \frac{1}{\sqrt{2}} \frac{i z v}{\rho^6}
   + O(1) \right) \left( \left( 1 + O(\rho^4) \right) \theta_0 + O(\rho^5) dz
  + O(\rho^5) d \ov{z} \right)   \\ & = & \left( \frac{1}{2\pi} \rho^{-2} + A +
  O(\rho) \right) \sqrt{2} dz + O(\rho^2) d\ov{z} + \left( - \frac{1}{\pi}
  \frac{1}{\sqrt{2}} \frac{z v}{\rho^6} + O(\rho) \right) \theta_0. \nonumber
\end{eqnarray}
Using \eqref{eq:varoconf}, the new connection form becomes
\begin{equation}\label{eq:o11GG}
    \o^1_1 = \hat{\o}^1_1 + 3 \left( \log G_p \right)_{,1} \hat{\th}^1 -
  3 \left( \log G_p \right)_{,\ou} \hat{\th}^{\ou}  + i
  \left( \D_b \left( \log G_p \right) + 8 \left( \log G_p \right)_{,1}
  \left( \log G_p \right)_{,\ou} \right) \hat{\th}.
\end{equation}
From \eqref{eq:Gpou} we find
\begin{eqnarray*}
  \left( \log G_p \right)_{,\ou} & = & \frac{(G_p)_{, \ou}}{G_p} =
  \frac{\frac{1}{2\pi} \frac{1}{\sqrt{2}} \frac{izv}{\rho^6} + O(1)}{
  \frac{1}{2\pi} \rho^{-2} + A + O(\rho)} = \frac{
  \frac{1}{\sqrt{2}} \frac{izv}{\rho^4} + O(\rho^2)}{
  1 + 2 \pi A \rho^2 + O(\rho^3)} \\
   & = & \frac{1}{\sqrt{2}} \frac{i z v}{\rho^4} + O(\rho^2),
\end{eqnarray*}
which implies
\begin{eqnarray}\label{eq:333} \nonumber
  & & 3 \left( \log G_p \right)_{,1} \hat{\th}^1 -
  3 \left( \log G_p \right)_{,\ou} \hat{\th}^{\ou} \\ & = &
  \left( - \frac{3i \ov{z} \ov{v}}{\sqrt{2} \rho^4} + 3 \sqrt{2} \pi A
  \frac{i \ov{z} \ov{v}}{\rho^2} + O(\rho^2) \right) \left( \sqrt{2} dz
  + O(\rho^4) dz + O(\rho^4) d\ov{z} + O(\rho^3) \th_0 \right) \\
  & + & \left( - \frac{3i z v}{\sqrt{2} \rho^4} + 3 \sqrt{2} \pi A
  \frac{i z v}{\rho^2} + O(\rho^2) \right) \left( \sqrt{2} d \ov{z}
  + O(\rho^4) dz + O(\rho^4) d\ov{z} + O(\rho^3) \th_0 \right). \nonumber
\end{eqnarray}
On the other hand, we also have that
\begin{eqnarray*}
  \D_b \left( \log G_p \right) + 8 \left( \log G_p \right)_{,1}
  \left( \log G_p \right)_{,\ou} &=& \left( \frac{(G_p)_{, 1}}{G_p} \right)_{, \ou}
  + \left( \frac{(G_p)_{, 1}}{G_p} \right)_{, \ou} \left( \frac{(G_p)_{, \ou}}{G_p}
  \right)_{, 1} + 8 \frac{(G_p)_{, 1}}{G_p} \frac{(G_p)_{, \ou}}{G_p}  \\
  & = & \frac{\D_b G_p}{G_p} + 6 \frac{(G_p)_{, \ou} (G_p)_{, 1}}{G_p^2}.
\end{eqnarray*}
%Moreover, using the above formulas, the $\hat{\th}$ coefficient of $\o^1_1$ satisfies
%\begin{equation}\label{eq:ccppuu}
%    \D_b \left( \log G_p \right) + 8 \left( \log G_p \right)_{,1}
%  \left( \log G_p \right)_{,\ou} = 3 i \frac{|z|^2}{\rho^4} + O(\rho^{-3}).
%\end{equation}
Letting $\stackrel{\circ}{\D_b}$ denote the standard sublaplacian in the
Heisenberg group, from the choice of CR normal coordinates we find
\begin{eqnarray*}
  \D_b & = & \stackrel{\circ}{\D_b} + O(\rho^4) \stackrel{\circ}{Z}_1 \stackrel{\circ}{Z}_1
  + O(\rho^4) \stackrel{\circ}{Z}_1 \stackrel{\circ}{Z}_{\ou}  +
  O(\rho^4) \stackrel{\circ}{Z}_{\ou} \stackrel{\circ}{Z}_1  +
  O(\rho^4) \stackrel{\circ}{Z}_{\ou} \stackrel{\circ}{Z}_{\ou} +
  O(\rho^5) \stackrel{\circ}{Z}_1 \pa_t + O(\rho^5) \stackrel{\circ}{Z}_{\ou} \pa_t
  \\ & + & O(\rho^5) \pa_t \stackrel{\circ}{Z}_1 + O(\rho^5) \pa_t \stackrel{\circ}{Z}_{\ou}
  + O(\rho^3) \stackrel{\circ}{Z}_1  +  O(\rho^3) \stackrel{\circ}{Z}_{\ou}
  + O(\rho^4) \pa_t.
\end{eqnarray*}
Using this expansion and the fact that $\stackrel{\circ}{\D_b} \rho^{-2} = 0$ in
$\H^1$ we obtain
$$
  \frac{\D_b G_p}{G_p} = \frac{\D_b \left( \frac{1}{2\pi \rho^2} + A + O(\rho)
  \right)}{\frac{1}{2\pi \rho^2} + A + O(\rho)} = O(\rho) \qquad \quad
  (\rho \hbox{ small}),
$$
and
\begin{eqnarray*}
  6 \frac{(G_p)_{, \ou} (G_p)_{, 1}}{G_p^2} & = & 6 \frac{\left( \frac{1}{2\pi}
  \frac{1}{\sqrt{2}} - \frac{i \ov{z} \ov{v}}{\rho^6} + O(1) \right) \left(
  \frac{1}{2\pi} \frac{1}{\sqrt{2}} + \frac{i z v}{\rho^6} + O(1) \right)}{\left(
  \frac{1}{2\pi \rho^2} + A + O(\rho)\right)^2} \\
  & = & 6 \frac{\left( \frac{1}{\sqrt{2}} \left( - \frac{i \ov{z} \ov{v}}{\rho^4}
  + O(\rho^2) \right) \right) \left( \frac{i z v}{\rho^4}
  + O(\rho^2) \right)}{\left( 1 + 2\pi A \rho^2 + O(\rho^3) \right)^2}
  = 3 \frac{|z|^2}{\rho^4} - 3 \frac{|z|^2}{\rho^4} 4 \pi A \rho^2 + O(\rho).
\end{eqnarray*}
Therefore we get
\begin{equation}\label{eq:2222}
    \D_b \left( \log G_p \right) + 8 \left( \log G_p \right)_{,1}
  \left( \log G_p \right)_{,\ou} = 3 \frac{|z|^2}{\rho^4} - 12
  \frac{|z|^2}{\rho^2} \pi A + O(\rho).
\end{equation}
From \eqref{eq:333} and \eqref{eq:2222}, taking Proposition \ref{p:CRcoord} into
account then one finds
\begin{eqnarray}\label{eq:starstar} \nonumber
  \o_1^1 & = & - 3 i \frac{\ov{z} \ov{v}}{\rho^4} dz - 3 i \frac{z v}{\rho^4} d \ov{z}
  + 3 i \frac{|z|^2}{\rho^4} \th_0
  \\ & + & \left( 3 \sqrt{2} \pi A i \frac{\ov{z} \ov{v}}{\rho^2}
   + O(\rho^2) \right) \sqrt{2} dz + \left( 3 \sqrt{2} \pi A i
   \frac{z v}{\rho^2} + O(\rho^2) \right) \sqrt{2} d\ov{z} \\ & + &
   \left( - 12 \pi A i \frac{|z|^2}{\rho^2} + O(\rho) \right) \th_0. \nonumber
\end{eqnarray}

\subsection{CR inversion}\label{ss:inv}

Next we want to express $\th$ and $\th^1$ in {\em inverted CR normal
coordinates} $(z_*,t_*)$. If $(z,t)$ are CR normal coordinates in
a neighborhood $\mathcal{U}$ of $p$, we define the inverted CR
normal coordinates as
\begin{equation}\label{eq:zsttst}
    z_* = \frac{z}{v}; \qquad  t_* = - \frac{t}{|v|^2}; \quad \qquad
  \hbox{ on } \mathcal{U} \setminus \{p\},
\end{equation}
where, as before, $v = t + i |z|^2$.

We have that
\begin{eqnarray*}
  v = t + i |z|^2 = - t_* |v|^2 + i |z_* v|^2
  = |v|^2 (- t_* + i |z_*|^2) = |v|^2 (- \overline{v}_*),
\end{eqnarray*}
with $v_* = t_* + i |z_*|^2$. The following two identities hold
$$
  1 = - \overline{v} \overline{v}_* \qquad \quad \hbox{ or } \qquad \quad
  1 = - v v_*,
$$
which imply
$$
  z = - \frac{z_*}{v_*}; \qquad \quad t = - \frac{t_*}{|v_*|^2}.
$$
We can write that
$$
  dz = \frac{\pa z}{\pa z_*} dz_* + \frac{\pa z}{\pa \overline{z}_*}
  d \overline{z}_* + \frac{\pa z}{\pa t_*} dt_*; \qquad \quad d t =
  \frac{\pa t}{\pa z_*} dz_* + \frac{\pa t}{\pa \overline{z}_*}
  d \overline{z}_* + \frac{\pa t}{\pa t_*} dt_*,
$$
so using elementary computations one finds
$$
  \frac{\pa z}{\pa z_*} = - \frac{t_*}{v_*^2}; \qquad \quad
  \frac{\pa z}{\pa \overline{z}_*} = i \frac{z_*^2}{v_*^2}; \qquad \quad
  \frac{\pa z}{\pa t_*} = \frac{z_*}{v_*^2};
$$
$$
  \frac{\pa t}{\pa z_*} = 2 \frac{t_* |z_*|^2 z_*^{\ou}}{|v_*|^4}; \qquad
  \quad \frac{\pa t}{\pa \overline{z}_*} = 2 \frac{t_* z_*^1 |z_*|^2}{|v_*|^4};
  \qquad \quad \frac{\pa t}{\pa t_*} = \frac{t_*^2 - |z_*|^4}{|v_*|^4}.
$$
These relations imply
\begin{equation}\label{eq:rel1}
    \theta_0 = \frac{(\theta_0)_*}{|v_*|^2} = \frac{(\theta_0)_*}{(\rho_*)^4}; \qquad
  \quad \rho = \frac{1}{\rho_*}, \quad (\rho_*)^4 = (t_*)^2 + |z_*|^4,
\end{equation}
where $(\theta_0)_* = dt_* + i z_* d \overline{z}_* - i \overline{z}_*
dz_*$, and also
\begin{equation}\label{eq:dzstst}
    dz = \frac{- t_* dz_* + i z_*^2 d \overline{z}_* + z_* dt_*}{v_*^2} =
  \frac{z_* (\theta_0)_* - \overline{v}_* dz_*}{v_*^2}.
\end{equation}
Now from \eqref{eq:th} we have that
\begin{eqnarray*}
  \th & = & G_p^2 \hat{\th} = \left( \frac{1}{(2\pi)^2}
\rho^{-4} + 2 \frac{1}{2\pi} A \rho^{-2} +
  O(\rho^{-1}) \right) \theta_0 + O(\rho) dz + O(\rho) d\ov{z} \\
   & = & \left( \frac{1}{(2\pi)^2} \rho_*^{4} + 2 \frac{1}{2\pi} A \rho_*^{2} + O(\rho_*)
  \right) \frac{(\theta_0)_*}{(\rho_*)^4} + O(\rho_*^{-1})
  \frac{z_* (\theta_0)_* - \overline{v}_* dz_*}{(v_*)^2}
  + O(\rho_*^{-1}) \frac{\overline{z}_* (\theta_0)_* -
v_* d\overline{z}_*}{(\overline{v}_*)^2} \\
   & = & \left( \frac{1}{(2\pi)^2} + \frac{1}{\pi}
A \rho_*^{-2} + O(\rho_*^{-3}) \right) (\theta_0)_*
  + O(\rho_*^{-3}) dz_*  + O(\rho_*^{-3}) d\overline{z}_*.
\end{eqnarray*}
Similarly from \eqref{eq:th1} we deduce
\begin{eqnarray*}
  \th^1 & = & \left( \frac{1}{2\pi} \rho^{-2} + A +
O(\rho) \right) \sqrt{2} dz + O(\rho^2) d\ov{z}
  + \left( - \frac{1}{\pi} \frac{1}{\sqrt{2}}
\frac{z v}{\rho^6} + O(\rho) \right) \theta_0 \\
   & = & \left( \frac{1}{2\pi} \rho_*^2 + A + O(\rho_*^{-1}) \right) \sqrt{2} \left(
   \frac{z_* (\theta_0)_* - \overline{v}_* dz_*}{(v_*)^2} \right) + O(\rho_*^{-2})
   \left( \frac{\overline{z}_* (\theta_0)_* - v_*
d\overline{z}_*}{(\overline{v}_*)^2} \right) \\ & + &
   \left( - \frac{1}{\pi} \frac{1}{\sqrt{2}}
\frac{z_* (\rho_*)^6}{(v_*)^2} + O(\rho_*^{-1}) \right)
   \frac{(\theta_0)_*}{\rho_*^4} \\
   & = & \left( \sqrt{2} \frac{A z_*}{v_*^2} + O(\rho_*^{-5}) \right) (\theta_0)_*
  + O(\rho_*^{-4}) d\overline{z}_* + \left(
   - \frac{1}{\sqrt{2} \pi} \rho_*^2 \frac{\overline{v}_*}{(v_*)^2} - A \sqrt{2}
   \frac{\overline{v}_*}{v_*^2} + O(\rho_*^{-3}) \right) dz_*.
\end{eqnarray*}
It is irrelevant to multiply the form $\th^1$ by a complex unitary
factor, which we can choose to be $- \frac{v_*^3}{\rho_*^6} := e^{i \var}$. In
this way the new $\th^1$, which we will denote as $\th^1_n$ is given by
$$
  \th^1_n = \left( - \sqrt{2} \frac{A z_* v_*}{\rho_*^6} + O(\rho_*^{-5})
  \right) (\theta_0)_* + O(\rho_*^{-4}) d\overline{z}_* + \left(
   \frac{1}{\sqrt{2} \pi} + A \sqrt{2} \rho_*^{-2} + O(\rho_*^{-3}) \right) dz_*.
$$
Under this complex rotation the connection form changes according to
the law
$$
  \o^1_1 \mapsto (\o^1_1)_n = \o^1_1 + i d \var.
$$
The function $\var$ is defined by
$$
  i \var = \log \left( - \rho_*^2 \frac{\ov{v}_*}{v_*^2} \right)
 = \log \left( \rho^{-2} v^2 \ov{v}^{-1} \right) = 2 \log v - 2 \log \rho
 - \log \ov{v}.
$$
By elementary computations one finds
$$
  \stackrel{\circ}{Z}_1 \var = \frac{1}{\sqrt{2}} \frac{3 \ov{z}
  \ov{v}}{\rho^4}; \qquad \qquad \pa_t \var = - 3 \frac{|z|^2}{\rho^4}.
$$
Hence, using Proposition \ref{p:CRcoord} one has
\begin{eqnarray}\label{eq:323232} \nonumber
  d \var &=& \var_{, 1} \hat{\th}^1 + \var_{, \ou} \hat{\th}^{\ou} +
    \var_{, 0} \hat{\th} \\
   & = & \frac{1}{\sqrt{2}} \frac{3 \ov{z}
  \ov{v}}{\rho^4} \sqrt{2} dz +  \frac{1}{\sqrt{2}} \frac{3 z v}{\rho^4}
   \sqrt{2} d\ov{z} - 3 \frac{|z|^2}{\rho^4} \th_0 + O(\rho^3) dz
  + O(\rho^3) d \ov{z} + O(\rho^2) \th_0.
\end{eqnarray}
By \eqref{eq:starstar} and \eqref{eq:323232} we then obtain
\begin{equation}\label{eq:o11n}
    (\o^1_1)_n =  6 \pi A i \frac{\ov{z} \ov{v}}{\rho^2} dz + 6 \pi A i
    \frac{zv}{\rho^2} d \ov{z} - 12 \pi A i \frac{|z|^2}{\rho^2} \th_0 +
    O(\rho^2) dz + O(\rho^2) d \ov{z} + O(\rho) \th_0.
\end{equation}
From \eqref{eq:zsttst} and \eqref{eq:dzstst} one then gets
$$
  6 \pi A i \frac{\ov{z} \ov{v}}{\rho^2} dz + 6 \pi A i \frac{zv}{\rho^2}
  d \ov{z} = - 6 \pi A \frac{\ov{z_*}(|z_*|^2+it_*)}{\rho_*^6} dz_* - 6 \pi
  A \frac{z_* (-|z_*|^2+it_*)}{\rho_*^6}
  d \ov{z}_* + 12 i \pi A \frac{|z_*|^2}{\rho_*^6} (\th_0)_*,
$$
and the coefficient of $(\th_0)_*$ cancels with
$$
  - 12 \pi A i \frac{|z|^2}{\rho^2} \th_0 = - 12 \pi A i \frac{|z_*|^2}{\rho_*^2} (\th_0)_*.
$$
In conclusion, using also \eqref{eq:rel1}, in inverted CR coordinates we obtain
\begin{equation}\label{eq:o11tstar}
    (\o^1_1)_n = \left[ - 6 \pi A \frac{\ov{z_*}(|z_*|^2+it_*)}{\rho_*^6} +
    O(\rho_*^{-4}) \right] dz_* + \left[ - 6 \pi A \frac{z_* (-|z_*|^2+it_*)}{\rho_*^6}
    + O(\rho_*^{-4}) \right] d \ov{z}_* + O(\rho_*^{-5}) (\th_0)_*.
\end{equation}
We can finally rescale by the constant factor $e^f = 2 \pi$ to get the
new forms (where we omit the normalization symbol)
\begin{equation}\label{eq:nt}
    \th = \left( 1 + 4 \pi A \rho_*^{-2} + O(\rho_*^{-3})
  \right) (\theta_0)_* + O(\rho_*^{-3}) dz_*  + O(\rho_*^{-3}) d\overline{z}_*;
\end{equation}
\begin{equation}\label{eq:nt1}
    \th^1 = \left( - 2 \sqrt{2} \pi A \frac{z_* v_*}{\rho_*^6} + O(\rho_*^{-5})
  \right) (\theta_0)_* + O(\rho_*^{-4}) d\overline{z}_* + \left(
   1 + 2 \pi A \rho_*^{-2} + O(\rho_*^{-3}) \right) \sqrt{2} dz_*.
\end{equation}
Notice that by \eqref{eq:varoconf} the connection form will stay invariant.

Later on, in Subsection \ref{ss:expboxovz}, we will also expand the adapted frame
associated to this triple $(\th, \th^1, \th^{\ou})$.

\subsection{Definition of p-mass and an integral formula}\label{ss:asyint}

Having in mind the previous expansions, and in particular formulas
\eqref{eq:nt}, \eqref{eq:nt1}, we give the following definition of asymptotically
flat manifold. Recall (see Subsection 1.1) that $\mathbb{H}^1$ denotes the
(3-dimensional) Heisenberg group with standard pseudohermitian structure
$(J_{0},\theta_{0}=\stackrel{\circ}{\theta})$. Let $B_{\rho_0}$ denote the Heisenberg ball
of radius $\rho_0$.

\begin{df}\label{d:aflat} A  three dimensional pseudohermitian manifold $(N,J,\th )$
is said to be asymptotically flat pseudohermitian if $N = N_0 \cup N_\infty$, with
$N_0$ compact and $N_\infty$ diffeomorphic to $\mathbb{H}^1\setminus B_{\rho_0}$
in which $(J,\th )$ is close to $(J_{0},\theta_{0})$ in the sense that
$$
  \th = \left( 1 + 4 \pi A \rho^{-2} + O(\rho^{-3})
  \right) \theta_0 + O(\rho^{-3}) dz  + O(\rho^{-3}) d\overline{z}; 
$$
$$
 \theta^{1} = O(\rho^{-3})\theta_{0} + O(\rho^{-4})d{\bar z}
               + (1 + 2{\pi}A\rho^{-2} + O(\rho^{-3})){\sqrt 2}dz
$$
%
%$$
%    \th^1 = O(\rho^{-3}) \theta_0 + O(\rho^{-4}) d\overline{z} + \left(
%   1 + 2 \pi A \rho^{-2} + O(\rho^{-3}) \right) \sqrt{2} dz,
%$$
for some unitary coframe $\theta^1$ and some $A \in \R$ in some system of coordinates
(called asymptotic coordinates). We also require that $R \in L^1(N)$.
\end{df}

\begin{rem}\label{r:regtor}
From Definition \ref{d:aflat} and the structure equations in Subsection \ref{ss:notprel}, 
we can then have
$$
    \th^1 = \left( - 2 \sqrt{2} \pi A \frac{z v}{\rho^6} + O(\rho^{-4})
  \right) \theta_0 + O(\rho^{-4}) d\overline{z} + \left(
   1 + 2 \pi A \rho^{-2} + O(\rho^{-3}) \right) \sqrt{2} dz,
$$
and
$$
   \o^1_1 = \left[ - 6 \pi A \frac{\ov{z}(|z|^2+it)}{\rho^6} +
    O(\rho^{-4}) \right] dz + \left[ - 6 \pi A \frac{z (-|z|^2+it)}{\rho^6}
    + O(\rho^{-4}) \right] d \ov{z} + O(\rho^{-5}) \th_0.
$$
We look for an expansion of the vector field $Z_1$ of the type
$$
  Z_1 = (1 + a) \stackrel{\circ}{Z}_1 + b \stackrel{\circ}{Z}_{\ou}
  + c \stackrel{\circ}{T}
$$
(recall the notation in \eqref{eq:standvf}). From the three relations
$$
  0 = \th(Z_1) = (1+a) O(\rho^{-3}) + b O(\rho^{-3}) + c
  \left(  1 + 4 \pi A \rho^{-2} + O(\rho^{-3}) \right);
$$
$$
  1 = \th^1(Z_1) = (1+a) \left( 1 + 2 \pi A \rho^{-2} + O(\rho^{-3})
  \right) + b O(\rho^{-4}) + c \left( - 2 \sqrt{2} \pi A  \frac{z v}{\rho^6}
  + O(\rho^{-3}) \right);
$$
$$
   0 = \th^{\ou}(Z_1) = (1+a) O(\rho^{-4}) + b \left( 1 + 2 \pi A
  \rho^{-2} + O(\rho^{-3}) \right) + c \left( - 2 \sqrt{2} \pi A
  \frac{\overline{z} \overline{v}}{\rho^6} + O(\rho^{-4}) \right)
$$
we deduce that
$$
  c = O(\rho^{-3}); \qquad \quad a = - 2 \pi A \rho^{-2}
  + O(\rho^{-3}); \qquad \quad b = O(\rho^{-4}).
$$
Similarly we obtain 
$$
  T=(1-4\pi A\rho^{-2}+O(\rho^{-3})) \stackrel{\circ}{T} + \left(2{\sqrt 2}\pi A\frac{zv}{\rho^6} + O(\rho^{-4})\right) \stackrel{\circ}{Z}_{1} + \left(2{\sqrt 2}\pi A\frac{{\bar z}{\bar v}}{\rho^6} + O(\rho^{-4}) \right) \stackrel{\circ}{Z}_{\bar 1}.
$$
By our conventions about the quantities of the type $O(\rho^k)$, we also have that
$$
  A_{11} = O(\rho^{-4}); \qquad R = \underbrace{Z_1(\o_1^1(Z_{\ou}))
  - Z_{\ou}(\o^1_1(Z_1))}_{O(\rho^{-4})} + i \underbrace{\o^1_1(T)}_{O(\rho^{-5})}
  + 2 \underbrace{\o^1_1(Z_1) \o_1^1(Z_{\ou})}_{O(\rho^{-6})}.
$$
Note that if we use $1+O(\rho^{-2})$ as the ${\sqrt 2}dz$-coefficient of $\th^1$ in
Definition 2.1, then only the real part of $O(\rho^{-2})$ (namely, $2 \pi A \rho^{-2}$)
is determined. But then the $\theta_0$-coefficient of $\o^1_1$ is not $O(\rho^{-5})$
(only $O(\rho^{-4})$). On the other hand, we may relax the conditions in Definition 2.1
for other purposes in the future.
\end{rem}

\noindent After this definition we are ready to introduce the notion
of p-mass using a variational characterization, in the same spirit of \eqref{eq:varmass}.
Considering a one-parameter family of CR structures $J(s)$ and using the notation
of Subsection \ref{sss:varJ},  we have that
$$
  \dot{J} = 2 E = 2 E_{11} \th^1 \otimes Z_{\ou} + 2 E_{\ou \ou}
  \th^{\ou} \otimes Z_1.
$$
Denoting by $R(s)$ the corresponding Tanaka-Webster curvature and using
\eqref{eq:varoA}-\eqref{eq:varW}  we then find
\begin{eqnarray*}
   & & \frac{d}{ds}_{|s = 0} \int_N R(s) \, \th \wedge d \th =
  \int_N \dot{R} \, \th \wedge d \th \\
  & = & \int_N \left[ i \left( E_{11},_{\ou \ou} - E_{\ou \ou},_{11}
  \right) - \left( A_{11} E_{\ou \ou} + A_{\ou \ou} E_{11} \right) \right]
  \th \wedge d \th \\
  & = & - \int_N d \left( E_{11},_1 \th \wedge \th^1 \right) + \hbox{conj.}
  - \int_N \left( A_{11} E_{\ou \ou} + \hbox{conj.} \right) \th \wedge d \th \\
  & = & - \oint_{\infty} E_{11,\ou} \, \th \wedge \th^1 + \hbox{conj.}
  - \int_N \left( A_{11} E_{\ou \ou} + \hbox{conj.} \right) \th \wedge d \th \\
  & = & \oint_{\infty} i \dot{\o}^1_1 \wedge \th - \int_N \left( A_{11}
  E_{\ou \ou} + \hbox{conj.} \right) \th \wedge d \th.
\end{eqnarray*}
This formula leads us to the following definition.

\begin{df}\label{d:mass} Let $N$ be an asymptotically flat manifold.
Then we define the p-mass of $(N,J,\th)$ as
$$
  m(J, \th) := i \oint_{\infty} \o^1_1 \wedge \th :=
  \lim_{\L \to + \infty} i \oint_{S_\L} \o^1_1 \wedge \th,
$$
where we have set $S_\L = \left\{ \rho = \L \right\}$.
\end{df}

\begin{rem}\label{r:mass} By the above computations one finds
$$
  \frac{d}{ds}_{| s = 0} \left( - \int_N R(s) \, \th \wedge d \th + m(J(s),\th)
   \right) = \int_N \left( A_{11} E_{\ou \ou} + A_{\ou \ou} E_{11}
  \right) \th \wedge d \th,
$$
which is the counterpart of \eqref{eq:varmass}, as desired. Apparently the
definition of the p-mass depends on a choice of asymptotic coordinates and a
(co)frame $\theta^1$. Whether the p-mass is independent of the choice of
asymptotic coordinates or an {\em admissible} (co)frame is an
interesting problem. But we will not pursue it in this paper.
\end{rem}

\noindent First, we show that $m(J,\theta)$ can be expressed in terms of the
constant $A$ appearing in Definition \ref{d:aflat}.

\begin{lem}\label{l:m=A} If $m(J,\th)$ is as in Definition \ref{d:mass}
and if $A$ is as in \eqref{eq:GpA}, then
$$
  m(J,\th) = 48 \pi^2 A.
$$
\end{lem}

\begin{pf}
By \eqref{eq:o11tstar} and some elementary estimates one finds
\begin{eqnarray*}
  m(J,\th) & = & \oint_\infty i \o^1_1 \wedge \th = - 3 i 2 \pi A \oint_\infty
  \rho^{-6} \left[ (|z|^2 \ov{z} + i \ov{z} t) dz - (|z|^2 z - i z t)
  d \ov{z} \right] \wedge \th_0 \\
  & = & - 6 i \pi A \oint_S |z|^2 (\ov{z} dz - z d \ov{z}) \wedge \th_0
  - 6 i \pi A \oint_S (i \ov{z} t dz + i z t d \ov{z}) \wedge \th_0,
\end{eqnarray*}
where we have set
$$
  S = \{ \rho = 1 \}.
$$
Using the relations
\begin{equation}\label{eq:relpol}
    z = r e^{i \var}; \qquad dz = e^{i \var} (dr + i r d
  \var); \qquad \ov{z} dz - z d \ov{z} = 2 i r^2 d \var,
\end{equation}
we get
\begin{eqnarray*}
  m(J,\th) & = & - 6 i \pi A \oint_S \left[ 2 i r^4 d \var +
  2 i t r dr  \right] \wedge \th_0 =
 12 \pi A \oint_S \left[ r^4 d \var \wedge dt
  + t r (dr \wedge dt + 2 r^2 dr \wedge d \var) \right] \\ & = &
  12 \pi A \oint_S \left( r^4 d \var \wedge dt + t r dr
  \wedge dt + 2 t r^3 dr \wedge d \var \right).
\end{eqnarray*}
On $S$ we have $4 r^3 dr + 2 t dt = 0$,
so the last formula becomes
$$
  m(J,\theta) = 12 \pi A \oint_S \left( r^4 d \var \wedge dt - t^2 dt \wedge d \var \right)
  = 12 \pi A \oint_S  d \var \wedge dt = 48 \pi^2 A.
$$
Therefore we obtain the conclusion.
\end{pf}

\

\noindent Let $N$ be a three dimensional asymptotically flat pseudohermitian
manifold: we derive next an integral formula for the p-mass.

\begin{pro}\label{p:genfor} Let $(N,J,\th)$ be an   asymptotically flat pseudohermitian
manifold. Let $\beta : N \to \C$ be such that
\begin{equation}\label{eq:betagen}
     \beta = \ov{z} + \beta_{-1} + O(\rho^{-2+\e}) \quad \hbox{ near }
     \infty,
\end{equation}
and
\begin{equation}\label{eq:betagen2}
    \Box_b \beta = O(\rho^{-4}),
\end{equation}
where $\beta_{-1}$ is a term with the homogeneity of $\rho^{-1}$ satisfying
\begin{equation}\label{eq:assasybeta}
    (\beta_{-1})_{, \ou} = - 2 \sqrt{2} \pi A
  \frac{1}{\rho^2} - \frac{\sqrt{2} A}{|z|^2 + i t}
\end{equation}
near infinity, where $\e \in (0,1)$. Then one has
\begin{equation}\label{eq:1318}
    \frac{2}{3} m(J,\th) = - \int_N |\Box_b \beta|^2 \th \wedge d \th +
    2 \int_N |\beta_{, \ou \ou}|^2 \th \wedge d \th  + 2 \int_N R |\beta_{, \ou}|^2
    \th \wedge d \th + \frac 12 \int_N \ov{\beta} P \beta \, \th \wedge d \th,
\end{equation}
where
$$
  P \beta := \ov{\Box}_b \Box_b \beta + 4 i (A_{11} \beta_{, \ou})_{, \ou}
$$
is the CR Paneitz operator.
\end{pro}

\begin{rem}\label{r:assbeta}  We will show in Section \ref{s:pfthm} that
it is indeed possible to find a solution of \eqref{eq:betagen2} satisfying
\eqref{eq:betagen} and condition \eqref{eq:assasybeta}.
\end{rem}

\begin{pf}
By the commutation rules \eqref{eq:comm} we have that
$$
  \beta_{, \ou 1 \ou} = \beta_{, \ou \ou 1} +
  i \beta_{, \ou 0} - R \beta_{, \ou}.
$$
Multiply by $\ov{\beta}_{, 1}$  to obtain
\begin{eqnarray}\label{eq:nome} \nonumber
  \beta_{, \ou 1 \ou} \ov{\beta}_{, 1} +  R \beta_{, \ou}
\ov{\beta}_{, 1} & = & \beta_{, \ou \ou 1}
  \ov{\beta}_{, 1}  + i \beta_{, \ou 0} \ov{\beta}_{, 1} \\
   & = & (\beta_{, \ou \ou} \ov{\beta}_{, 1})_{, 1} - \beta_{, \ou \ou} \ov{\beta}_{, 1 1}
   + i (\beta_{, \ou 0} \ov{\beta})_{, 1} - i \beta_{, \ou 0 1} \ov{\beta}.
\end{eqnarray}
Integrating by parts we get
$$
  \int_N  (\beta_{, \ou \ou} \ov{\beta}_{, 1})_{, 1} \th \wedge d \th + \hbox{conj.}
  = i \oint_\infty (Z_{\ou} \beta_{, \ou}) \ov{\beta}_{, 1} \th^{\ou} \wedge \th + \hbox{conj.}
  + i \oint_\infty \o^1_1(Z_{\ou}) |\beta_{, \ou}|^2 \th^{\ou} \wedge \th + \hbox{conj.}.
$$
Regarding the first integral on the right-hand side of the
last formula we have that
$$
  i \oint_\infty (Z_{\ou} \beta_{, \ou})
\ov{\beta}_{, 1} \th^{\ou} \wedge \th + \hbox{conj.}
  = i \oint_\infty Z_{\ou} \beta_{, \ou} d \ov{z} \wedge
(dt - i \ov{z} dz) + \hbox{conj.},
$$
while for the second using the definition of $m(J,\th)$ we find
$$
    i \oint_\infty \o^1_1(Z_{\ou}) |\beta_{, \ou}|^2 \th^{\ou}
    \wedge \th + \hbox{conj.}  = \frac 12 m(J,\th).
$$
For the second divergence term in \eqref{eq:nome} and recalling
the asymptotics of $\o^1_1$ in Definition \ref{d:aflat} we get
\begin{eqnarray}\label{eq:nome2} \nonumber
  & & \int_N i (\beta_{, \ou 0} \ov{\beta})_{, 1} \th \wedge d \th +
  \hbox{conj.}  = - \oint_\infty \beta_{, \ou 0} \ov{\beta} \th^{\ou} \wedge
  \th + \hbox{conj.} \\ & = & - \oint_\infty (T \beta_{, \ou})
  \ov{\beta} \th^{\ou} \wedge \th + \hbox{conj.} + \oint_\infty \o^1_1(T)
  (\ov{\beta}_{, 1} \beta \th^1 \wedge \th - \beta_{, \ou} \ov{\beta} \th^{\ou}
  \wedge \th) \\ & = & - \oint_\infty (T \beta_{, \ou})
  \ov{\beta} \th^{\ou} \wedge \th + \hbox{conj.}. \nonumber
\end{eqnarray}
We have the Paneitz operator appearing if we use the formulas (the boundary terms
vanish after integration by parts)
\begin{eqnarray*}
  - \int_N \frac{1}{4} (P \beta) \ov{\beta} \th \wedge d \th + \hbox{conj.}
  & = &  - \int_N  (\beta_{, \ou 1 1} +
  i A_{11} \beta_{, \ou})_{, \ou} \ov{\beta} \th \wedge d \th + \hbox{conj.}   \\
   & = & - \int_N \beta_{, \ou 1} \ov{\beta}_{, \ou 1}   \th \wedge d \th
    + \int_N \left( - i (A_{11} \beta_{, \ou})_{, \ou}
  \ov{\beta} + \hbox{conj.} \right) \th \wedge d \th,
\end{eqnarray*}
and (from \eqref{eq:comm} and the identity $i \ov{\beta}_{,0}=\ov{\beta}_{,1{
\ov{1}}}-\ov{\beta}_{,{\ov{1}} 1}$)
\begin{eqnarray*}
 - i \int_N \beta_{, \ou 0 1} \ov{\beta} \th \wedge d \th + \hbox{conj.} & = &
  - i \int_N  \left( \beta_{, \ou 1 0} \ov{\beta} + (\beta_{, \ou} A_{11})_{, \ou}
  \ov{\beta} \right) \th \wedge d \th + \hbox{conj.} \\
   & = & \int_N \left( |\beta_{, \ou 1}|^2 - \frac{1}{4} (P \beta)
  \ov{\beta}  \right) \th \wedge d \th  + \hbox{conj.}.
\end{eqnarray*}
We will show next the following two identities
\begin{equation}\label{eq:claim1}
    i \oint_\infty Z_{\ou} \beta_{, \ou} d \ov{z} \wedge
    (dt - i \ov{z} dz) + \hbox{conj.} = \frac{7}{12} m(J,\th);
\end{equation}
\begin{equation}\label{eq:claim2}
     - \oint_\infty (T \beta_{, \ou}) z \sqrt{2} d \ov{z} \wedge
     (dt - i \ov{z} dz) + \hbox{conj.} = - \frac{5}{12} m(J,\th).
\end{equation}
To prove \eqref{eq:claim1} we notice that (modulo $O(\rho^{-3+\epsilon})$)
\begin{eqnarray*}
  \beta_{, \ou} &=& (1 - 2 \pi A \rho^{-2}) \stackrel{\circ}{Z}_{\ou}(\ov{z})
  + (\beta_{-1})_{, \ou} = \frac{1}{\sqrt{2}} - \sqrt{2} \pi A \rho^{-2}
  - 2 \sqrt{2} \pi A \frac{1}{\rho^2} - \frac{\sqrt{2} A}{|z|^2 + i t} \\
   & = & \frac{1}{\sqrt{2}} - 3 \sqrt{2} \pi A \frac{1}{\rho^2}
   - \frac{\sqrt{2}A}{|z|^2+it}.
\end{eqnarray*}
By elementary computations one finds
$$
  Z_{\ou} \beta_{, \ou} = - 3 \pi A z \frac{(it-|z|^2)}{\rho^{6}} + \frac{2Az(|z|^2-it)^2}{\rho^8} =
  \pi A z \frac{\mathfrak{F} + i t \mathfrak{G}}{\rho^{6}} + \frac{2Az(|z|^2-it)^2}{\rho^8},
$$
where we have set
$$
  \mathfrak{F} = 3 |z|^2; \qquad \quad \mathfrak{G} = -3.
$$
Then, by a scaling argument, we obtain
\begin{eqnarray}\label{eq:nome3} \nonumber
   & & i \oint_\infty Z_{\ou} \beta_{, \ou}  d \ov{z} \wedge (dt - i \ov{z} dz)
    + \hbox{conj.} \\
  & = & i \pi A \oint_S \left[ (\mathfrak{F} + i t \mathfrak{G})
    z d \ov{z} \wedge (dt - i \ov{z} dz) - (\mathfrak{F} - i t \mathfrak{G})
    \ov{z} d z \wedge (dt + i z d \ov{z}) \right]
    \\ & + & 2 A i \oint_S \left[ (|z|^2-it)^2
    z d \ov{z} \wedge (dt - i \ov{z} dz) - (|z|^2 + i t)^2
    \ov{z} d z \wedge (dt + i z d \ov{z}) \right]\nonumber
\end{eqnarray}
(recall that $S = \{ \rho = 1 \}$). By straightforward computations, and in
particular the fact that
\begin{equation}\label{eq:relpol2}
     z d \ov{z} + \ov{z} dz = 2 r dr; \qquad \qquad
     dz \wedge d \ov{z} = 2 i r d \var \wedge dr
\end{equation}
we find
\begin{eqnarray}\label{eq:nome4}
    i \oint_\infty Z_{\ou} \beta_{, \ou}  d \ov{z} \wedge (dt - i \ov{z} dz)
    + \hbox{conj.} & = & 2 \pi A \int_S \mathfrak{F} |z|^2 d \var \wedge dt - 2 \pi A
    \int_S \mathfrak{G} t^2 d \var \wedge dt \\ & + & 8 \pi A \int_{-1}^1 \sqrt{1-t^2} dt. \nonumber
\end{eqnarray}
Substituting for $\mathfrak{F}$ and $\mathfrak{G}$ we then get
\begin{eqnarray*}
  i \oint_\infty Z_{\ou} \beta_{, \ou}  d \ov{z} \wedge (dt - i \ov{z} dz)
    + \hbox{conj.} & = & 24 \pi^2 A +
    8 \pi A \int_{-1}^1 \sqrt{1-t^2} dt =
   28 \pi^2 A = \frac{7}{12} m(J,\th),
\end{eqnarray*}
which is \eqref{eq:claim1}.

To prove \eqref{eq:claim2} instead we use the following formula, which can be obtained from the 
expansion of $T$ in Remark \ref{r:regtor} 
\begin{eqnarray*}
  T \beta_{, \ou} & = & 3 \sqrt{2} \pi A  \frac{t}{\rho^{6}}
    + \frac{i A \sqrt{2}}{(|z|^2+it)^2}  +  O(\rho^{-5+\varepsilon}) \\
  & = & - \sqrt{2} \pi A i \frac{(\tilde{\mathfrak{F}} +
i t \tilde{\mathfrak{G}})}{\rho^{6}}
  +  \frac{iA\sqrt{2}(|z|^2-it)^2}{\rho^8} + O(\rho^{-5+\varepsilon}),
\end{eqnarray*}
where, using a notation similar to the previous one,
$$
   \tilde{\mathfrak{F}} = 0; \qquad \qquad \tilde{\mathfrak{G}} = 3.
$$
Using again a scaling and elementary computations (some of which similar
to those for the proof of \eqref{eq:claim1}) one gets
\begin{eqnarray*}
   & & - \oint_\infty (T \beta_{, \ou}) z \sqrt{2} d \ov{z} \wedge (dt - i \ov{z} dz)
   + \hbox{conj.} \\
   & = & 2 i \pi A \oint_S \left[ (\tilde{\mathfrak{F}} + i t \tilde{\mathfrak{G}})
    z d \ov{z} \wedge (dt - i \ov{z} dz) - (\tilde{\mathfrak{F}} - i t \tilde{\mathfrak{G}})
    \ov{z} d z \wedge (dt + i z d \ov{z}) \right]
     - \frac{1}{12} m(J,\th).
\end{eqnarray*}
Exactly as for \eqref{eq:nome3} and \eqref{eq:nome4} we deduce
\begin{eqnarray*}
  & & - \oint_\infty (T \beta_{, \ou}) z \sqrt{2} d \ov{z} \wedge (dt - i \ov{z} dz)
   + \hbox{conj.} \\ & = & 4 \pi A \int_S \tilde{\mathfrak{F}} |z|^2 d \var \wedge dt - 4 \pi A
    \int_S \tilde{\mathfrak{G}} t^2 d \var \wedge dt -
    \frac{1}{12} m(J,\th)
   \\ & = & - 16 \pi^2 A -
    \frac{1}{12} m(J,\th) = - \frac 13 m(J,\th)
   -\frac{1}{12} m(J,\th),
\end{eqnarray*}
which proves also \eqref{eq:claim2}.

Using all the above formulas and the fact that $P$ is real we then obtain
\begin{eqnarray*}
  & & - 2 \int |\beta_{, \ou 1}|^2 \th \wedge d \th + 2 \int_N R |\beta_{, \ou}|^2
  \th \wedge d \th \\ & = & \frac{2}{3} m(J,\th) - 2 \int_N |\beta_{, \ou \ou}|^2
  \th \wedge d \th + 2 \int_N |\beta_{, \ou 1}|^2 \th \wedge d \th - \frac 12
  \int_N \ov{\beta} (P \beta)  \, \th \wedge d \th.
\end{eqnarray*}
Since $\Box_b \beta = - 2 \beta_{, \ou 1}$, we finally get the conclusion.
\end{pf}

\section{Proof of Theorem \ref{t:pm}}\label{s:pfthm}

\noindent We are now in position to prove our main theorem. For doing this,
we need to introduce some notation and preliminary facts. For two functions
$f_1, f_2 \in L^1(\H^1)$ we define the convolution operator
$$
  (f_1 * f_2) (Z) = \int_{\H^1} f_1(Z W^{-1}) f_2(W) dW = \int_{\H^1}
  f_1(W) f_2(W^{-1} Z) dW; \qquad \quad
  Z = (z,t), \quad W = (w,s),
$$
with
$$
  W^{-1} Z = \left( z - w, t - s - 2 \hbox{ Im } \ov{z} w \right);
  \qquad \quad dW = \left( \stackrel{\circ}{\theta}
  \wedge \, d \stackrel{\circ}{\theta} \right) (W),
$$
see the notation in Subsection \ref{ss:notprel}.

By the solvability theory of $\stackrel{\circ}{\Box}_b$ in $\H^1$
(see Chapter 10 in \cite{CS}) we have that
\begin{equation}\label{eq:decompheis}
\stackrel{\circ}{\Box}_b \mathcal{K} = \mathcal{K} \stackrel{\circ}{\Box}_b
  = Id - \mathcal{S} \qquad \quad \hbox{ on } L^2(\H^1),
\end{equation}
where $\mathcal{S}$ is the Szeg\"{o} projection
$$
  \mathcal{S} h = \lim_{\e \to 0} \frac{1}{4 \pi^2}  h * \eta_\e^{-2};
  \qquad \qquad  \eta_\e = |z|^2 + \e^2 - i t,
$$
and where
$$
  \mathcal{K} \, h = h * \Phi; \qquad \qquad \Phi =
  \frac{1}{8 \pi^2} \log \left( \frac{|z|^2
  -  i t}{|z|^2 + i t} \right) (|z|^2 - i t)^{-1}.
$$
Define now
\begin{equation}\label{eq:deff}
   \mathfrak{f} = 4 \pi A \frac{\ov{z} (|z|^2 + i t)}{\rho^6},
\end{equation}
and let $\chi : \H^1 \to \R$ be a cutoff function, weakly monotone in $\rho$,
satisfying
$$
  \left\{
    \begin{array}{ll}
      \chi(z,t) = 0 & \hbox{in a neighborhood of } (0,0); \\
      \chi(z,t) = 1 & \hbox{near infinity.}
    \end{array}
  \right.
$$

\

\noindent Our next goal is to find a function $\beta_{-1}$
which has the homogeneity of $\rho^{-1}$
near infinity and such that
\begin{equation}\label{eq:eqbeta-1}
    \stackrel{\circ}{\Box}_b \beta_{-1} = -
  \chi \mathfrak{f} \qquad \hbox{ near infinity}.
\end{equation}
This function will be useful to solve \eqref{eq:betagen2} and to determine the
asymptotics of a function $\beta$ as in Proposition
\ref{p:genfor}. We can exhibit indeed two functions $\hat{g}$ and $\tilde{g}$ which
satisfy
\begin{equation}\label{eq:eqgi}
    \stackrel{\circ}{\Box}_b g = - 4  \pi A \frac{\ov{z} (|z|^2 + i t)}{\rho^6}
    \qquad \quad \hbox{ in } \H^1
\end{equation}
pointwise almost everywhere, but which are not continuous. The two functions are
\begin{equation}\label{eq:defg1}
  \hat{g}(z,t) = \left\{
             \begin{array}{ll}
              - 4 i \pi A  \left( \frac{\rho^2}{v z}
 - \frac 1 z \right)  & \hbox{ for } t > 0; \\
           - 4 i \pi A  \left( \frac{\rho^2}{v z}
 + \frac 1 z \right)   & \hbox{ for } t < 0,
             \end{array}
           \right.
\end{equation}
and
\begin{equation}\label{eq:defg2}
    \tilde{g}(z,t) = - 4 i \pi A \frac{\rho^2}{v z} \qquad \quad \hbox{ for } z \neq 0,
\end{equation}
where we have set
\begin{equation}\label{eq:defv}
    v = v(z,t) = t + i |z|^2.
\end{equation}
The function $\hat{g}$ is not continuous on the plane $\{t = 0\}$, while
$\tilde{g}$ is not continuous on the axis $\{z = 0\}$. By this reason,
$\hat{g}$ and $\tilde{g}$ should be considered spurious solutions of \eqref{eq:eqgi}:
however, it turns out that they are useful in some integral estimates, see
in particular Subsections \ref{ss:exbeta} and \ref{ss:exbeta2} below.

\

\noindent The proof of our main theorem relies on the next proposition, where
a solution of $\Box_b \beta = 0$ with a precise asymptotic behavior is found.

\begin{pro}\label{p:exbeta} If $N$ is an asymptotically flat
pseudohermitian manifold of dimension $3$, there exists a solution $\beta$
of
$$
\Box_b \beta = O(\rho^{-4}) \qquad \quad \hbox{ on } N
$$
satisfying \eqref{eq:betagen}. More precisely, if
$\e \in (0,1)$ and if $\mathfrak{f}$ is as in \eqref{eq:deff}, we can choose
$\beta_{-1}$ so that
\begin{equation}\label{eq:asybeta-1}
    \beta_{-1} = \mathcal{K} (- \chi \mathfrak{f}) + O(\rho^{-2+\e})
    \qquad \hbox{ for } \rho \hbox{ large, }
\end{equation}
and
\begin{equation}\label{eq:beta-1ou}
  (\beta_{-1})_{, \ou} =  - 2 \sqrt{2} \pi A
  \frac{1}{\rho^2} - \frac{\sqrt{2} A}{|z|^2 + i t} 
\end{equation}
near infinity. Moreover, if $N$ is the blow-up of an
embeddable compact three dimensional CR manifold with $\mathcal{Y}(J)>0$, then there exists
a $C^{\infty}$-smooth solution $\beta$ of $\Box_b{\beta}=0$ satisfying \eqref{eq:betagen}.
\end{pro}

\

\noindent The proof of Proposition \ref{p:exbeta} is given in the
next three subsections.

\subsection{Expansion of $\Box_b \ov{z}$ near infinity}\label{ss:expboxovz}

To study equation \eqref{eq:betagen2}, we  begin finding a correction to $\ov{z}$
by evaluating first $\Box_b$ on $\ov{z}$. The following asymptotic
expansion holds.

\begin{lem}\label{l:expboxbovz} If $N$ is an asymptotically flat
pseudohermitian manifold of dimension $3$, then
$$
  \Box_b \ov{z} = 4 \pi A \frac{\overline{z}}{\rho^6}
  (|z|^2 + i t) + O(\rho^{-4}) \qquad \hbox{ near infinity}.
$$
\end{lem}

\begin{pf} We want to express $Z_1$ and $\o^1_1$ in asymptotic
coordinates, recalling that
$$
   \stackrel{\circ}{Z}_1 = \frac{1}{\sqrt{2}} \left(\frac{\pa}{\pa z^1}
   + i z^{\ou} \frac{\pa}{\pa t} \right); \qquad \quad \stackrel{\circ}{\o}^1_1 = 0.
$$
Recall also that, from the  expansions in Section \ref{s:afman}
$$
  \th = \left( 1 + 4 \pi A \rho^{-2} + O(\rho^{-3})
  \right) \theta_0 + O(\rho^{-3}) dz  + O(\rho^{-3}) d\overline{z};
$$
$$
   \th^1 = \left( - 2 \sqrt{2} \pi A \frac{z v}{\rho^6} + O(\rho^{-5})
  \right) \theta_0 + O(\rho^{-4}) d\overline{z} + \left(
   1 + 2 \pi A \rho^{-2} + O(\rho^{-3}) \right) \sqrt{2} dz.
$$
Recalling also the estimates on the functions $a, b, c$ from Remark \ref{r:regtor} 
we also have
\begin{eqnarray*}
  Z_1 Z_{\ou} + \o^1_1(Z_1) Z_{\ou} & = & \left[ (1+a) \stackrel{\circ}{Z}_1
  + b \stackrel{\circ}{Z}_{\ou} + c \stackrel{\circ}{T} \right]
  \left[ (1+a) \stackrel{\circ}{Z}_{\ou} + b \stackrel{\circ}{Z}_1 + c
  \stackrel{\circ}{T} \right] + \o^1_1(Z_1) Z_{\ou}  \\
   & = & \stackrel{\circ}{Z}_1 \stackrel{\circ}{Z}_{\ou} + (\stackrel{\circ}{Z}_1 a)
   \stackrel{\circ}{Z}_{\ou} + 2 a \stackrel{\circ}{Z}_1 \stackrel{\circ}{Z}_{\ou}
  + \o^1_1(Z_1) Z_{\ou} + \mathcal{R}_{\Box_b},
\end{eqnarray*}
where the remainder term $\mathcal{R}_{\Box_b}$ satisfies
\begin{eqnarray*}
  \mathcal{R}_{\Box_b} & = & O(\rho^{-4}) \stackrel{\circ}{Z}_1 \stackrel{\circ}{Z}_{\ou}
  + O(\rho^{-4}) \stackrel{\circ}{Z}_{\ou}^2 + O(\rho^{-4}) \stackrel{\circ}{Z}_1^2
  + O(\rho^{-8}) \stackrel{\circ}{Z}_{\ou} \stackrel{\circ}{Z}_1 + O(\rho^{-3})
  Z_1 \stackrel{\circ}{T} + O(\rho^{-5}) \stackrel{\circ}{Z}_{\ou} \stackrel{\circ}{T}
  \\ & + & O(\rho^{-6}) \stackrel{\circ}{T}^2 + O(\rho^{-7}) \stackrel{\circ}{T}
  \stackrel{\circ}{Z}_1 + O(\rho^{-3})
  \stackrel{\circ}{T} \stackrel{\circ}{Z}_{\ou} + O(\rho^{-5}) \stackrel{\circ}{Z}_1 +
  O(\rho^{-5}) \stackrel{\circ}{Z}_{\ou} + O(\rho^{-4}) \stackrel{\circ}{T}.
\end{eqnarray*}
Recall that, by Definition \ref{d:aflat} one has
$$
  \o^1_1 = \left[ - \ov{g}_{-3} \sqrt{2} +
    O(\rho^{-4}) \right] dz + \left[ g_{-3} \sqrt{2}
    + O(\rho^{-4}) \right] d \ov{z} + O(\rho^{-5}) \th_0,
$$
where
$$
   g_{-3} = \frac{\pi A \sqrt{2}}{\rho^6} z (3 |z|^2 - 3 i t).
$$
%
%
%
%
%
%
%
%from which, comparing the coefficients of $dz \wedge d \ov{z}$ in
%\eqref{eq:dtheta1} and \eqref{eq:dth111}, we deduce
%\begin{eqnarray*}
%% \nonumber to remove numbering (before each equation)
%  2 g_{-3} & = & 2 \pi A \sqrt{2} \rho^{-6} z^2 \ov{z} - 4 \pi \sqrt{2} A i z v \rho^{-6}
%  = \frac{2 \pi A \sqrt{2}}{\rho^6} z ( |z|^2 -  2i v )  \\
%  & = & \frac{2 \pi A \sqrt{2}}{\rho^6} z \left( |z|^2 - 2 i (t + i |z|^2)
%  \right) = \frac{2 \pi A \sqrt{2}}{\rho^6} z (3 |z|^2 - 3 i t).
%\end{eqnarray*}
Noticing that $a = - 2 \pi A \rho^{-2} + O(\rho^{-3})$ implies
$$
  \stackrel{\circ}{Z}_1 a + O(\rho^{-4})
  = - 2 \pi A \stackrel{\circ}{Z}_1 \rho^{-2} + O(\rho^{-4})  = 2
  \pi A \frac{1}{\sqrt{2}} \frac{i \overline{z} \overline{v}}{\rho^6}
  + O(\rho^{-4}),
$$
we obtain the expansion
\begin{eqnarray*}
  \mathcal{E} & = & - 4 a \stackrel{\circ}{Z}_1 \stackrel{\circ}{Z}_{\ou} 
  - 2 (\stackrel{\circ}{Z}_1 a) \stackrel{\circ}{Z}_{\ou} + 2 \overline{g}_{-3}
  \stackrel{\circ}{Z}_{\ou} - 2 \mathcal{R}_{\Box_b} \\ & = & 8 \pi A \rho^{-2} \stackrel{\circ}{Z}_1
  \stackrel{\circ}{Z}_{\ou} - 4 \pi A \frac{1}{\sqrt{2}}
  \frac{i \overline{z} \overline{v}}{\rho^6} \stackrel{\circ}{Z}_{\ou} +
  2 \overline{g}_{-3} \stackrel{\circ}{Z}_{\ou} - 2 \mathcal{R}_{\Box_b} \\
   & = & 8 \pi A \rho^{-2} \stackrel{\circ}{Z}_1 \stackrel{\circ}{Z}_{\ou}
  + 2 \sqrt{2} \pi A \frac{\overline{z}}{\rho^6} (3 |z|^2 + 3 i t - i
  \overline{v}) \stackrel{\circ}{Z}_{\ou} - 2 \mathcal{R}_{\Box_b} \\ & = &
  8 \pi A \rho^{-2} \stackrel{\circ}{Z}_1
  \stackrel{\circ}{Z}_{\ou} + 4 \sqrt{2} \pi A \frac{\overline{z}}{\rho^6}
  (|z|^2 + i t) \stackrel{\circ}{Z}_{\ou} - 2 \mathcal{R}_{\Box_b}.
\end{eqnarray*}
Applying this operator to the function $\ov{z}$ we have
that
\begin{equation}\label{eq:expanboxovz}
    \mathcal{E} \overline{z} = 8 \pi A \rho^{-2} \stackrel{\circ}{Z}_1 (1) +
   4 \pi A \frac{\overline{z}}{\rho^6} (|z|^2 + i t) + O(\rho^{-4}) =
  4 \pi A \frac{\overline{z}}{\rho^6} (|z|^2 + i t) + O(\rho^{-4}),
\end{equation}
which gives the desired conclusion. \end{pf}

\subsection{Existence of a solution to $\Box_b \beta = 0$
and proof of \eqref{eq:asybeta-1}}\label{ss:exbeta}

\noindent In this Section we will prove the existence part under the stronger
assumption that $N$ is a blow-up as in the statement. For the general situation,
see Remark \ref{r:moregen}.

We have the following result concerning approximate orthogonality of $\mathfrak{f}$
to the Szeg\"{o} projection.

\begin{lem}\label{l:fperp} If $\mathfrak{f}$ is defined as in \eqref{eq:deff} then
$$
  (\mathcal{S} (\chi \mathfrak{f})) (z,t) = O(\rho^{-4}) \qquad
  \hbox{ as } \rho \to + \infty.
$$
\end{lem}

\begin{pf} Let $\hat{g}(z,t)$ be as in \eqref{eq:defg1}. Then one has
$$
  \stackrel{\circ}{\Box}_b \hat{g} = - \mathfrak{f}
  \qquad \hbox{ for } t \neq 0.
$$
Since
$$
\stackrel{\circ}{Z}_{\ou} \left( \frac 1 z \right) = 0; \qquad \qquad
\stackrel{\circ}{Z}_{\ou} \left( i \frac{\rho^2}{v z} \right)
= \rho^{-2}
$$
we have that $\stackrel{\circ}{Z}_{\ou} \hat{g}$ extends smoothly to $t = 0$
except for $\rho = 0$.  We consider the function
$$
  \hat{\mathfrak{f}} = - \stackrel{\circ}{\Box}_b (\chi \hat{g}),
$$
which coincides with $\mathfrak{f}$ outside a
compact set containing the origin of $\H^1$.

We then compute
\begin{eqnarray*}
  - (\mathcal{S} \hat{\mathfrak{f}})(Z) & = & - \frac{1}{4 \pi^2} \lim_{\e \to 0}
  \int_{\H^1} \hat{\mathfrak{f}}(W) \eta_\e^{-2}(W^{-1} Z) dW =
  - \frac{1}{2 \pi^2} \lim_{\e \to 0} \int_{\H^1} (\chi \hat{g})_{\ou 1}(W)
  \eta_\e^{-2}(W^{-1} Z) dW \\
  & = & - \frac{1}{2 \pi^2} \lim_{\e \to 0} \int_{\H^1} \left( (\chi \hat{g})_{\ou}(W)
  \eta_\e^{-2}(W^{-1} Z) \right)_{, 1} dW
  + \frac{1}{2 \pi^2} \lim_{\e \to 0} \int_{\H^1} (\chi \hat{g})_{\ou}(W)
  (\eta_\e^{-2}(W^{-1} Z))_{, 1} dW \\ &  = &
   - \frac{1}{2 \pi^2} \lim_{\e \to 0} \int_{\H^1} \left( (\chi \hat{g})_{\ou}(W)
  \eta_\e^{-2}(W^{-1} Z) \right)_{, 1} dW.
\end{eqnarray*}
%where we have set
%$$
%  \Psi(z,t) = \frac{2}{\pi} (|z|^2 - it)^{-2}.
%$$
In the above equalities we  used the fact that $\eta_\e^{-2}$ is a conjugate
CR function in $W$.
Integrating by parts we then get
\begin{eqnarray*}
  & & \int_{\H^1} \stackrel{\circ}{\Box}_b (\chi \hat{g})
  \eta_\e^{-2}(W^{-1} Z) dW \\ & = &
  \frac{i}{2 \pi^2} \lim_{\e \to 0} \left[
  \oint_{t = 0^+} (\chi \hat{g})_{\ou}(W) \eta_\e^{-2}(W^{-1} Z)
  (\th \wedge \th^{\ou})(W) - \oint_{t = 0^-} (\chi \hat{g})_{\ou}(W)
  \eta_\e^{-2}(W^{-1} Z) (\th \wedge \th^{\ou})(W) \right]
 \\ & = & \frac{1}{2 \pi^2} \lim_{\e \to 0} \left[ \oint_{t = 0^+}
 (\chi \hat{g})_{\ou}(W) \eta_\e^{-2}(W^{-1} Z)
  \frac{\ov{z}}{\sqrt{2}} dw \wedge d \ov{w} -
 \oint_{t = 0^-} (\chi \hat{g})_{\ou}(W) \eta_\e^{-2}(W^{-1} Z)
 \frac{\ov{z}}{\sqrt{2}} dw \wedge d \ov{w} \right].
\end{eqnarray*}
Since $(\chi \hat{g})_{\ou}$ is smooth outside a compact set, we have vanishing
of the difference of the boundary integrands outside a compact set. We also have that
$$
  \eta_\e^{-2}(W^{-1} Z) = O(\rho^{-4})(Z) \quad \hbox{ for } W \hbox{ in
  any given compact set of } \H^1.
$$
The last two formulas then imply
$$
  (\mathcal{S} \hat{\mathfrak{f}})(Z) = O(\rho^{-4})(Z).
$$
By the same reason, since $\mathfrak{f}$ and $\hat{\mathfrak{f}}$ coincide outside
a compact set we have that as well
$$
  (\mathcal{S} \mathfrak{f})(Z) = O(\rho^{-4})(Z),
$$
which is the desired conclusion. \end{pf}

\

\noindent Below, for $\mu \in \R$ we denote by $\mathcal{E}(\rho^{\mu})$ the set of smooth
functions $u$ on an asymptotically flat pseudohermitian manifold $N$ for which near infinity
one has
$$
  |\mathcal{Z}^{(\a)} u| \leq C_{\a,u}  \, \rho^{\mu - |\a|},
$$
where $\a$ is a multi-index and  $\mathcal{Z}^{(\a)}$ is a composition of $Z_1$
and $Z_{\ou}$ derivatives of order $|\a|$, and where $C_{\a,u}$ is a
positive constant depending on $\a$ and $u$. 
The main tool to prove the existence part in Proposition \ref{p:exbeta} is the
next result (notice that in \cite{HY} a different convention for $\theta$ is used).
Define
$$
\Box_{b,1}:=G_{p}^{2}\Box_{b}; \quad m_{1}:=G_{p}^{-2}\theta\wedge d\theta .
$$ 
Denote the space of square integrable functions with respect to the volume form
$m_{1}$ by  $L^2(m_{1})$.

\begin{thm}\label{t:boxb} (Theorem 1.3 and Theorem 1.4 in \cite{HY}) Suppose $N$ is the blow-up of
an embeddable compact three dimensional CR manifold with $\mathcal{Y}(J) > 0$. Then 
$$
  \Box_{b,1} : \hbox{Dom}(\Box_{b,1}) \subseteq L^2(m_{1}) \to L^2(m_{1}) 
$$
has closed range, and the following decomposition holds 
$$
  \Box_{b,1} \, \mathcal{K}_N + \mathcal{S}_N = Id \qquad \quad \hbox{ on } L^2(m_{1}), 
$$
where $\mathcal{K}_N$ denotes a partial inverse of $\Box_{b,1}$ and $\mathcal{S}_N$ the Szeg\"o 
projection onto the kernel of $\Box_{b,1}$. Moreover, for every $\e \in (0,2)$, $\mathcal{K}_N$ and 
$\mathcal{S}_N$ can be extended continuously as maps 
$$
  \mathcal{K}_N : \mathcal{E}(\rho^{2-\e}) \to \mathcal{E}(\rho^{-\e}); \qquad 
  \qquad \mathcal{S}_N : \mathcal{E}(\rho^{2-\e}) \to \mathcal{E}(\rho^{2-\e}) 
$$
and hence there holds
$$
\Box_{b,1} \, \mathcal{K}_N + \mathcal{S}_N = Id \qquad \quad \hbox{ on }\mathcal{E}(\rho^{2-\e}).  
$$
\end{thm}

\begin{rem} 
 Note that $\Box_b$ may not have $L^2$ closed range with respect to the volume form
 $\theta\wedge d\theta$. This causes difficulty to solve the $\Box_b$ equation 
 directly. In \cite{HY} the authors introduce the above weighted Kohn Laplacian
 $\Box_{b,1}$ and show that $\Box_{b,1}$ has $L^2$ closed range with respect to
 the weighted volume form $m_1$.        
\end{rem}

\

\noindent We apply Theorem \ref{t:boxb} with
$$
  h = - \Box_{b,1} (\chi \ov{z} + \mathcal{K} (- \chi \mathfrak{f})).
$$
By Lemmas \ref{l:expboxbovz} and \ref{l:fperp} (in fact, $O(\rho^{-4})$ in these two lemmas can be replaced by $\mathcal{E}(\rho^{-4})$) we have that
$h \in \mathcal{E}(\rho^{\e })$ and hence we can find $u= \mathcal{K}_{N}h \in \mathcal{E}(\rho^{-2+\e})$
such that  $\Box_{b,1}u = h$ (note that $\mathcal{S}_{N}h=0$, see Theorem 1.5 in \cite{HY}). Therefore the function
$$
   \beta = \chi \ov{z} + \mathcal{K} (- \chi \mathfrak{f}) + u
$$
is a solution of $\Box_{b,1}\beta=0$, hence $\Box_b \beta = 0$ (see also the deduction in the introduction of \cite{HY}).
By scaling arguments, one finds that $\mathcal{K}
(\chi \mathfrak{f})$ is the sum of a function with the homogeneity of $\rho^{-1}$ and
a function of order $\rho^{-2+\e}$, and therefore we deduce \eqref{eq:asybeta-1}.

\subsection{Proof of \eqref{eq:beta-1ou}}\label{ss:exbeta2}

Since
$$
  (\tilde{g})_{\ou} =  - 2 \sqrt{2} \pi A
  \frac{1}{\rho^2}  \qquad \quad (z \neq 0)
$$
extends continuously to $\H^1 \setminus \{0\}$, to show \eqref{eq:beta-1ou}
it is sufficient to prove that
\begin{equation}\label{eq:diffder}
    (\mathcal{K} (- \chi \mathfrak{f}))_{\ou} + \frac{\sqrt{2} A}{|z|^2 + i t}
    = (\tilde{g})_{\ou}
    + O(\rho^{-3}) \qquad \hbox{ for $z \neq 0$ and $\rho$ large}.
\end{equation}
To obtain this estimate we use a distributional argument. First of all,
we remark that $\frac 1 z$ is well defined, as a distribution, by
the formula
$$
  \langle \frac 1 z,  \var \rangle = - \int_{\H^1} \log |z|^2 \frac{\pa
  \ov{\var}}{\pa z} \stackrel{\circ}{\th} \wedge d \stackrel{\circ}{\th},
$$
for any  test function  $\var \in C^\infty_c(\H^1)$. We have now the following result.

\begin{lem}\label{l:Boxbdistr} If $\tilde{g}$ is defined as in
\eqref{eq:defg2} then one has
\begin{equation}\label{eq:ld1}
    \stackrel{\circ}{\Box}_b \tilde{g} = - \mathfrak{f} - 4 \sqrt{2}
    i \pi^2 A \frac{\rho^2}{v} Z_1
    \d_{\{z = 0\}} \qquad \hbox{ in the distributional sense of }
    \H^1  \setminus \{0\}.
\end{equation}
\end{lem}

\begin{pf}
We have
\begin{equation}\label{eq:distr1}
     Z_{\ou} \tilde{g} = - 4 i \pi A Z_{\ou} \left( \frac{\rho^2}{v} \right) \frac 1 z
  - 4 i \pi A \frac{\rho^2}{v} Z_{\ou} \left( \frac 1 z \right).
\end{equation}
We next notice that, for any  function $\var \in C^\infty_c(\H^1)$
\begin{eqnarray*}
  \langle Z_{\ou} \frac 1 z, \var \rangle & = & \langle \frac 1 z,
  Z_{\ou}^* \var \rangle = \langle \frac 1 z, - Z_1 \var \rangle =
  \frac{1}{\sqrt{2}} \int_{\H^1} \log |z|^2 \ov{\left[ \pa_{\ov{z}}
  (\pa_z \var + i \ov{z} \pa_t \var) \right]} \stackrel{\circ}{\th}
  \wedge d \stackrel{\circ}{\th} \\ & = & \frac{1}{\sqrt{2}} \int_{\H^1} \log
  |z|^2 \ov{\frac 14 \D_{\R^2} \var} \stackrel{\circ}{\th} \wedge d \stackrel{\circ}{\th}
  = \frac{1}{4 \sqrt{2}} \langle \D_{\R^2} \log |z|^2, \var \rangle \\
  & = & \frac{1}{4 \sqrt{2}} \cdot 4 \pi \langle \d_{\{z = 0\}}, \var \rangle
  = \frac{\pi}{\sqrt{2}} \langle  \d_{\{z = 0\}}, \var \rangle.
\end{eqnarray*}
Therefore by \eqref{eq:distr1} we have
$$
  Z_{\ou} \tilde{g} = -4 i \pi A Z_{\ou} \left( \frac{\rho^2}{v} \right) \frac 1 z
  + 4 \pi A i \pi \frac{1}{\sqrt{2}} \frac{\rho^2}{v} \d_{\{z = 0\}}
  \qquad ((z,t) \neq (0,0)).
$$
Using the (distributional) formulas
$$
   Z_1 \frac 1z = - \frac{1}{\sqrt{2}} \frac{1}{z^2}; \qquad \quad
   Z_1 \left( \frac{\rho^2}{v} \right) = - \frac{i}{\sqrt{2}}
   \frac{z(t - i |z|^2)}{\rho^6}
$$
and the fact that $\ov{z} \d_{\{z = 0\}} = 0$ we find
\begin{eqnarray*}
  - 2 Z_1 Z_{\ou} \tilde{g} & = & - \mathfrak{f} - 4 \sqrt{2} \pi^2 i A Z_1
  \left( \frac{\rho^2}{v} \right) \d_{\{z = 0\}} - 4 \sqrt{2} \pi^2 i A
  \frac{\rho^2}{v} Z_1  \d_{\{z = 0\}} \\ & = & - \mathfrak{f} - 4 \sqrt{2} \pi^2
  i A \frac{\rho^2}{v} Z_1  \d_{\{z = 0\}} \quad \qquad ((z,t) \neq (0,0)),
\end{eqnarray*}
which is the desired conclusion. \end{pf}

\

\noindent The previous lemma implies clearly that
$$
  \Box_b (\chi \tilde{g}) = \chi \left( - \mathfrak{f} - 4 \sqrt{2} i \pi^2
  A \frac{\rho^2}{v} Z_1 \d_{\{z = 0\}} \right) + \tau,
$$
where $\tau$ is a distribution with compact support $A_\tau$.

We next consider a sequence of smooth nonincreasing cutoff functions
$\chi_j = \chi_j(\rho)$ such that
$$
  \chi_j(\rho) = 1 \quad \hbox{ for } \rho \leq j; \qquad
  \quad \chi_j(\rho) = 0 \quad \hbox{ for } \rho \geq j + 1,
$$
and another sequence of nondecreasing cutoff functions $\tilde{\chi}_j
= \tilde{\chi}_j(|z|)$ satisfying
$$
  \tilde{\chi}_j(|z|) = 1 \quad \hbox{ for } |z| \geq \frac{1}{j}; \qquad
  \quad \tilde{\chi}_j(|z|) = 0 \quad \hbox{ for } |z| \leq \frac{1}{2j}.
$$
We then define the sequence
$$
  \tilde{g}_j(z,t) = \tilde{\chi}_j(|z|) \chi_j(\rho) \chi(z,t) \tilde{g}(z,t).
$$
Since the $\tilde{g}_j$'s converge to $\tilde{g}$ locally in $L^1$ we have that
\begin{equation}\label{eq:boxbtgj}
    \stackrel{\circ}{\Box}_b \tilde{g}_j = - \chi_j \chi \mathfrak{f}
    + \zeta_j + \xi_j,
\end{equation}
where $\zeta_j$ is a sequence of functions supported in $\{ 1 \leq
\rho \leq j + 1\} \cap \{ |z| \leq 1/j\}$ satisfying
\begin{equation}\label{eq:distrz}
    \zeta_j \stackrel{\mathcal{D}'}{\to} - 4 \sqrt{2}
    i  \pi^2 A \frac{\rho^2}{v} Z_1 \d_{\{z = 0\}},
\end{equation}
on every open bounded set of $\H^1$ which does not intersect $\{\rho \geq 2\}$,
and where $\xi_j$ is supported in $\{ 1 \leq \rho \leq 2\} \cup
\{ j \leq \rho \leq j+1 \}$, with order $O(\rho^{-3})$.
For $z \neq 0$ we now compute
\begin{eqnarray}\label{eq:3terms} \nonumber
  \stackrel{_{(Z)}}{Z}_{\ou} (\mathcal{K} (\stackrel{\circ}{\Box}_b
  \tilde{g}_j)) (Z) & = &
  \int_{\H^1} (\stackrel{\circ}{\Box}_b \tilde{g}_j) (W)
  \stackrel{_{(Z)}}{Z}_{\ou} \Phi(W^{-1} Z)
  dW = \int_{\H^1} (\chi_j \chi \mathfrak{f}) (W) \stackrel{_{(Z)}}{Z}_{\ou}
  \Phi(W^{-1} Z) dW \\ & + & \int_{\H^1} \zeta_j (W) \stackrel{_{(Z)}}{Z}_{\ou}
  \Phi(W^{-1} Z) dW +
  \int_{\H^1} \xi_j (W) \stackrel{_{(Z)}}{Z}_{\ou} \Phi(W^{-1} Z) dW.
\end{eqnarray}

\begin{lem}\label{l:3form} The following three formulas hold
\begin{equation}\label{eq:num1}
    \int_{\H^1} (\chi_j \chi \mathfrak{f}) (W) \stackrel{_{(Z)}}{Z}_{\ou}
  \Phi(W^{-1} Z) dW = Z_{\ou} \mathcal{K} (\chi \mathfrak{f}) (Z) +
  o_j(1) O(\rho^{-3});
\end{equation}
\begin{equation}\label{eq:num2}
   \int_{\H^1} \zeta_j (W) \stackrel{_{(Z)}}{Z}_{\ou} \Phi(W^{-1} Z) dW
   = \frac{\sqrt{2} A}{|z|^2 + i t} + O(\rho^{-3}) + o_j(1);
\end{equation}
\begin{equation}\label{eq:num3}
   \int_{\H^1} \xi_j (W) \stackrel{_{(Z)}}{Z}_{\ou}
   \Phi(W^{-1} Z) dW = O(\rho^{-3}) + o_j(1).
\end{equation}
\end{lem}

\begin{pf}  To prove \eqref{eq:num1}  we notice that clearly
\begin{eqnarray*}
  \int_{\H^1} (\chi_j \chi \mathfrak{f}) (W) \stackrel{_{(Z)}}{Z}_{\ou}
  \Phi(W^{-1} Z) dW & = & \int_{\H^1} (\chi \mathfrak{f}) (W)
  \stackrel{_{(Z)}}{Z}_{\ou} \Phi(W^{-1} Z) dW \\
  & + & \int_{\H^1} ((1-\chi_j) \chi \mathfrak{f}) (W)
  \stackrel{_{(Z)}}{Z}_{\ou} \Phi(W^{-1} Z) dW,
\end{eqnarray*}
so it is sufficient to show
$$
  \int_{\H^1} ((1-\chi_j) \chi \mathfrak{f}) (W) \stackrel{_{(Z)}}{Z}_{\ou}
  \Phi(W^{-1} Z) dW = o_j(1) O(\rho^{-3}).
$$
We first notice that
\begin{equation}\label{eq:decphizw}\nonumber
  \stackrel{_{(Z)}}{Z}_{\ou} \Phi(W^{-1} Z) = \frac{\sqrt{2}}{8 \pi^2}
   \frac{w-z}{(\rho(W^{-1} Z))^4} = O(\rho^{-3}(W^{-1} Z)),
\end{equation}
which implies
$$
  \left| \int_{\H^1} ((1-\chi_j) \chi \mathfrak{f}) (W)
  \stackrel{_{(Z)}}{Z}_{\ou} \Phi(W^{-1} Z) \right|
  \leq C \int_{\{\rho(W) \geq j\}} (\rho(W))^{-3}
  (\rho(W^{-1} Z))^{-3} dW.
$$
By a scaling argument one shows that
$$
  \int_{\{\rho(W) \geq j\}} (\rho(W))^{-3}
  (\rho(W^{-1} Z))^{-3} dW \leq C (\rho(Z))^{-2}
  \int_{\left\{ \rho(W) \geq \frac{j}{\rho(Z)} \right\}}
  (\rho(W))^{-6} dW \leq \frac Cj (\rho(Z))^{-3},
$$
giving \eqref{eq:num1}. Recalling the properties of the support of
$\zeta_j$, we can write that
\begin{eqnarray*}
  \int_{\H^1} \xi_j (W) \stackrel{_{(Z)}}{Z}_{\ou} \Phi(W^{-1} Z)
   dW & = & \int_{\{ 1 \leq \rho(W) \leq 2 \}} \xi_j (W)
   \stackrel{_{(Z)}}{Z}_{\ou} \Phi(W^{-1} Z) dW \\
  & + & \int_{\{ j \leq \rho(W) \leq j+1 \}} \xi_j (W)
  \stackrel{_{(Z)}}{Z}_{\ou} \Phi(W^{-1} Z) dW.
\end{eqnarray*}
An argument similar to the one for \eqref{eq:num1} gives the same
estimate for the second integral in the last formula. To control the first
one we use \eqref{eq:decphizw} to find
$$
  \left| \int_{\{ 1 \leq \rho(W) \leq 2 \}} \xi_j (W)
   \stackrel{_{(Z)}}{Z}_{\ou} \Phi(W^{-1} Z) dW \right|
   \leq C \int_{\{ 1 \leq \rho(W) \leq 2 \}} (\rho(W^{-1} Z))^{-3} dW
   = O((\rho(Z))^{-3}),
$$
which gives \eqref{eq:num3}.

It remains to show \eqref{eq:num2}. In this case, we divide the integral
into the three regions
$$
  \mathcal{A}_1 = \{ 1 \leq \rho \leq 4 \}; \qquad \quad
  \mathcal{A}_2 = \{ 4 \leq \rho \leq j \}; \qquad \quad
  \mathcal{A}_3 = \{ j \leq \rho \leq j+1 \}.
$$
Using scaling arguments similar to the previous ones, \eqref{eq:distrz},
and the fact that
$$
  \stackrel{_{(W)}}{Z}_{1} \stackrel{_{(Z)}}{Z}_{\ou} \Phi(W^{-1} Z) =
  - \frac{1}{8 \pi^2} \frac{1}{(|z-w|^2 + i (t - s - 2 \hbox{Im}(\ov{z}w)))^2},
$$
we easily obtain
$$
  \left| \int_{\mathcal{A}_1} \zeta_j (W) \stackrel{_{(Z)}}{Z}_{\ou}
  \Phi(W^{-1} Z) dW \right| = O((\rho(Z))^{-3}); \qquad \quad
  \left| \int_{\mathcal{A}_3} \zeta_j (W) \stackrel{_{(Z)}}{Z}_{\ou}
  \Phi(W^{-1} Z) dW \right| = o_j(1).
$$
On the other hand, from \eqref{eq:distrz} and some elementary estimates
one finds
$$
  \int_{\mathcal{A}_2} \zeta_j (W) \stackrel{_{(Z)}}{Z}_{\ou}
  \Phi(W^{-1} Z) dW = - 4 \sqrt{2} i \pi^2 A \int_{\mathcal{A}_2}
  \frac{\rho^2}{v}  \stackrel{_{(W)}}{Z}_{1} \d_{\{w = 0\}}
  \stackrel{_{(Z)}}{Z}_{\ou} \Phi(W^{-1} Z) dW + o_j(1).
$$
Using the facts that
$$
  Z_1 \left( \frac{\rho^2}{s + i |w|^2} \right) = - i \frac{1}{\sqrt{2}}
  \frac{\ov{w} (s - i |w|^2)}{\rho^6} = 0; \qquad
   \frac{\rho^2}{s + i |w|^2} =
   \left\{
     \begin{array}{ll}
       1, & \hbox{ for } s > 0; \\
       -1, & \hbox{for } s < 0,
     \end{array}
   \right.
   \qquad \qquad \hbox{ if } w = 0,
$$
we deduce
\begin{eqnarray*}
  \int_{\mathcal{A}_2} \zeta_j (W) \stackrel{_{(Z)}}{Z}_{\ou}
  \Phi(W^{-1} Z) dW & = & - 4 \sqrt{2} i \pi^2 A
  \int_{\mathcal{A}_2} \frac{\rho^2(W)}{s + i |w|^2} \left( \stackrel{_{(W)}}{Z}_1
   \d_{\{w = 0\}} \right) \stackrel{_{(Z)}}{Z}_{\ou}
   \Phi(W^{-1} Z) dW + o_j(1) \\ & = & 4 \sqrt{2} i \pi^2 A
    \int_{\mathcal{A}_2} \d_{\{w = 0\}}
  \left\{ \frac{s}{|s|} \stackrel{_{(W)}}{Z}_1 \stackrel{_{(Z)}}{Z}_{\ou}
  (\Phi(W^{-1} Z)) \right\} dW + o_j(1).
\end{eqnarray*}
Using also the identity
$$
  \stackrel{_{(W)}}{Z}_{1} \stackrel{_{(Z)}}{Z}_{\ou} \Phi(W^{-1} Z) =
  - \frac{1}{8 \pi^2} \frac{1}{(|z|^2 + i(t-s))^2}, \qquad \qquad \hbox{ for }
  w = 0,
$$
we then find
\begin{eqnarray*}
  \int_{\mathcal{A}_2} \zeta_j (W) \stackrel{_{(Z)}}{Z}_{\ou}
  \Phi(W^{-1} Z) dW & = & - \frac{\sqrt{2}}{2} i A \int_{4}^j
  \frac{1}{(|z|^2 + i(t-s))^2} ds \\
  & + & \frac{\sqrt{2}}{2} i A \int_{-j}^{-4}
  \frac{1}{(|z|^2 + i(t-s))^2} ds + o_j(1) \\ & = &
  \frac{\sqrt{2} A}{|z|^2 + i t} + O(\rho^{-4}) + o_j(1),
\end{eqnarray*}
which gives \eqref{eq:num2}.
\end{pf}

\

\noindent Summing up, by \eqref{eq:3terms} and Lemma \ref{l:3form}
we obtained that
$$
  Z_{\ou} (\mathcal{K} (\stackrel{\circ}{\Box}_b
  \tilde{g}_j)) = Z_{\ou} \mathcal{K} (- \chi \mathfrak{f})
  + \frac{\sqrt{2} A}{|z|^2 + i t} + O(\rho^{-3}) + o_j(1).
$$
Since $\mathcal{K} \stackrel{\circ}{\Box}_b = Id - \mathcal{S}$ on
distributions with compact support and since the image of the Szeg\"{o} projection consists
of CR functions, by Lemma \ref{l:fperp} we deduce that for $z \neq 0$
$$
  Z_{\ou} \tilde{g}_j = Z_{\ou} \tilde{g}_j - Z_{\ou} (\mathcal{S} \tilde{g}_j)
  = Z_{\ou} \mathcal{K} (- \chi \mathfrak{f})
  + \frac{\sqrt{2} A}{|z|^2 + i t} + O(\rho^{-3}) + o_j(1).
$$
Letting $j \to + \infty$, this means
$$
  Z_{\ou} \tilde{g} = Z_{\ou} \mathcal{K} (\chi \mathfrak{f})
  + \frac{\sqrt{2} A}{|z|^2 + i t} + O(\rho^{-3})  \qquad \quad
  \hbox{ for } z \neq 0 \hbox{ and } \rho \hbox{ large}.
$$
This is precisely \eqref{eq:diffder}, concluding the proof of \eqref{eq:beta-1ou}.

\begin{rem}\label{r:moregen} Under the more general assumption of Proposition \ref{p:exbeta},
it is sufficient to apply all the above arguments except for Theorem \ref{t:boxb},
which allows to pass from a solution of $\Box_b \beta = O(\rho^{-4})$ to a solution
of $\Box_b \beta = 0$.
\end{rem}

%
%
%\begin{rem}\label{r:uniq}
%The solution $\beta$ is unique up to a CR function: therefore its $\ou$ derivative
%at infinity is uniquely determined as in \eqref{eq:assasybeta}.
%\end{rem}

\subsection{The role of positivity of the Paneitz operator}\label{ss:P>0}

We prove here that the positivity of $P$ on the compact manifold implies
positivity of the last term in \eqref{eq:1318} as well, up to a multiple of
$m(J,\th)$ (which can be reabsorbed into the left-hand side). We have the following proposition, 
which also exploits some result from \cite{HY}.

\begin{pro}\label{p:P>0} Let $N$ be the blow-up of a three dimensional CR manifold
with positive Tanaka-Webster invariant and non-negative CR Paneitz operator. Let $\beta$
be as in Proposition \ref{p:genfor}: then one has
\begin{equation}\label{3.9.1}
    \frac 12 \int_N{\bar\beta}P{\beta}{\theta}{\wedge}d{\theta} =
\frac 12 \int_N (\ov{\beta}_{-1} + \ov{\beta}_{-2+\e} + u)
  P (\beta_{-1} + \beta_{-2+\e} + \ov{u}) \, \theta \wedge d \theta
  - \frac{2}{3} m(J,\th)
\end{equation}
where $u=O(\rho^{-3})$ and
\begin{equation}\label{3.9.2}
    \frac{4}{3} m(J,\th) \geq - \int_N |\Box_b \beta|^2 \th \wedge d \th
  + 2 \int_N |\beta_{, \ou \ou}|^2 \th \wedge d \th + 2 \int_N R
  |\beta_{, \ou}|^2  \th \wedge d \th.
\end{equation}
Moreover, if $\theta$ is chosen as in \eqref{eq:bbuu} and $\beta$ is chosen so that
$\Box_b \beta = 0$ (this can be achieved by Proposition \ref{p:exbeta}), then $R = 0$ and
the above inequality is reduced to
\begin{equation}\label{3.9.3}
    \frac 43 m(J,\theta) \geq
2\int_N|\beta_{,{\bar 1}{\bar 1}}|^{2}{\theta}{\wedge}d{\theta} \geq 0.
\end{equation}
\end{pro}

\begin{pf} We set
$$
  f = \ov{\pa}_b z \in \mathcal{E}(\rho^{-4}).
$$
By Theorem 4.5 in \cite{HY}, we can find $u\in\mathcal{E}(\rho^{-3})$ so that $\ov{\pa}_b u = \ov{\pa}_b z$ 
(note that in \cite{HY}, the coordinates that authors use can be changed to ours by a CR inversion).
Let us write $\beta$ as
$$
  \beta = \ov{z} + \beta_{-1} + \beta_{-2+\e} = (\ov{z} - \ov{u})
  + \beta_{-1} + \beta_{-2+\e} + \ov{u},
$$
so
$$
  \ov{\beta} = (z-u) + \ov{\beta}_{-1} + \ov{\beta}_{-2+\e} + u.
$$
We write the Paneitz operator $P \psi = 4 (\psi_{\ou 1 1} + i A_{11}
\psi_{\ou})_{\ou}$ as
$$
  P = P_3 \circ \ov{\pa}_b \qquad \quad \hbox{ or }
  \quad \qquad  P = \pa_b^* \circ \tilde{P}_3.
$$
where 
$$
P_3 (a_{\bar 1}\theta^{\bar 1})=4(a_{{\bar 1},11}+iA_{11}a_{\bar 1})_{,{\bar 1}},
$$
and
$$
\tilde{P}_3 (\psi )= 4 (\psi_{\ou 1 1} + i A_{11}
\psi_{\ou})\theta^{1}.
$$
Since $P$ is real we have that $\int_N \ov{\beta} P \beta \, \th \wedge
d \th = \int_N \beta  P \ov{\beta} \, \th \wedge d \th \in \R$, and hence
\begin{eqnarray*}
   & & \int_N \left[ (\ov{z} - \ov{u}) + \beta_{-1} + \beta_{-2+\e} +
   \ov{u} \right] P (\ov{\beta}_{-1} + \ov{\beta}_{-2+\e} + u) \, \th
   \wedge d \th \\ & = & \int_N \left[ (z-u) + \ov{\beta}_{-1} +
   \ov{\beta}_{-2+\e} + u \right] P (\beta_{-1} + \beta_{-2+\e}
   + \ov{u}) \, \th \wedge d \th.
\end{eqnarray*}
Integrating by parts, we move a $\ou$ derivative to the first term in the
integrand of the right-hand side, which annihilates $z-u$. In this way we get
\begin{eqnarray*}
  & & \int_N \left[ (\ov{z} - \ov{u}) + \beta_{-1} + \beta_{-2+\e}
  + \ov{u} \right] P (\ov{\beta}_{-1} + \ov{\beta}_{-2+\e} + u) \, \th
   \wedge d \th \\ & = &
  - \int_N (\ov{\beta}_{-1} + \ov{\beta}_{-2+\e} + u)_{\ou}
  \tilde{P}_3 (\beta_{-1} + \beta_{-2+\e} + \ov{u}) \, \th
   \wedge d \th + \oint_\infty z i
  \tilde{P}_3  \beta_{-1} \, \theta \wedge \theta^1.
\end{eqnarray*}
For the last term we used the decay at infinity of the lower order terms
in $\beta$, giving
$$
  \oint_\infty i ((z-u) + \ov{\beta}_{-1} + \ov{\beta}_{-2+\e} + u)
  \tilde{P}_3 (\beta_{-1} + \beta_{-2+\e} + \ov{u}) \, \theta \wedge \theta^1
  = \oint_\infty z i \tilde{P}_3 \beta_{-1} \, \theta \wedge \theta^1.
$$
Using the Stokes theorem again, we find
\begin{eqnarray*}
  \oint_\infty z i \tilde{P}_3 \beta_{-1} \, \theta \wedge \theta^1 & = &
  \int_N ((z-u) + \ov{\beta}_{-1} + \ov{\beta}_{-2+\e} + u) P (\beta_{-1}
  + \beta_{-2+\e} + \ov{u}) \, \theta \wedge d \theta  \\
  & + & \int_N (\ov{\beta}_{-1} + \ov{\beta}_{-2+\e} + u)_{\ou} \tilde{P}_3
  (\beta_{-1} + \beta_{-2+\e} + \ov{u}) \, \theta \wedge d \theta.
\end{eqnarray*}
From the last formulas we deduce
\begin{eqnarray}\label{eq:miao1} \nonumber
  \int_N \beta P \ov{\beta} \, \th \wedge d \th & = & \int_N
  (\ov{\beta}_{-1} + \ov{\beta}_{-2+\e} + u) P
  (\beta_{-1} + \beta_{-2+\e} + \ov{u}) \, \th \wedge d \th \\
  & + &  \oint_\infty z i \tilde{P}_3 \beta_{-1} \, (dt
  - i \ov{z} d z + i z d \ov{z}) \wedge \sqrt{2} dz.
\end{eqnarray}
Using the decay of $A_{11}$ at infinity we find
\begin{equation}\label{eq:miao2}
     \oint_\infty \tilde{P}_3 \beta_{-1} \cdot z i (dt - i \ov{z} d z + i z
  d \ov{z}) \wedge \sqrt{2} dz = \oint_\infty 4 (\beta_{-1})_{\ou 1 1}
  (i \sqrt{2} z dt \wedge dz - \sqrt{2} z^2 d \ov{z} \wedge d z).
\end{equation}
By our choice of $\beta_{-1}$ we have that
$$
  - 2 (\beta_{-1})_{\ou 1} = \stackrel{\circ}{\Box}_b \beta_{-1}
  = - 4 \pi A \frac{\ov{z}}{\rho^6} \left\{ |z|^2 + i t \right\}
  + O(\rho^{-4}),
$$
and moreover
$$
   (\pa_z + i \ov{z} \pa_t) \left[ \frac{\ov{z}}{\rho^6} \left( |z|^2
   + i t \right) \right] = - 3 \ov{z}^2 (|z|^2 + i t)^2 \rho^{-6},
$$
so we find
$$
  4 (\beta_{-1})_{\ou 1 1} = - 12 \sqrt{2} A \pi \ov{z}^2
  \frac{(|z|^2 + i t)^2}{\rho^{10}}.
$$
Inserting this expression into the boundary formula we get
\begin{eqnarray*}
  & & 4 \oint_\infty (\beta_{-1})_{\ou 1 1} (i \sqrt{2} z dt \wedge dz - \sqrt{2} z^2
  d \ov{z} \wedge d z) \\ & = & - 12 \sqrt{2} \oint_\infty \pi A \left[
   |z|^4 -  t^2 + 2 i t |z|^2 \right] \rho^{-10} \ov{z}^2 (i \sqrt{2} z dt \wedge dz - \sqrt{2}
 z^2 d  \ov{z} \wedge d z) \\ & = & - 24 \pi A \oint_S
   \left[ |z|^4 -  t^2 + 2 i t |z|^2 \right] \ov{z}^2 (i  z dt
 \wedge dz - z^2 d \ov{z} \wedge d z).
\end{eqnarray*}
Using the fact that $r^4 + t^2 = 1$ on $S$ and the relations
$$
  dz \wedge d \ov{z} = - 2 i r dr \wedge d \var; \qquad \quad
  dt \wedge dz = i r e^{i \var} dt \wedge d \var,
$$
one finds
\begin{eqnarray*}
  & & 4 \oint_\infty (\beta_{-1})_{\ou 1 1} (i \sqrt{2} z dt \wedge dz - \sqrt{2} z^2
  d \ov{z} \wedge d z) \\ & = & - 2 4 \pi A \int_S \left[ |z|^4 - t^2 + 2 i t |z|^2
  \right] r^2 \left( - r^2 dt \wedge d \var - 2 i r^3 dr \wedge d \var \right)
  \\ & = & - 24 \pi A \int_S \left[ |z|^4 - t^2 + 2 i t |z|^2
  \right] \left( r^4 d \var \wedge dt + 2 i r^5 d \var \wedge dr \right)
  \\ & = & - 24 \pi A \int_S \left[ |z|^4 - t^2 + 2 i t |z|^2
  \right] \left( r^4 d \var \wedge dt - i r^2 t d \var \wedge dt \right).
\end{eqnarray*}
The terms which are odd in $t$ vanish after integration, so the last
expression becomes
$$
   - 48 \pi^2 A \int_{-1}^1 (1-t^2) dt =  - 64 \pi^2 A
  = -  \frac{4}{3} m(J,\th).
$$
In conclusion, from \eqref{eq:1318}, \eqref{eq:miao1}, \eqref{eq:miao2}
and the last formula we find
\begin{eqnarray*}
  \frac{2}{3} m(J,\th) & = & -\int_N|\Box_b(\beta)|^{2}{\theta}{\wedge}d{\theta}
  + 2 \int_N |\beta_{, \ou \ou}|^2 \theta \wedge d \theta
  + \frac 12 \int_N \ov{\beta} P \beta \theta \wedge d \theta
  \\ & = & -\int_N|\Box_b(\beta)|^{2}{\theta}{\wedge}d{\theta}+
  2 \int_N |\beta_{, \ou \ou}|^2 \theta \wedge d \theta \\
   & + & \frac 12 \int_N (\ov{\beta}_{-1} + \ov{\beta}_{-2+\e} + u)
  P (\beta_{-1} + \beta_{-2+\e} + \ov{u}) \, \theta \wedge d \theta
  - \frac{2}{3} m(J,\th).
\end{eqnarray*}
Viewed as a function on $M$, $\beta_{-1} + \beta_{-2+\e} + \ov{u}$,
vanishes at $p$, is of class $\mathfrak{S}^{2,2}(M)$, and is smooth on $M
\setminus \{p\}$. Therefore, by the conformal invariance of $P$ and
its non-negativity we have
\begin{eqnarray*}
   & & \int_N (\ov{\beta}_{-1} + \ov{\beta}_{-2+\e} + u) P_\theta (\beta_{-1} +
  \beta_{-2+\e} + \ov{u}) \, \theta \wedge d \theta \\
   & = & \int_{M \setminus
  \{p\}} (\ov{\beta}_{-1} + \ov{\beta}_{-2+\e} + u) P_{\hat{\th}}
  (\beta_{-1} + \beta_{-2+\e} + \ov{u}) \, \hat{\theta} \wedge d
  \hat{\theta} \geq 0.
\end{eqnarray*}
Notice that the non-negativity of $P$ is assumed on smooth functions, but by approximation and integration 
by parts one also gets a sign condition on the integral of $v  P_{\hat{\th}} \ov{v}$ for $v$ of 
class $\mathfrak{S}^{2,2}(M)$. 

Summing up we then obtain
\begin{eqnarray*}
  \frac{4}{3} m(J,\th) & = & -\int_N|\Box_b(\beta)|^{2}{\theta}{\wedge}d{\theta}+
  2 \int_N |\beta_{, \ou \ou}|^2 \theta \wedge d \theta + 2 \int_N R
  |\beta_{, \ou}|^2  \th \wedge d \th \\
  & + & \frac 12 \int_{M \setminus \{p\}} (\ov{\beta}_{-1} + \ov{\beta}_{-2+\e} + u)
  P_{\hat{\th}} (\beta_{-1} + \beta_{-2+\e} + \ov{u}) \, \hat{\theta}
  \wedge d \hat{\theta},
\end{eqnarray*}
which implies \eqref{3.9.2}. Formula \eqref{3.9.3} then follows immediately. \end{pf}

\subsection{Conclusion of the proof of Theorem \ref{t:pm}}\label{ss:concl}

The first statement of the theorem follows immediately from
Proposition \ref{p:P>0}.

Suppose now that $m(J,\th) = 0$. Then, from \eqref{3.9.3}, Proposition \ref{p:exbeta}
and \eqref{3.9.1} we have the following identities
\begin{equation}\label{3.5.1}
    \beta_{, \ou \ou} \equiv 0; \quad \qquad \beta_{, \ou 1} \equiv 0;
  \quad \qquad P \beta \equiv 0.
\end{equation}
The first two relations imply that $|\beta_{, \ou}|^2$ is constant: in particular,
from the behavior of $\beta_{, \ou}$ at infinity we deduce that $|\beta_{, \ou}|^2 \equiv
\frac 12$. We also have then
$$
  R \equiv 0.
$$
From $P \beta = 0$ we also deduce that $A_{11,\ou} \equiv 0$. Let us now show that the
torsion vanishes identically: in order to do this, we follow the arguments of Section
3 in \cite{SY1}. Since the ideas here require rather straightforward modifications,
we will be rather sketchy in this part.

We consider the flow $\var_s$ generated by the Reeb vector field $T$ of $N$, and we let
$$
J_{(s)} = \var_s^* J  \qquad \quad (\dot{J} = L_T J = 2 A_{J,\th}),
$$
where $A_{J,\th} = - i A_{11} \th^1 \otimes Z_{\ou} + $conj.. By the variation
formulas \eqref{eq:varoA}, \eqref{eq:varW} we have that
\begin{equation}\label{eq:RRAA}
    \frac{d}{ds} R_{J_{(s)},\th} = - 2 |A_{11}|^2 + i (A_{11,\ou \ou} -
  A_{\ou \ou, 11}); \qquad \quad \frac{d}{ds} (A_{11})_{J_{(s)},\th} = i A_{11,0}.
\end{equation}
Since for an asymptotically flat manifold we have $A_{11} = O(\rho^{-4})$ at
infinity (see Remark \ref{r:regtor}) and since $A_{11,\ou} \equiv 0$ at
$s = 0$, from the second equation in \eqref{eq:RRAA} we derive that
$$
    A_{11,\ou \ou} = O(s \rho^{-8}) \quad \hbox{ for } s \hbox{ small and }
    \rho \hbox{ large}.
$$
Then the first equation in \eqref{eq:RRAA} implies that for $s$ small
\begin{equation}\label{eq:estRs}
    R_{J_{(s)},\th} \geq \left\{
                          \begin{array}{ll}
                            - C s & \hbox{ everywhere on $N$};  \\
                            - C \frac{s}{\rho^8} & \hbox{ at infinity}.
                          \end{array}
                        \right.
\end{equation}
In particular, it turns out that the negative part of $R_{J_{(s)},\th}$
is arbitrarily small as $s \to 0$. Similarly to Lemma 3.1 in \cite{SY1}
one has the Sobolev type inequality
$$
  \left( \int_N u^4 \th \wedge d \th \right)^{\frac 12} \leq C \int_N
  |\n_b u|^2 \th \wedge d \th; \qquad C > 0, u \in C^\infty_c(N).
$$
This implies that
\begin{eqnarray*}
  \int_N \left( 4 |\n_b u|^2 + R_{J_{(s)},\th} u^2 \right) \th \wedge d \th &
  \geq & 4 \int_N  |\n_b u|^2 \th \wedge d \th - \left( \int_N R_{J_{(s)},\th}^2
  \th \wedge d \th \right)^{\frac 12} \left( \int_N u^4 \th \wedge d \th
   \right)^{\frac 12} \\
   & \geq & (4 + o_s(1)) \int_N  |\n_b u|^2 \th \wedge d \th \qquad \hbox{ as }
   s \to 0.
\end{eqnarray*}
Therefore, the operator $- 4 \D_b + R_{J_{(s)},\th}$ is coercive for
$s$ small. Estimates similar to \eqref{eq:estRs} imply that $|R_{J_{(s)},\th}|
\leq C \rho^{-8}$ at infinity,  and hence we can find a solution $v_s$ of
$$
  - 4 \D_b v_s + R_{J_{(s)},\th} v_s = R_{J_{(s)},\th} \qquad \quad \hbox{ on } N,
$$
which decays to zero at infinity. More precisely, using estimates similar to
those of Lemma 3.2 in \cite{SY1}, together with Lemma \ref{l:greenheis}
one finds that
$$
  v_s = \frac{1}{32 \pi \rho^2} \int_N (R_{J_{(s)},\th} - R_{J_{(s)},\th}
   v_s) \th \wedge d \th + O(\rho^{-3}) \qquad \quad \hbox{ at infinity}.
$$
The function $u_s = 1 - v_s$ then satisfies
\begin{equation}\label{eq:us}
    - 4 \D_b u_s + R_{J_{(s)},\th} u_s = 0 \qquad \quad \hbox{ on } N,
\end{equation}
which means that the manifold $(N, J_{(s)}, u_s^2 \th)$ has zero Tanaka-Webster
curvature and is asymptotically flat. By the asymptotic formula for $v_s$
we deduce immediately
$$
  u_s = 1 - \frac{1}{32 \pi \rho^2} \int_N R_{J_{(s)},\th} u_s \th \wedge d \th
+ O(\rho^{-3}) \qquad \quad \hbox{ at infinity},
$$
and hence by Lemma \ref{l:m=A} (see also the coefficient of $\th$ in Definition
\ref{d:aflat}) we have that
\begin{equation}\label{eq:massss}
    m(J_{(s)}, u_s^2 \th) = - \frac 34 \int_N
  R_{J(s),\th} u_s \theta \wedge d \theta.
\end{equation}
Differentiating the last formula in $s$ and using the fact that $u_0 \equiv 1$,
together with \eqref{eq:RRAA},  we find
$$
  \frac{d}{ds}|_{s=0} m(J_{(s)}, u_s^2 \th) = \frac 32 \int_N |A_{11}|^2 \theta \wedge
  d \theta > 0 \qquad \hbox{ in case } A_{11} \not \equiv 0 \hbox{ at } s = 0.
$$
This implies that for $s$ negative we obtain a manifold which is asymptotically
flat, with zero Tanaka-Webster curvature and with negative mass.

We show that this is a contradiction to the first part of the theorem. In fact,
the flow $\varphi_s$ on the  manifold $N$ corresponds to the one parameter
family of diffeomorphisms $\hat{\varphi}_{(s)}$ on $(M, \hat{J}, \hat{\th})$ generated
by the vector field $(G_p)^{-2} \left[ \hat{T} + 2i (\log G_p)_{, 1} \hat{Z}_{\ou}
- 2i (\log G_p)_{, \ou} Z_1 \right]$, where $G_p$ is the Green's function of $L_b$ on $M$
for $s = 0$, and $\hat{T}$ the corresponding Reeb vector field. Letting $\hat{J}_{(s)}$
denote $(\hat{\varphi}_{(s)})^* \hat{J}$, then one has
clearly that $(G_p)_{(s)} = u_s G_p$ (since $u_s G_p \hat{\th}$ has zero Tanaka-Webster
curvature, and the correct asymptotics near the singularity). Hence we obtain
$(N, J_{(s)}, u_s^2 \th)$ as blow-up of the compact manifold $(M, \hat{J}_{(s)},
\hat{\th})$, through the Green's function $(G_p)_{(s)}$.

The contradiction will follow if we show that $(M, \hat{J}_{(s)}, \hat{\th})$
satisfies the same assumptions of Theorem \ref{t:pm}. The Tanaka-Webster class of
$(M, \hat{J}_{(s)})$ is indeed clearly positive, since it is a perturbation of a
structure with positive class. Concerning the Paneitz operator, notice that we
are pulling back the structure via the contact diffeomorphisms $\hat{\varphi}_{(s)}$,
so the positivity of $P_{\hat{J}_{(s)},\hat{\th}}$ follows from a change of
variable and its conformal invariance. Therefore, we proved that $m(J,\th) = 0$ also
implies $A_{11} \equiv 0$.

\

\noindent To show that $N$ coincides with Heisenberg, first we define a map
from a neighborhood of infinity $\mathcal{U}$ in $N$ to a neighborhood
of infinity $\mathcal{V}$ in $\H^1$.
From \eqref{3.5.1} we find that
$$
 d(\beta_{, \ou} \th^{\ou}) = \beta_{, \ou 1} \th^1 \wedge
 \th^{\ou}  + \beta_{, \ou 0} \th \wedge \th^{\ou} = 0,
$$
which implies that  $\beta_{,1{\bar 1}}=\beta_{,11{\bar 1}}=0$,  $\beta_{,11}=0$, $|\beta_{,1}|^2=constant=0$, 
and hence $d \beta = \beta_{, \ou} \th^{\ou}$. Taking $z = \ov{\beta}$ we have
$$
  d (\th - i z d \ov{z} + i \ov{z} d z) = i \th^1 \wedge \th^{\ou}
  - 2 i dz \wedge d \ov{z} = 0.
$$
Note that $dz \wedge d \ov{z} = |\beta_{, \ou}|^2 \th^1 \wedge \th^{\ou}$ and
$|\beta_{,\ou}|^2 = \frac 12$ hence, if we take $\mathcal{U}$ simply connected
(for example, $\{ \rho \geq R \}$, $R \gg 1$, in inverted CR coordinates), there
exists a function $\tilde{t}$ such that
$$
  d\tilde{t} = \th - i z d \ov{z} + i \ov{z} dz.
$$
So we get a pseudohermitian isomorphism between $\mathcal{U}$ and its image in $\H_1$,
$\mathcal{V}$, if we send $q \in N$ into
$$
  q \mapsto (z(q),t(q)) = \left(\ov{\beta}(q),
   \int_{q_0}^q  d\tilde{t} \right),
$$
where we are taking curves connecting $q_0$ to $q$ inside $\mathcal{U}$.

We call $\Psi : \mathcal{V} \to \mathcal{U}$ (sets which we can assume to be
connected by arcs) the inverse of this map: next, we want to extend $\Psi$
globally on $\H^1$. Taking $q_0 \in \mathcal{V}$ and $q \in
\H^1$ arbitrary, we can find a curve
$$
  \Gamma : [0,1] \to \H^1, \qquad \quad \G(0) = q_0, \quad \G(1) = q.
$$
We can then consider the curve $\tilde{\G} : [0,1] \to N$ obtained composing $\G$ with $\Psi$
initially, and then prolonging it through the local pseudohermitian isomorphisms obtained before.
In fact, given any such flat structure $(M,J,\theta)$, any point $p$ in M, and any
unit contact vector $e_1$ at $p$, there exists a unique local pseudohermitian isomorphism
mapping a neighborhood of the origin in $\H^1$ carrying $0$ to $p$, and $\partial_x$ to $e_1$
(hence $\partial_y$ to $e_2= J e_1$, and $\partial_t$ to $T$) onto a neighborhood
of $p$ in $M$. Then using the group structure, the same is true if we substitute
$0$ above by any point $q$ in the Heisenberg group, and $\partial_x$, by any unit
contact vector at $q$.
`

We show that this procedure defines a map $\tilde{\Psi} : \H^1 \to N$, showing that
$\tilde{\G}(1)$ is independent of the choice of $\G$ (once $q$ is fixed).
In fact, if we choose two different curves $\G_0, \G_1$, we can then consider the new
(closed) curve $\G_2 = \G_1^{-1} * \G_0 : [0,1] \to \H^1$, with $\G_2(0) =
\G_2(1) = q_0$. We are done if we show that $\tilde{\G}_2(1) = \tilde{\G}_2(0)$.

By the simple connectedness of $\H^1$, there is a homotopy $\ov{\G} : [0,1] \times [0,1]
\to \H^1$ such that $\ov{\G}(\cdot,0) \equiv q_0$, $\ov{\G}(t,1) = \G_2(t)$ and such
that $\ov{\G}(0,s) = \ov{\G}(1,s) = q_0$ for all $s$. We call $\tilde{\ov{\G}}(\cdot,s)$
the curve constructed with the above procedure starting from $\ov{\G}(\cdot, s)$.
To prove the claim, it is sufficient to show that
$$
\tilde{\ov{\G}}(1,s) = \tilde{\G}(0) = \Psi(q_0) \qquad \hbox{ for all } s \in [0,1].
$$
This is clearly true for $s$ close to $0$, since we have a local pseudohermitian isomorphism between a
neighborhood of $q_0$ and a neighborhood of $\Psi(q_0)$, and since
$\tilde{\ov{\G}}(\cdot,s)$ is contained in the first neighborhood for $s$ small.
We call
$$
   I = \left\{ r \in (0, 1] \; : \; \tilde{\ov{\G}}(1,s) = \tilde{\G}(0)
   \hbox{ for every } s \in [0,r] \right\}; \qquad \quad \ov{s} = \sup I.
$$
Our goal is equivalent to prove that $\ov{s} = 1$. It is clear that $\ov{s}
\in I$, so we obtain the claim if we show that $I$ is open in $[0,1]$.
Assume by contradiction that $\ov{s} < 1$: we can find finitely many $t_i$'s
$$
  0 = t_0 < t_1 < \cdots < t_{k-1} < t_k = 1
$$
such that in each interval $[t_{i-1}, t_{i+1}]$ a neighborhood $\mathcal{V}_i$
of $\ov{\G}(\cdot,\ov{s})$ is mapped homeomorphically onto a neighborhood
$\mathcal{U}_i$ of $\tilde{\ov{\G}}(\cdot,\ov{s})$, for $i = 1, \dots k-1$.
We can then find $\e$ so small that
$$
  \ov{\G}(t,s) \in \mathcal{V}_i \qquad \hbox{ for } t \in [t_{i-1},
  t_{i+1}], s \in [\ov{s}, \ov{s}+\e] \hbox{ and } i = 1, \dots, k-1.
$$
Since $\Psi(q_0) \in \mathcal{U}_{k-1}$ and since all curves $\ov{\G}(\cdot,s)
\in \mathcal{V}_{k-1}$ for $t \in [t_{k-1}, 1]$ and $s \in [\ov{s}, \ov{s}+\e]$,
all ending in $q_0$, we have also that $\tilde{\ov{\G}}(1,s) = \Psi(q_0)$
for $s \in [\ov{s}, \ov{s}+\e]$, so we obtain the desired claim, namely the
well definiteness of $\tilde{\Psi}$ as a local pseudohermitian isomorphism, and therefore
as a covering map.

If $\mathcal{V}$ is as above, and if $q \in \mathcal{V}$, then we can connect
$q_0$ to $q$ with a curve $\Gamma$ inside $\mathcal{V}$. This implies that
$\tilde{\Psi}$ coincides with $\Psi$ on $\mathcal{V}$, and hence extends it.
By construction, $\tilde{\Psi}$ is also proper, and hence by Corollary 1.4
in \cite{AP}, Chapter 3, it is also injective. It is indeed easy to prove
its surjectivity: given any point $\tilde{p} \in N$, we can find a curve
$\tilde{\G} : [0,1] \to N$ such that $\tilde{\G}(0) = \Psi(q_0)$ and
such that $\tilde{\G}(1) = \tilde{p}$. Then, if we consider the
curve $\Gamma : [0,1] \to \H^1$ inverting $\tilde{\Psi}$ along
$\tilde{\G}$ (see Lemma 1.11 in Chapter 3 of \cite{AP}), we have that
$\tilde{\Psi}(\G(1)) = \tilde{\G}(1) = \tilde{p}$. Note that $\H^1 \setminus
\mathcal{V}$ is compact. In conclusion, we have shown that $\tilde{\Psi}$ realizes
a pseudohermitian isomorphism between $\H^1$ and $N$. This concludes the proof.

\

\begin{rem} If $J$ is  spherical we have a different proof of
the vanishing of $A_{11}$. In fact, if $\mathfrak{Q}$ stands for the Cartan tensor
(see Subsection \ref{ss:CRcoord}), for  spherical structures one has
$$
  0 = \mathfrak{Q}_{11} = \frac 16 R_{,11} + \frac i2 R A_{11} - A_{11,0}
  - \frac 23 i A_{11,\ou 1}.
$$
Then (see Section \ref{s:app}), since $A_{11,\ou} = 0$ and $R = 0$ we find that
$$
  0 = i A_{11,0} = A_{11,1\ou} - A_{11,\ou 1} - 2 A_{11} R,
$$
which implies $A_{11,1\ou} = 0$. By integration we obtain
$$
 0 = - \int_N A_{11,1\ou} A_{\ou \ou} \theta \wedge d \theta = \int_N
 A_{11,1} A_{\ou \ou, \ou} \theta \wedge d \theta = \int_N
 |A_{11,1}|^2 \theta \wedge d \theta,
$$
which implies $A_{11,1} = 0$. Differentiating we get $|A_{11}|^2 = const.$
and the constant is zero at infinity, which implies $A_{11} \equiv 0$.
\end{rem}

\section{Some examples}\label{s:ex}

In this section we describe some examples: first we find a structure close to
the spherical one with nonpositive Paneitz operator. Then, we derive an example
of manifold with positive Tanaka-Webster class but with negative mass. We finally 
describe a CR structure on $S^2 \times S^1$ with non-negative Paneitz operator.

\subsection{Three dimensional CR manifolds with nonpositive Paneitz
operator}\label{ss:pneg}

Recall that from the commutation rules \eqref{eq:comm} one has
\begin{equation}
\Box _{b}u= - \Delta _{b}u+iTu=-(u_{,1\overline{1}}+u_{\overline{,1}%
1})+iu_{,0}=-2u_{,\overline{1}1}=-2Z_{1}(Z_{\overline{1}}u)-2\o _{\overline{1%
}}^{\overline{1}}(Z_{1})Z_{\overline{1}}u.  \label{4.1.1}
\end{equation}%
Recall also that, along a deformation $J_{(s)}$ of the CR structure
the following relations hold (see Subsection \ref{sss:varJ})
\begin{equation}
\dot{Z}_{1}=-iE_{11}Z_{\overline{1}};\qquad \qquad \dot{\theta}^{1}=-iE_{\bar{%
1}\bar{1}}\theta ^{\overline{1}},  \label{2.1.2.2}
\end{equation}%
\begin{equation}
\dot{\omega}_{1}^{1}=i\left( A_{11}E_{\overline{1}\overline{1}}+A_{\overline{%
1}\overline{1}}E_{11}\right) \theta -i(E_{11,\bar{1}}\theta ^{1}+E_{\bar{1}%
\bar{1},1}\theta ^{\overline{1}});\qquad \quad \dot{A}_{11}=iE_{11,0}.
\label{2.1.2.3}
\end{equation}%
\begin{equation}
- \dot{\Delta}_{b}=2i\left( E_{11}Z_{\overline{1}}Z_{\overline{1}}+E_{11,%
\overline{1}}Z_{\overline{1}}\right) +\hbox{ conj.}.  \label{2.1.2.6}
\end{equation}%
We derive next the first and second variations of the Paneitz operator near
the standard three dimensional pseudohermitian sphere. We can
express the Paneitz operator $P$ as follows
\begin{equation}
P\psi :=\overline{\Box }_{b}\Box _{b}\psi +4i(A_{11}\psi _{,\bar{1}})_{, \bar{1}}  \label{6.1.1}
\end{equation}
for a smooth real function $\psi .$ By (\ref{4.1.1}) the
expression (\ref{6.1.1}) is equivalent to
\begin{equation}
P\psi = 4(\psi _{,\bar{1}11}+iA_{11}\psi _{,\bar{1}})_{,\bar{1}}\qquad \text{ or}
\label{6.1.2} \qquad P\psi =
\Delta _{b}^{2}\psi +T^{2}\psi +4\func{\hbox{Im }}(A_{\bar{1}\bar{1}}\psi
_{,1})_{,1}.
\end{equation}
Along the above deformation $J_{(s)}$ , let $P_{(s)}$ denote the Paneitz operator
associated to $(J_{(s)},$ $\theta ).$ Let $(\cdot ,\cdot )$ denote the $L^{2}$ inner
product with respect to the volume form $\theta \wedge d\theta .$ We then find
\begin{equation}
\frac{d}{ds}(P_{(s)}\psi ,\psi )  \label{6.1.3}
 = 2(\dot{\Delta}_{b}\psi ,\Delta _{b}\psi )  -
4\func{\hbox{Im }}\int_{M}\dot{A}_{\bar{
1}\bar{1}}(\psi _{,1})^{2}\th \wedge d \th -
8\func{\hbox{Im }}\int_{M}A_{\bar{1}\bar{1}}(\dot{\psi}_{,1})\psi _{,1} \th \wedge d \th
\end{equation}
by (\ref{6.1.2}). From (\ref{2.1.2.2}) we have
\begin{equation}
(\dot{\psi}_{,1})=-iE_{11}\psi _{,\bar{1}}.  \label{6.1.4}
\end{equation}
Substituting (\ref{2.1.2.3}), (\ref{2.1.2.6}), and (\ref{6.1.4})
into (\ref{6.1.3}), we obtain
\begin{eqnarray}
\frac{d}{ds}(P_{(s)}\psi ,\psi )  \label{6.1.5}
&=&4\int_{M}E_{\bar{1}\bar{1}}\psi _{,1} \left( -i\psi _{,1\bar{1}1}-i\psi _{,\bar{
1}11}-\psi _{,10}+A_{11}\psi _{,\bar{1}}\right) \th \wedge d \th +\hbox{ conj.}   \\
&=&-8i\int_{M}E_{\bar{1}\bar{1}}\psi _{,1}(\psi _{,\bar{1}11}+iA_{11}\psi
_{,\bar{1}}) \th \wedge d \th +\hbox{ conj.}  \nonumber
\end{eqnarray}
using an integration by parts and the commutation relations (\ref{eq:comm}).
Take $\psi $ $\in $ $KerP_{(0)},$ i.e.
\begin{equation}
4(\psi _{,\bar{1}11}+iA_{11}\psi _{,\bar{1}})_{,\bar{1}}=0.  \label{6.1.6}
\end{equation}%
\noindent If $J_{(0)}$ has transverse symmetry, then we can find a contact
form $\theta _{(0)}$ such that $A_{11}$ $=$ $0$ with respect to $(J_{(0)},$ $
\theta _{(0)})$ (see \cite{L2}). By the transformation law \eqref{eq:transfP} we have
that $KerP_{(0)}$ and $(P_{(s)}\psi ,\psi )$ are independent of the choice
of contact form. So we may just take $\theta _{(0)}$ as the contact
form in (\ref{6.1.6}) to get $\psi _{,\bar{1}11\bar{1}}$ $=$ $0.$
Multiplying this by $\psi _{,1\bar{1}}$ and integrating, we get
\begin{eqnarray}
0 &=&\int_M \psi _{,\bar{1}11\bar{1}}\psi _{,1\bar{1}}\theta _{(0)}\wedge
d\theta _{(0)}  \label{6.1.7} = -\int_M \psi _{,\bar{1}11}
\psi_{,1\bar{1}\bar{1}}\theta _{(0)}\wedge
d\theta _{(0)}  \nonumber \\
&=&-\int_M |\psi _{,\bar{1}11}|^{2}\theta _{(0)}\wedge d\theta _{(0)}
\nonumber
\end{eqnarray}
using an integration by parts and the fact that $\psi$ is real. It follows that
$\psi _{,\bar{1}11}$ $=$ $0$, so from (\ref{6.1.5}) we have
\begin{equation}
\frac{d}{ds}|_{s=0}(P_{(s)}\psi ,\psi )=0  \label{6.1.8}
\end{equation}
for $\psi $ $\in $ $KerP_{(0)}$ if $J_{(0)}$ has transverse
symmetry. For example we can take $J_{(0)}$ $=$ $\hat{J},$ the standard pseudohermitian structure
compatible with $\hat{\xi}$ on $M$ $=$ $S^{3},$ where $(S^{3},$ $\hat{\xi})$
is the standard contact $3$-sphere. Note that from the above argument, $\psi
$ $\in $ $KerP_{(0)}$ implies that $\psi $ is CR-pluriharmonic if $J_{(0)}$
has transverse symmetry. So by Proposition 3.4 in \cite{L2}, we have
\begin{equation}
\psi _{,\bar{1}11}+iA_{11}\psi _{,\bar{1}}=0  \label{6.1.9}
\end{equation}
with respect to $J_{(0)}$ coupled with any contact form.
Starting with (\ref{6.1.5}) we compute%
\begin{equation}
\frac{d^{2}}{ds^{2}}(P_{(s)}\psi ,\psi )|_{s=0}  \label{6.1.10}
=-8i\int_{S^3}E_{\bar{1}\bar{1}}\psi _{,1}\frac{d}{ds}\left(\psi _{,\bar{1}
11}+iA_{11}\psi _{,\bar{1}}\right)|_{s=0} \theta _{(0)}\wedge
d\theta _{(0)} +\hbox{ conj.}
\end{equation}
by (\ref{6.1.9}). From (\ref{2.1.2.2}) and (\ref{2.1.2.3}) we obtain
\begin{equation}
\frac{d}{ds}\psi _{,\bar{1}1}=i[(E_{\bar{1}\bar{1}}\psi
_{,1})_{,1}-(E_{11}\psi _{,\bar{1}})_{,\bar{1}}].  \label{6.1.11}
\end{equation}
We then have%
\begin{equation}
\frac{d}{ds}\psi _{,\bar{1}11} = \left(\frac{d}{ds}Z_{1}\right) \left(\psi _{,\bar{1}
1}\right)+Z_{1} \left(\frac{d}{ds}\psi _{,\bar{1}1}\right)  \label{6.1.12} =
-iE_{11}\psi_{,\bar{1}1\bar{1}}+i(E_{\bar{1}\bar{1}}\psi _{,1})_{,11}
-i(E_{11}\psi _{,\bar{1}})_{,\bar{1}1}
\end{equation}
by (\ref{2.1.2.2}) and (\ref{6.1.11}). Substituting (\ref{6.1.12}),
the conjugate of (\ref{6.1.4}), and $\dot{A}_{11}$ $=$ $iE_{11,0}$ into
(\ref{6.1.10}) we get
\begin{eqnarray}
\frac{d^{2}}{ds^{2}}(P_{(s)}\psi ,\psi )|_{s=0}  \label{6.1.13}
&=&8\int_{S^3}E_{\bar{1}\bar{1}}\psi _{,1}\left[(E_{\bar{1}\bar{1}}\psi
_{,1})_{,11}-(E_{11}\psi _{,\bar{1}})_{,\bar{1}1} \right.  \\
&& \left.-E_{11}\psi _{,\bar{1}1\bar{1}}+iE_{11,0}\psi _{,\bar{1}}+iA_{11}E_{\bar{1}
\bar{1}}\psi _{,1}\right] \theta _{(0)}\wedge
d\theta _{(0)}+\hbox{ conj.}.  \nonumber
\end{eqnarray}
In view of the conjugate of (\ref{6.1.9}) and the commutation
relations (\ref{eq:comm}), we compute%
\begin{equation}
\psi _{,\bar{1}1\bar{1}} = \psi _{,1\bar{1}\bar{1}}-i\psi _{,0\bar{1}}
\label{6.1.14}
=iA_{\bar{1}\bar{1}}\psi _{,1}-i\psi _{,\bar{1}0}-i\psi _{,1}A_{\bar{1}
\bar{1}}
=-i\psi _{,\bar{1}0}.
\end{equation}
Substituting (\ref{6.1.14}) into (\ref{6.1.13}), integrating by
parts, and making use of the first commutation relation in (\ref{eq:comm}):
\begin{equation*}
i(E_{\bar{1}\bar{1}}\psi _{,1})_{,0} =(E_{\bar{1}\bar{1}}\psi _{,1})_{,1
\bar{1}} -(E_{\bar{1}\bar{1}}\psi _{,1})_{,\bar{1}1}+(E_{\bar{1}\bar{1}}\psi
_{,1})R,
\end{equation*}
hence we obtain%
\begin{eqnarray}
\frac{d^{2}}{ds^{2}}(P_{(s)}\psi ,\psi )|_{s=0}  \label{6.1.15}
&=&8\int_{S^3}\left(-[(E_{\bar{1}\bar{1}}\psi _{,1})_{,1}]^{2}+2|(E_{\bar{1}\bar{1}
}\psi _{,1})_{,1}|^{2}-|(E_{\bar{1}\bar{1}}\psi _{,1})_{,\bar{1}}|^{2} \right.
 \\
&& \left. -R|E_{\bar{1}\bar{1}}\psi _{,1}|^{2}+iA_{11}(E_{\bar{1}\bar{1}}\psi
_{,1})^{2}\right) \theta _{(0)}\wedge
d\theta _{(0)}+\hbox{ conj.}  \nonumber \\
&=&8\int_{S^3}\left(4|(E_{\bar{1}\bar{1}}\psi _{,1})_{,1}|^{2}-[(E_{\bar{1}\bar{1}
}\psi _{,1})_{,1}]^{2}-[(E_{11}\psi _{,\bar{1}})_{,\bar{1}}]^{2} \right. \nonumber
\\
&& \left. -2|(E_{\bar{1}\bar{1}}\psi _{,1})_{,\bar{1}}|^{2}-2R|E_{\bar{1}\bar{1}%
}\psi _{,1}|^{2}\right) \theta _{(0)}\wedge
d\theta _{(0)}.  \nonumber
\end{eqnarray}
In the last equality of \eqref{6.1.15} we have used $A_{11} = 0$. We now choose
$E_{11}$ $\neq $ $0$ such that
\begin{equation}
E_{11,\bar{1}}=0  \label{6.1.16}
\end{equation}
Assuming (\ref{6.1.16}), we reduce (\ref{6.1.15}) to%
\begin{eqnarray}
\frac{d^{2}}{ds^{2}}(P_{(s)}\psi ,\psi )|_{s=0}  \label{6.1.17}
&=&8\int_{S^3}\left( 4|E_{\bar{1}\bar{1}}\psi _{,11}|^{2}-[E_{\bar{1}\bar{1}}\psi
_{,11}]^{2}-[E_{11}\psi _{,\bar{1}\bar{1}}]^{2} \right.  \\
&& \left. -2|(E_{\bar{1}\bar{1}}\psi _{,1})_{,\bar{1}}|^{2}-2R|E_{\bar{1}\bar{1}
}\psi _{,1}|^{2}\right) \theta _{(0)}\wedge d\theta _{(0)}.\text{ }  \nonumber
\end{eqnarray}
Now on ($S^{3},$ $J_{(0)}$ $=$ $\hat{J},$ $\theta _{(0)}$ $=$ $
\hat{\theta})$, where
\begin{equation}\label{eq:hattheta}
    \hat{\theta} = i (\ov{\partial} - \partial) (|z_1|^2 + |z_2|^2)
\end{equation}
is the standard contact form on ($S^{3},$ $\hat{\xi}),$ we have that $R$
is a positive constant. Note also that $Z_{1}$ $=$ $\frac{1}{\sqrt{2}} \left( \bar{z}_{2}
\frac{\partial}{\partial z_{1}}-
\bar{z}_{1}\frac{\partial }{\partial z_{2}} \right),$ $Z_{\bar{1}}$ $=$ $
\frac{1}{\sqrt{2}} \left( z_{2}\frac{\partial }{\partial \bar{z}_{1}}-z_{1}
\frac{\partial }{\partial \bar{z}_{2}} \right).$ To make (\ref{6.1.17}) negative,
we choose $\psi $ $=$ $z_{1}+\bar{z}_{1}$ where $(z_{1},$
$z_{2})$ $\in $ $S^{3}$ $\subseteq $ $\C^{2}.$ It follows that%
\begin{equation}
\psi _{,1} = \frac{1}{\sqrt{2}} \bar{z}_{2}; \qquad \text{ }\psi _{,11}=0;
\qquad  \label{6.1.18}
\psi _{,1\bar{1}} = - \frac{1}{2} z_{1}; \qquad
\text{ }\psi _{,\bar{1} \bar{1}}=0.
\end{equation}
Observing that the first three terms of the integrand in (\ref{6.1.17})
vanish due to the fact that $\psi _{,11}$ $=$ $0$ by (\ref{6.1.18}), we have%
\begin{eqnarray}
\frac{d^{2}}{ds^{2}}(P_{(s)}\psi ,\psi )|_{s=0}  \label{6.1.19}
&=& 8 \int_{S^{3}}\left( -\left|E_{\bar{1}\bar{1},\bar{1}}\bar{z}_{2}-E_{\bar{1}\bar{1}
} \frac{1}{\sqrt{2}} z_{1}\right|^{2}- R|E_{\bar{1}\bar{1}}\bar{z}_{2}|^{2}\right)
\hat{\theta}\wedge d\hat{\theta}<0.
\end{eqnarray}
For example, we can take $E_{11}$ to be a nonzero constant in (\ref
{6.1.19}), which satisfies (\ref{6.1.16}). In view of (\ref{6.1.19}), (\ref
{6.1.8}), and $(P_{(0)}\psi ,\psi )$ $=$ $0,$ we can find $\varepsilon $ $>$
$0$ small so that
\[
(P_{(s)}\psi ,\psi )<0
\]
for $0$ $<$ $s$ $\leq $ $\varepsilon .$ In conclusion we obtain the following
result.

\begin{pro}\label{p:P<0} There exist CR structures
on $S^3$ arbitrarily close to the spherical one for which the CR Paneitz
operator is not non-negative.
\end{pro}

\begin{rem}\label{r:embed}
In \cite{ANSIU} and \cite{ROS} some examples of non-embeddable structures
were given, where $E_{11}$ was taken to be a nonzero constant (the associated deformed
CR structures are indicated by their type $(0,1)$ vector fields $Z_{\ou} + \frac{iE_{{\bar 1}{\bar 1}}s}{\sqrt{1+|E_{{\bar 1}{\bar 1}}s|^{2}}}Z_{1})$, satisfying
in particular \eqref{6.1.16}. Observe that the tangency condition
for the embeddability reads
$$
  E = B_J(f) := (f_{, 11} + i A_{11} f) \th^1 \otimes Z_{\ou} +
  (\overline{f}_{, \ou \ou} - i A_{\ou \ou} \overline{f}) \th^{\ou} \otimes Z_1,
$$
for some complex-valued function $f$ (\cite{CHENG}). Since on the sphere the
torsion vanishes, the condition simplifies as
\begin{equation}\label{eq:tang}
    E = f_{, 11} \th^1 \otimes Z_{\ou} + \overline{f}_{, \ou \ou} \th^{\ou} \otimes Z_1
\end{equation}
for some complex-valued function $f$.  For other embeddability criteria using
a Fourier representation (or a suitable normal form) or a spectrum of $\Box_b$, see \cite{BLA} and \cite{BE}.

We notice that an $E$ satisfying \eqref{6.1.16} cannot correspond to any nonzero
$f$ in \eqref{eq:tang}. In fact if this would be the case, we would have that
$f_{, 11 \ou} = 0$. Multiplying this equation by $\overline{f}_{, \ou}$ and integrating
by parts we would get $\int_M |f_{, 11}|^2 \theta \wedge d \theta = 0$,  which would also
imply $E_{11} = 0$. Therefore we expect our examples to be non-embeddable ones.
The same comments apply to the construction in the next subsection.
\end{rem}

\begin{rem} The statement of Proposition \ref{p:P<0} can also be deduced by applying
the results in \cite{BLA} and \cite{CHCY1}.
\end{rem}

\subsection{First and second variations of the mass near the standard
sphere}\label{ss:12mass}

In this subsection we perturb the pseudohermitian structure of the standard sphere
in order to understand the variation of the mass. We begin with a deformation
in $\H^1$, which we will eventually transfer to a compact setting.

\

\noindent We consider a deformation $J(s)$ of the standard pseudohermitian structure of $\H^1$,
namely such that $J(0)$ is as in \eqref{eq:J0}. Recall the variation formulas for $\D_b$
and $R$ in Subsection \ref{sss:varJ}. Assuming the deformation of $J$
decays sufficiently fast near infinity, we can find a scalar flat conformal
contact form $\th_s = u_s^2 \stackrel{\circ}{\th}$, with $u_s$ satisfying
\eqref{eq:us}, as in Subsection \ref{ss:concl}.  Then the mass $m(J_{(s)}, u_s^2
\th)$ will be still given by formula \eqref{eq:massss}. Differentiating this
formula twice in $s$, taking into account that $u \equiv 1$
and $R \equiv 0$ for $s = 0$ we obtain, at $s = 0$
\begin{equation}\label{eq:ddotm}
    \ddot{m}(J, \th) = - \frac 34 \int_{\H^1} \ddot{R} \stackrel{\circ}{\th}
    \wedge d \stackrel{\circ}{\th} - \frac 32 \int_{\H^1} \dot{R} \dot{u}
   \stackrel{\circ}{\th} \wedge d \stackrel{\circ}{\th}.
\end{equation}
We choose a deformation $J(s)$ so that $E_{11}$ (see the notation in Subsection
\ref{sss:varJ}) is a CR function, namely $E_{11, \ou} = 0$. We can take for example
$$
   E_{11}(z,\ov{z},t) = \frac{1}{\left( t + i (|z|^2 + 1) \right)^k},
$$
with $k$ integer and large (to have a fast decay). By \eqref{eq:varW} we have then
$\dot{R} = 0$, which by \eqref{eq:ddotm} implies
\begin{equation}\label{eq:masssss2}
    \ddot{m}(J, \th) = - \frac 34 \int_{\H^1} \ddot{R} \stackrel{\circ}{\th}
    \wedge d \stackrel{\circ}{\th}.
\end{equation}
Reasoning as for \eqref{eq:ddotD}, one can check that on $\H^1$
\begin{eqnarray}\label{eq:ddotR}\nonumber
  \ddot{R} & = & - 6 |E_{11,\ou}|^2 - 2 |E_{11,1}|^2 - 3 E_{11,\ou 1}
  E_{\ou \ou} - 3 E_{\ou \ou, 1 \ou} E_{11} - E_{11} E_{\ou \ou, \ou 1}
  - E_{\ou \ou} E_{11,1 \ou}  \\
   & - & i E_{11,0} E_{\ou \ou} + i E_{\ou \ou, 0} E_{11}. \nonumber
\end{eqnarray}
Therefore, by our choice of $E_{11}$ one has
$$
   \ddot{R} =  - 2 |E_{11,1}|^2 - E_{11} E_{\ou \ou, \ou 1}- E_{\ou \ou}
  E_{11,1 \ou} - i E_{11,0} E_{\ou \ou} + i E_{\ou \ou, 0} E_{11}.
$$
Using the commutation rules \eqref{eq:comm} and again the fact that
$E_{11}$ is a CR function we obtain
$$
   \ddot{R} =  - 2 |E_{11,1}|^2 - 2 E_{11} E_{\ou \ou, \ou 1}
  - 2 E_{\ou \ou} E_{11,1 \ou}.
$$
Then the above formula \eqref{eq:masssss2} implies
$$
  \ddot{m}(J,\th) = \frac 32 \int_{\H^1} \left( |E_{11,1}|^2 + E_{11}
  E_{\ou \ou, \ou 1} + E_{\ou \ou} E_{11,1 \ou} \right)
  \stackrel{\circ}{\th} \wedge d \stackrel{\circ}{\th},
$$
so integrating by parts we obtain
$$
  \ddot{m}(J,\th) = - \frac 32 \int_{\H^1} |E_{11,1}|^2
  \stackrel{\circ}{\th} \wedge d \stackrel{\circ}{\th} < 0.
$$
We can transport the latter example on $S^3$ using the Cayley transform 
$\varpi : S^3 \setminus p \to \H^1$, $p = (0,1) \in \C^2$, defined as
$$
  \varpi(z_1,z_2) = \left( \frac{z_1}{1+z_2}, \hbox{Re } \left(
  i \frac{1-z_2}{1+z_2} \right) \right),
$$
where $(z_1,z_2)$ are standard coordinates in $\C^2$. It turns out that
$$
  G_p(z_1,z_2) = \frac{1}{\pi} \left( \frac{(1+z_2)
  (1+\ov{z}_2)}{|z_1|^4 - (z_2 - \ov{z}_2)^2} \right)^{\frac 12}.
$$
As in Subsection \ref{ss:concl}, we can obtain $(\H^1, J_{(s)}, u_s^2 \th)$
as blow up of the structure $(S^3, \hat{J}_{(s)}, \hat{\th})$ at the point
$p$ through the Green's function $(G_p)_s =(u_s \circ \varpi) G_p$.
In conclusion, we derived the following result.

\begin{pro}\label{p:example} There exist compact three dimensional pseudohermitian manifolds
of positive Tanaka-Webster class such that, once the contact form is blown-up as in
\eqref{eq:bubu} at proper points, the resulting asymptotically flat pseudohermitian
manifold has negative mass.
\end{pro}

\begin{rem}\label{r:1stvarm} For an arbitrary variation $E_{11}$ we have
$$
   \dot{R} = i\left( E_{11,\overline{1}\overline{1}}-E_{\overline{1}\overline{1}%
,11}\right),
$$
and
$$
  \dot{m}(J, \th) = - \frac 34 \int_{\H^1} \dot{R} \stackrel{\circ}{\th}
    \wedge d \stackrel{\circ}{\th} = - \frac 34 i \int_{\H^1}
   \left( E_{11,\overline{1}\overline{1}}-E_{\overline{1}
  \overline{1},11}\right) \stackrel{\circ}{\th}
    \wedge d \stackrel{\circ}{\th}.
$$
Assuming sufficiently fast decay to zero of $E_{11}$ and its derivatives,
integrating by parts we obtain $\dot{m}(J, \th) = 0$, namely vanishing
of the first variation of the p-mass under arbitrary perturbations (with
sufficient decay and regularity).
\end{rem}

\subsection{A CR structure on $S^2 \times S^1$ with non-negative Paneitz operator and non-vanishing torsion} \label{ss:s2s1}

Consider the Heisenberg group $\H^1$ with it origin removed:  endowing it with the standard 
CR structure $J_0$ and with the  contact form 
$$
  \check{\th} = \rho^{-2} \stackrel{\circ}{\th}, 
$$
the transformation $(z,t) \mapsto (2z, 4t)$ defines a pseudohermitian isomorphism.  Introducing then the equivalence relation 
$$
  (z_1, t_1) \sim (z_2, t_2)   \qquad \quad \hbox{ if and only if } \qquad \quad  z_1 = 2^k z_2 
  \hbox{ and } t_1 = 2^{2k} t_2 \hbox{ for some } k \in \Z,
$$
we obtain a pseudohermitian structure on $S^2 \times S^1$. Using \eqref{eq:varAconf}, one finds that 
$$
  \check{A}_{11}  = i |z|^2 \left(|z|^2 + i t \right)^2 \rho^{-6}, 
$$
so the torsion of $(S^2 \times S^1, J_0, \check{\th})$  is non-zero. Given any smooth 
function $\varphi$ we would like to prove the non-negativity of $\int_{S^2 \times S^1} 
\varphi P_{(J_0,\check{\th})} \varphi \, \check{\th} \wedge d \check{\th}$. 

Using the conformal covariance of $P$, see \eqref{eq:transfP}, we obtain
\begin{equation}\label{eq:pppp}
    \int_{S^2 \times S^1} 
    \varphi P_{(J_0,\check{\th})} \varphi \, \check{\th} \wedge d \check{\th} = 
    \int_{\{\rho \in [1,2] \}} \varphi P_{(J_0,\stackrel{\circ}{\th})} \varphi 
    \, \stackrel{\circ}{\th} \wedge d \stackrel{\circ}{\th}. 
\end{equation}
Extend next $\varphi$ by periodicity in $\log \rho$ to a smooth function $\tilde{\varphi}$ 
on $\H^1 \setminus \{0\}$:  for any positive integer $n$ we obtain that 
\begin{equation}\label{eq:nn}
    \int_{\{\rho \in [1,2]\}} \varphi P_{(J_0,\stackrel{\circ}{\th})} \varphi 
    \, \stackrel{\circ}{\th} \wedge d \stackrel{\circ}{\th} = \frac{1}{n} 
    \int_{\{\rho \in [1,2 ^n]\}} \tilde{\varphi} P_{(J_0,\stackrel{\circ}{\th})} \tilde{\varphi}
        \, \stackrel{\circ}{\th} \wedge d \stackrel{\circ}{\th}. 
\end{equation}
We next consider a smooth cut-off function $\chi$, depending only on $\rho$, such that 
$$
  \begin{cases}
  \chi(\rho) = 0 & \hbox{ for } \rho \leq \frac{5}{8}; \\ 
  \chi(\rho) = 1 & \hbox{ for } \rho \geq \frac{7}{8}, 
  \end{cases}
$$
and define the function $\tilde{\varphi}_n : \H \to \R$ as 
$$
    \tilde{\varphi}_n := \chi(\rho) \left[ 1 - \chi \left(  2^{-(n+1)} \rho \right) \right] \tilde{\varphi}.    
$$
By the non-negativity of $P_{(J_0,\stackrel{\circ}{\th})}$ (notice that the torsion vanishes identically) 
then we clearly have 
\begin{eqnarray*}
  0 & \leq & \int_{\H^1}  \tilde{\varphi}_n P_{(J_0,\stackrel{\circ}{\th})} \tilde{\varphi}_n 
           \, \stackrel{\circ}{\th} \wedge d \stackrel{\circ}{\th}   = 
          \int_{\{\rho \in [1,2 ^n]\}} \tilde{\varphi} P_{(J_0,\stackrel{\circ}{\th})} \tilde{\varphi}
                      \, \stackrel{\circ}{\th} \wedge d \stackrel{\circ}{\th} \\ & + &  \int_{\{\rho \in [1/2,1]\}} 
                      \tilde{\varphi}_n P_{(J_0,\stackrel{\circ}{\th})} \tilde{\varphi}_n 
                               \, \stackrel{\circ}{\th} \wedge d \stackrel{\circ}{\th}  
                      + \int_{\{\rho \in [2 ^n,2^{n+1}]\}} \tilde{\varphi}_n P_{(J_0,\stackrel{\circ}{\th})} 
                      \tilde{\varphi}_n  \, \stackrel{\circ}{\th} \wedge d \stackrel{\circ}{\th}. 
\end{eqnarray*}
By dilation invariance the last two terms in the above formula are uniformly bounded: therefore by \eqref{eq:pppp}, 
\eqref{eq:nn} we obtain that 
$$
  \int_{S^2 \times S^1} 
      \varphi P_{(J_0,\check{\th})} \varphi \, \check{\th} \wedge d \check{\th}  
      = \lim_{n \to + \infty}  \frac{1}{n} \int_{\H^1}  \tilde{\varphi}_n P_{(J_0,\stackrel{\circ}{\th})} \tilde{\varphi}_n 
      \, \stackrel{\circ}{\th} \wedge d \stackrel{\circ}{\th} \geq 0, 
$$ 
which is the desired conclusion.

\section{Proof of Theorem \ref{t:y}}\label{s:y}

The main goal of this section is to prove the inequality $\mathcal{Y}(J) <
\mathcal{Y}_0$ for a manifold $M$ as in Theorem \ref{t:y}. We first discuss the
existence of the Green's function and provide some estimates near the singularity.
Then we define suitable test functions for which the Sobolev quotient is below the
value of the standard pseudohermitian sphere.

\subsection{Green's function of $L_b$ and
test functions for the Tanaka-Webster quotient}\label{ss:green}

Given a compact CR manifold $M$ of positive Tanaka-Webster class and a
point $p \in M$, let us consider the equation
\begin{equation}\label{eq:deltap}
    - 4 \D_b G_p + R G_p = 16 \d_p,
\end{equation}
where $\d_p$ is the Dirac delta at the point $p$. We will use the
classical method of the parametrix to construct the Green's kernel.

\begin{lem}\label{l:greenheis} In the Heisenberg group one has
$$
  \stackrel{\circ}{\D}_b \frac{1}{\rho^2} = - 8 \pi \d_0.
$$
\end{lem}

\begin{pf} One can easily check that
$\stackrel{\circ}{\D}_b \frac{1}{\rho^2} = 0$ outside of
the origin of $\H^1$. We compute then the flux of the subgradient of $\frac{1}{\rho^2}$
on the boundary of any open set containing the origin. We can take for example the set
$$
  \{ \rho \leq 1 \}, \quad \hbox{ with boundary } S = \left\{ \rho = 1 \right\}.
$$
Let $f = \frac{1}{\rho^2}$: since $\stackrel{\circ}{\D}_b f = f_{, 1 \ou}
+ f_{, \ou 1}$, integrating by parts we get
\begin{equation}\label{eq:Dmass}
    \int_{\{ \rho \leq 1 \}} \stackrel{\circ}{\D}_b f \stackrel{\circ}{\th}
    \wedge d \stackrel{\circ}{\th} = \int_{\{ \rho \leq 1 \}} f_{, 1 \ou}
    \stackrel{\circ}{\th} \wedge d \stackrel{\circ}{\th}  +\hbox{conj.} =
    - i \oint_{S} f_{, 1} \stackrel{\circ}{\th^1} \wedge \stackrel{\circ}{\th}
    +\hbox{conj.}.
\end{equation}
By elementary computations one finds
$$
  f_{, 1} = - \frac{1}{\sqrt{2}} \ov{z} \frac{|z|^2 + i t}{\rho^6}; \qquad
  \qquad f_{, \ou} = - \frac{1}{\sqrt{2}} z \frac{|z|^2 - i t}{\rho^6},
$$
which implies
\begin{eqnarray*}
  \int_{\{ \rho \leq 1 \}} \stackrel{\circ}{\D}_b f \stackrel{\circ}{\th}
    \wedge d \stackrel{\circ}{\th} & = & \frac{i}{\sqrt{2}} \oint_S
    (|z|^2 + i t) \sqrt{2} \ov{z} dz \wedge (dt + i z d \ov{z}) +\hbox{conj.} \\
  & = & i \oint_S \left[ (|z|^2 + it) \ov{z} dz \wedge (dt + i z d \ov{z})
  - (|z|^2 - it) z d \ov{z} \wedge (dt - i \ov{z} dz) \right] \\
  & = & i \oint_S \left\{ |z|^2 \left[ (\ov{z} dz - z d \ov{z}) \wedge
  dt \right] + i t (2 i |z|^2 dz \wedge d \ov{z}) \right\}.
\end{eqnarray*}
Using \eqref{eq:relpol} and \eqref{eq:relpol2} we then deduce
$$
  \int_{\{ \rho \leq 1 \}} \stackrel{\circ}{\D}_b f \stackrel{\circ}{\th}
    \wedge d \stackrel{\circ}{\th} = - 2 \oint_S \left( |z|^4 d \var
    \wedge dt + t^2 d \var \wedge dt \right) = - 2 \oint_S d \var
    \wedge dt = - 8 \pi,
$$
so we obtain the conclusion. \end{pf}

\

\noindent To find a solution of \eqref{eq:deltap}, using CR normal
coordinates at $p$, we are going to consider a function $G_p$ of the
form
$$
  G_p = \frac{1}{2\pi} \rho^{-2} + w; \qquad \quad \rho^4 = |z|^4 + t^2,
$$
where $w$ has to be suitably chosen. First of all, we evaluate
$\D_b$ on $\rho^{-2}$. Using the fact that
$$
  |Z_1 \rho^i| \leq C_i \rho^{i-1}, \qquad
  |T \rho^i| \leq C_i \rho^{i-2}; \qquad \quad i \in \Z, i < 0,
$$
together with \eqref{errorXY}, we obtain
$$
  \D_b (\rho^{-2}) = - 8 \pi \d_p + g; \qquad \quad g = O(1).
$$
Therefore, we are reduced to solving the following equation for $w$
$$
  - 4 \D_b w + R w = g - \frac{1}{2 \pi} R \rho^{-2}.
$$
From the choice of CR normal coordinates (see Subsection \ref{ss:CRcoord})
it follows that $R = O(\rho^2)$, which implies that $w$ satisfies
$$
  - 4 \D_b w + R w = \tilde{g}; \qquad \quad \tilde{g} = O(1).
$$
The operator on the left-hand side is subelliptic, while the right-hand
side is in $L^q(M)$ for any $q > 1$. The regularity theory in \cite{FS}
then implies that $w \in \mathfrak{S}^{2,q}(M)$ for any $q > 1$. By the
Folland-Stein embeddings this also implies that $w \in C^{1,\g}(M)$
for any $\g \in (0,1)$. In conclusion we obtain the following
result.

\begin{pro}\label{p:exG}
Suppose $M$ is of positive Tanaka-Webster class, and that $p \in M$. Then
the conformal sublaplacian admits a Green's function $G_p$ with pole
at $p$, and in CR normal coordinates one has
$$
  G_p = \frac{1}{2 \pi \rho^2} + w,
$$
where $w$ is a function of class $C^{1,\g}$ for any $\g \in (0,1)$.
\end{pro}

\noindent By Theorem \ref{t:pm} we also deduce the following result.

\begin{pro}\label{p:Apos} Let $M$, $p$ and $G_p$ be as in Proposition \ref{p:exG}.
Then the following expansion holds
$$
  G_p = \frac{1}{2 \pi \rho^2} + A + \hat{w},
$$
where $A \geq 0$ and $\hat{w} \in C^{1,\g}$ for any $\g \in (0,1)$, and
$\hat{w}(p) = 0$. If $A = 0$, then $(M,J,\hat{\th})$ is pseudohermitian equivalent to the
standard three dimensional pseudohermitian sphere.
\end{pro}

\

\noindent To find suitable test functions we follow Schoen's idea, see
\cite{Sc}, and glue the Green's function to a {\em standard bubble},
which for the Heisenberg group is given by
$$
  \o_\l(z,t) = \l \left[ \l^4 t^2 + (1+\l^2 |z|^2)^2 \right]^{-\frac
  12} = \frac 1 \l \left( t^2 + |z|^4 + \frac{2}{\l^2} |z|^2 +
  \frac{1}{\l^4} \right)^{- \frac 12}; \qquad \l > 0.
$$
This function satisfies
\begin{equation}\label{eq:constR}
  - \stackrel{\circ}{\D}_b \o_\l = \o_\l^3
\end{equation}
and $\H^1$, endowed with the contact form $\o_\l^2 \stackrel{\circ}{\th}$
has constant positive Tanaka-Webster curvature (equal to $\frac 14$). Moreover, $\o_\l$
realizes the Tanaka-Webster quotient on $\H^1$ in the sense that
\begin{equation}\label{eq:best}
    \frac{\int_{\H^1} |\n_b \o_\l|^2 \stackrel{\circ}{\th} \wedge d
    \stackrel{\circ}{\th}}{\left( \int_{\H^1} |\o_\l|^4 \stackrel{\circ}{\th}
    \wedge d \stackrel{\circ}{\th} \right)^{\frac 12}} = \inf_{v \in \
    C^\infty_c(\H^1)} \frac{\int_{\H^1} |\n_b v|^2 \stackrel{\circ}{\th} \wedge d
    \stackrel{\circ}{\th}}{\left( \int_{\H^1} |v|^4 \stackrel{\circ}{\th}
    \wedge d \stackrel{\circ}{\th} \right)^{\frac 12}} = \mathcal{Y}_0.
\end{equation}
We now choose $\rho_0 \gg \frac 1 \l$, and define a cutoff function
$\psi$ which satisfies
$$
  \psi(z,t) = 1 \quad \hbox{ in } \{ \rho \leq \rho_0 \}; \qquad \qquad \psi(z,t)
  = 0 \quad \hbox{ in } \mathbb{H}^1 \setminus \{ \rho \geq 2 \rho_0 \}.
$$
Next, consider the following test function
\begin{equation}\label{eq:testfn}
  u(z,t) = \left\{
             \begin{array}{ll}
               \o_\l(z,t) & \hbox{ in } \{ \rho \leq \rho_0 \}; \\
               \e_0 \left( \tilde{G}_p(z,t) - \psi \tilde{w}(z,t) \right) + \psi(z,t)
      \frac 1 \l \var(z,t)
       & \hbox{ in } \{ \rho \leq 2 \rho_0 \} \setminus \{ \rho \leq \rho_0 \}; \\
               \e_0 \tilde{G}_p(z,t) & \hbox{ in } M \setminus \{ \rho \leq 2 \rho_0 \},
             \end{array}
           \right.
  \end{equation}
where we have set
$$
  \var(z,t) = \left( t^2 + |z|^4 + \frac{2}{\l^2} |z|^2
  + \frac{1}{\l^4} \right)^{- \frac 12} - \rho^{-2};
  \qquad \quad \tilde{G}_p = 2 \pi G_p; \qquad \tilde{w} = 2 \pi \hat{w}.
$$
The constant $\e_0$ has to be chosen so that the function $u$ is
continuous, namely so that
\begin{equation}\label{eq:e0}
  \e_0  = \frac{1}{\l (1 + \tilde{A} \rho_0^2)}; \qquad \quad
  \tilde{A} = 2 \pi A.
\end{equation}

\subsection{Estimates on the Tanaka-Webster quotient $\mathcal{Y}(J)$}\label{s:est}

We now evaluate the Sobolev-type quotient
\begin{equation}\label{eq:sobq}\tag{$Q$}
  \frac{\int_M \left( |\n_b u|^2 + \frac 14 R u^2 \right)
  \th \wedge d \th}{\left(\int_{M} |u|^4 \th \wedge d \th\right)^{\frac 12}}
\end{equation}
on the  test function in \eqref{eq:testfn}. We divide the estimate into four parts.

\

\noindent Recall from Subsection \ref{ss:bu} that for CR normal coordinates one has
\begin{equation}\label{errorXY}
    \left\{
    \begin{array}{ll}
      \th = \left( 1 + O(\rho^4) \right) \theta_0
     + O(\rho^5) d z + O(\rho^5) d \ov{z}; &  \\
      \th^1 = \sqrt{2} \left( 1 + O(\rho^4) \right) d z
     + O(\rho^4) d \ov{z} + O(\rho^3) d \theta_0,  &
    \end{array}
  \right.
\end{equation}
which implies
\begin{equation}\label{eq:CRvol}
  \th \wedge d \th = (1 + O(\rho^4)) \stackrel{\circ}{\theta}
  \wedge d  \stackrel{\circ}{\theta},
\end{equation}
and
\begin{equation}\label{eq:vvff}
    \left\{
    \begin{array}{ll}
      Z_1  = \left( 1 + O(\rho^4) \right) \stackrel{\circ}{Z}_1
      + O(\rho^4) \stackrel{\circ}{Z}_{\ou} + O(\rho^5)
      \stackrel{\circ}{T}; &  \\
      T = \left( 1 + O(\rho^4) \right)  \stackrel{\circ}{T}
     + O(\rho^3) \stackrel{\circ}{Z}_1 + O(\rho^3)
     \stackrel{\circ}{Z}_{\ou}. &
    \end{array}
  \right.
\end{equation}
Since $u$ is real, its squared subgradient is given by
$$
  |\n_b u|^2  = 2 (Z_1 u) (Z_{\ou} u).
$$

\subsubsection{Integral estimate of $|\n_b u|^2$}

First of all, we estimate $\int_{ \{ \rho \leq\rho_0 \} } |\n_b u|^2
\th \wedge d \th$ in terms of the flat subgradient. We can write
\begin{eqnarray}\label{eq:decgrad}\nonumber
  & & \int_{ \{ \rho \leq\rho_0 \} } |\n_b u|^2 \th \wedge d \th \\
  & = & 2 \int_{ \{ \rho \leq\rho_0 \} }
  (\stackrel{\circ}{Z}_1 u)  (\stackrel{\circ}{Z}_{\ou} u)   \stackrel{\circ}{\theta}
  \wedge d  \stackrel{\circ}{\theta}
  + 2 \int_{ \{ \rho \leq\rho_0 \} } ( (Z_1 u) (Z_{\ou} u) - (\stackrel{\circ}{Z}_1 u)
  (Z_{\ou} u) )
  \th \wedge d \th \\ & + & 2 \int_{ \{ \rho \leq\rho_0 \} } (\stackrel{\circ}{Z}_1 u)
  (\stackrel{\circ}{Z}_{\ou} u)  (\th \wedge d \th -  \stackrel{\circ}{\theta}
  \wedge d  \stackrel{\circ}{\theta}).  \nonumber
\end{eqnarray}
To evaluate the last term in the above expression we use
\eqref{eq:CRvol}, %and since $u$ is symmetric in $|z|$ in
%$B_{\rho_0}$, we can limit ourselves to the $(x,t)$ plane
and the fact that
\begin{equation}\label{eq:grad2}
  2 (\stackrel{\circ}{Z}_1 u) (\stackrel{\circ}{Z}_{\ou} u)  =
  \frac{|z|^2 \l^6}{\left( (1 + \l^2 |z|^2)^2 + \l^4 t^2 \right)^2},
\end{equation}
to obtain ($\tilde{\rho} = \l \rho$)
\begin{equation}\label{eq:diffvolgrad2}
  \left| \int_{ \{ \rho \leq\rho_0 \} } (\stackrel{\circ}{Z}_1 u)
  (\stackrel{\circ}{Z}_{\ou} u)  (\th \wedge d \th -  \stackrel{\circ}{\theta}
  \wedge d  \stackrel{\circ}{\theta}) \right| \leq C \int_0^{\rho_0}
  \frac{\rho^2 \l^6 \rho^4}{(1 + \l^4 \rho^4)^2} \rho^3 d \rho \leq
  \frac{C}{\l^2} \int_0^{\l \rho_0} \tilde{\rho} d \tilde{\rho} \leq
  C \frac{\rho_0^2}{\l^2}.
\end{equation}
To evaluate instead the second term in the right hand side of
\eqref{eq:decgrad}, we use \eqref{errorXY} to find
\begin{eqnarray} \label{eq:ciao} \nonumber
  \int_{ \{ \rho \leq\rho_0 \} } \left| (Z_1 u) (Z_{\ou} u) -
  (\stackrel{\circ}{Z}_1 u) (Z_{\ou} u) \right|
  \th \wedge d \th & \leq & C \int_{ \{ \rho \leq\rho_0 \} }
  |\stackrel{\circ}{\n}_b \o_\l| \left( \rho^4 |\stackrel{\circ}{\n}_b
  \o_\l| + \rho^5 |\stackrel{\circ}{T} \o_\l| \right) \stackrel{\circ}{\th}
  \wedge d \stackrel{\circ}{\th} \\ & + &  C \int_{ \{ \rho \leq\rho_0 \} }
  \left( \rho^8 |\stackrel{\circ}{\n}_b \o_\l|^2 + \rho^{10} |\stackrel{\circ}{T}
  \o_\l|^2 \right) \stackrel{\circ}{\th} \wedge d \stackrel{\circ}{\th}.
\end{eqnarray}
The first terms in the right hand side of the first line in
\eqref{eq:ciao} and in the second line in \eqref{eq:ciao} can be
treated exactly as for \eqref{eq:diffvolgrad2}. For the remaining
ones, we use the explicit expression
$$
 \stackrel{\circ}{T} \o_\l = - \frac{t \l^5}{\left( (1 + \l^2 |z|^2)^2 + \l^4 t^2
 \right)^{\frac 32}},
$$
to find
$$
  \int_{ \{ \rho \leq\rho_0 \} } |\stackrel{\circ}{\n}_b \o_\l| \; \rho^5
  \; |\stackrel{\circ}{T} \o_\l| \; \stackrel{\circ}{\th} \wedge d
  \stackrel{\circ}{\th} \, \leq C \int_0^{\rho_0} \frac{\rho^8 \l^8}{(1
  + \l^4 \rho^4)^{\frac 52}} \rho^3 d \rho \leq \frac{C}{\l^4} \int_0^{\l \rho_0}
  \frac{\tilde{\rho}^{11}}{(1 + \tilde{\rho}^4)^{\frac 52}} d \tilde{\rho} \leq
  \frac{C}{\l^2} \rho_0^2;
$$
$$
  \int_{ \{ \rho \leq\rho_0 \} }  \rho^{10} |\stackrel{\circ}{T} \o_\l|^2
  \stackrel{\circ}{\th} \wedge d \stackrel{\circ}{\th} \,
  \leq C \int_0^{\rho_0} \frac{\rho^{14} \l^{10}}{(1
  + \l^4 \rho^4)^3} \rho^3 d \rho \leq \frac{C}{\l^8} \int_0^{\l \rho_0}
  \frac{\tilde{\rho}^{17}}{(1 + \tilde{\rho})^{12}} d \tilde{\rho} \leq
  C \frac{\rho_0^6}{\l^2}.
$$
In conclusion, using \eqref{eq:decgrad}, \eqref{eq:diffvolgrad2},
\eqref{eq:ciao} and the last two formulas we obtain that
\begin{equation}\label{eq:estfinalgrad}
  \left| \int_{ \{ \rho \leq\rho_0 \} } |\n_b u|^2 \th \wedge d \th -
  \int_{ \{ \rho \leq\rho_0 \} } |\stackrel{\circ}{\n}_b u|^2
  \stackrel{\circ}{\th} \wedge d \stackrel{\circ}{\th} \right| \leq
   C \frac{\rho_0^2}{\l^2}.
\end{equation}

\subsubsection{Integral estimate of $R u^2$}

We estimate next the term $\int_{ \{ \rho \leq\rho_0 \} } R u^2 \th \wedge d
\th$. For doing this, we first expand in Taylor series the Tanaka-Webster curvature
$R$, using real CR normal coordinates. From the relations (see Subsection
\ref{ss:CRcoord})
\begin{equation}\label{eq:CrWebster}
    R = 0; \qquad R_{, 1} = 0; \qquad \D_b R = R_{, 1 \ov{1}}
    + R_{, \ov{1} 1} = 0; \qquad R_{, 0} = 0,
\end{equation}
and the fact that for a smooth function $f$ one has
$$
  \stackrel{\circ}{\D}_b f = f_{xx} + f_{yy} + 2y f_{xt} - 2x f_{yt} + 4
  |z|^2 f_{tt} \qquad \quad \hbox{ at } p,
$$
we can write
\begin{equation}\label{eq:expebster}
     R(z,t) = \frac{1}{2} \left( R_{xx} x^2 + R_{yy} y^2 + 2 R_{xy} xy \right)
  + O(\rho^3), \qquad \hbox{ with } \quad R_{yy} = - R_{xx}.
\end{equation}
By this reason and by the axial symmetry of $\o_\l$ the terms of
order $\rho^2$, integrated, will vanish. Similarly, by symmetry,
also the cubic ones will vanish, and therefore by \eqref{eq:CRvol}
we deduce
\begin{eqnarray} \label{eq:webin} \nonumber
  \left| \int_{ \{ \rho \leq\rho_0 \} } R u^2 \th \wedge d \th \right| & = &
  \left| \int_{ \{ \rho \leq\rho_0 \} } R \o_\l^2 \stackrel{\circ}{\th} \wedge d \stackrel{\circ}{\th}
  \right| + O \left( \int_{ \{ \rho \leq\rho_0 \} } \rho^6 \o_\l^2 \stackrel{\circ}{\th}
  \wedge d \stackrel{\circ}{\th} \right) \leq C \int_{ \{ \rho \leq\rho_0 \} } \rho^4 \o_\l^2
  \stackrel{\circ}{\th} \wedge d \stackrel{\circ}{\th} \\ & \leq & C \int_0^{\rho_0}
  \frac{\l^2}{(1 + \l^4 \rho^4)} \rho^4 d \rho \leq \frac{C}{\l^3}
  \int_0^{\l \rho_0} \frac{\tilde{\rho}^4}{(1 + \tilde{\rho})^4} d
  \tilde{\rho} \leq \frac{C \rho_0}{\l^2}.
\end{eqnarray}

\subsubsection{Estimate of the boundary term}

From the Stokes theorem, for a smooth domain $\O \subseteq \H^1$ we get
\begin{eqnarray}\label{eq:intparts} \nonumber
  \int_\O |\stackrel{\circ}{\n}_b u|^2 \stackrel{\circ}{\th} \wedge d \stackrel{\circ}{\th}
  & = & \int_\O \stackrel{\circ}{Z}_1 u \stackrel{\circ}{Z}_{\ou} u
  \stackrel{\circ}{\th} \wedge d \stackrel{\circ}{\th}
  + \int_\O \stackrel{\circ}{Z}_{\ou} u \stackrel{\circ}{Z}_1 u
  \stackrel{\circ}{\th} \wedge d \stackrel{\circ}{\th} \\
  & = & i \oint_{\pa \O} u (\stackrel{\circ}{Z}_{\ou} u)
  \stackrel{\circ}{\th^{\ou}} \wedge \stackrel{\circ}{\th}
  - i \oint_{\pa \O} u (\stackrel{\circ}{Z}_1 u)
  \stackrel{\circ}{\th^1} \wedge \stackrel{\circ}{\th}
  - \int_\O u \stackrel{\circ}{\D}_b u \stackrel{\circ}{\th} \wedge d
  \stackrel{\circ}{\th}.
\end{eqnarray}
Next we take $\O = \{ \rho < \rho_0 \} $ and $u = \o_\l$, focusing our
attention to the boundary integral. From elementary computations we find
$$
  \stackrel{\circ}{Z}_1 u = - \frac{1}{\sqrt{2}} \frac{\ov{z} (1+\l^2|z|^2 +
  i t \l^2)\l^3}{\left[ \l^4 t^2 + (1+\l^2|z|^2)^2 \right]^{\frac 32}}; \qquad
  \qquad \stackrel{\circ}{Z}_{\ou} u = - \frac{1}{\sqrt{2}} \frac{z
  (1+\l^2|z|^2 - i t \l^2)\l^3}{\left[ \l^4 t^2 + (1+\l^2|z|^2)^2
  \right]^{\frac 32}},
$$
hence after some manipulations the above boundary integral becomes
\begin{equation}\label{eq:bdintbub}
  - i \oint_{ S_{\rho_0} } \frac{z (1+\l^2|z|^2 - i t \l^2)
  \l^4}{\left[ \l^4 t^2 + (1+\l^2|z|^2)^2 \right]^2} d \ov{z}
  \wedge (dt - i \ov{z} dz) + \hbox{conj.}.
\end{equation}
On the other hand, the function $\o_\l$ satisfies \eqref{eq:constR},
therefore from the above computations we get
\begin{eqnarray}\label{eq:contrbub}\nonumber
  \int_{ \{ \rho \leq\rho_0 \} } |\stackrel{\circ}{\n}_b u|^2
  \stackrel{\circ}{\th} \wedge d \stackrel{\circ}{\th} & = &
  \int_{ \{ \rho \leq\rho_0 \} } \o_\l^4
  \stackrel{\circ}{\th} \wedge d \stackrel{\circ}{\th} \\
  & - & i
  \oint_{ S_{\rho_0} } \frac{z (1+\l^2|z|^2 - i t \l^2)
  \l^4}{\left[ \l^4 t^2 + (1+\l^2|z|^2)^2 \right]^2} d \ov{z}
  \wedge (dt - i \ov{z} dz) + \hbox{conj.} \\ & = &
  - 2 \oint_S \frac{\l^4 (|z|^2 + \l^2 \rho^4)}{\left[
  \l^4 t^2 + (1+\l^2|z|^2)^2 \right]^2} d \var \wedge dt. \nonumber
\end{eqnarray}
We next evaluate the integrals on $M \setminus \{ \rho \leq \rho_0 \} $. We have
\begin{eqnarray}\label{eq:intext} \nonumber
  & & \int_{M \setminus \{ \rho \leq \rho_0 \} } \left( |\n_b u|^2 + \frac 14 R u^2 \right)
  \th \wedge d \theta = \e_0^2 \int_{M \setminus \{ \rho \leq \rho_0 \} } \left(
  |\n_b \tilde{G}_p|^2 + \frac 14 R \tilde{G}_p^2 \right) \th \wedge d \theta
  \\ & + & \frac{1}{\l^2} \int_{\{ \rho \leq 2 \rho_0 \} \setminus \{ \rho \leq \rho_0 \} }
  \nonumber |\n_b (\var \psi)|^2 \th \wedge d \theta + \frac 14
  \int_{\{ \rho \leq 2 \rho_0 \} \setminus \{ \rho \leq \rho_0 \} } R (u^2 - \e_0^2 \tilde{G}_p^2)
  \th \wedge d \theta \\ & + &
  \e_0^2 \int_{\{ \rho \leq 2 \rho_0 \} \setminus
\{ \rho \leq \rho_0 \} } \left( |\n_b (\psi \tilde{w})|^2
  - 2 \n_b \tilde{G}_p \cdot \n_b (\psi \tilde{w}) \right) \th \wedge d \theta \\ & + &
  2 \frac{\e_0}{\l} \int_{\{ \rho \leq 2 \rho_0 \} \setminus \{ \rho \leq \rho_0 \} } \left( \n_b
  \tilde{G}_p \cdot \n_b (\var \psi) - \n_b (\psi \tilde{w}) \cdot \n_b (\var \psi)
  \right) \th \wedge d \theta.  \nonumber
\end{eqnarray}
We derive now some estimates on the above terms. First of all we
notice that we can take
\begin{equation}\label{eq:scal}
    \psi(z,t) = \psi_0 \left( \frac{z}{\rho_0}, \frac{t}{\rho_0^2}
\right),
\end{equation}
for a smooth fixed $\psi_0$ which, using \eqref{eq:vvff}, implies
$$
  Z_1 \psi = \frac{1}{\rho_0} \pa_z \psi_0 \left( \frac{z}{\rho_0}, \frac{t}{\rho_0^2}
\right) + i \frac{\ov{z}}{\rho_0^2} \pa_t \psi_0 \left( \frac{z}{\rho_0},
\frac{t}{\rho_0^2} \right) + O(\rho_0^3) + O(\rho_0) = O \left( \frac{1}{\rho_0} \right).
$$
A similar estimate holds true for $Z_{\ou} \psi$ and therefore, since
$\tilde{w}(z,t) = O(\rho)$  we have
$$
  |\n_b (\psi \tilde{w})| \leq C  \qquad \quad \hbox{ in } \{ \rho \leq 2 \rho_0 \}
  \setminus \{ \rho \leq \rho_0 \}  \qquad \hbox{ for some fixed } C > 0.
$$
Using also the fact that $\l \rho_0 \gg 1$, from a Taylor expansion
we find
\begin{eqnarray}\label{eq:bau} \nonumber
    \var & = & \rho^{-2} \left[ \left( 1 + \frac{2|z|^2}{\l^2 \rho^4} +
  \frac{1}{\l^4 \rho^4} \right)^{-\frac 12} - 1 \right] = \frac{1}{\rho^2} O \left(
  \frac{2|z|^2}{\l^2 \rho^4} + \frac{1}{\l^4 \rho^4} \right)
  \\ & = & \frac{1}{\l^4 \rho^6} O \left( 1 + \l^2 |z|^2 \right) = O \left(
  \frac{1}{\l^2 \rho^4} \right).
\end{eqnarray}
By straightforward computations we also get
$$
  \stackrel{\circ}{Z}_1 \var = \frac{1}{\sqrt{2}} \frac{\ov{z} \rho}{\l^2}
  \frac{\l^6 \rho^6}{(\l^4 t^2 + (1 + \l^2 |z|^2)^2)^{\frac 32}}
  \left\{ \frac{\l^2(|z|^2 + i t)}{\l^6 \rho^6} \left[ (\l^4 t^2 +
   (1 + \l^2 |z|^2)^2)^{\frac 32} - \l^6 \rho^6 \right] - 1 \right\},
$$
and similarly
$$
  \stackrel{\circ}{Z}_{\ou} \var = \frac{1}{\sqrt{2}} \frac{z \rho}{\l^2}
  \frac{\l^6 \rho^6}{(\l^4 t^2 + (1 + \l^2 |z|^2)^2)^{\frac 32}}
  \left\{ \frac{\l^2(|z|^2 - i t)}{\l^6 \rho^6} \left[ (\l^4 t^2 +
   (1 + \l^2 |z|^2)^2)^{\frac 32} - \l^6 \rho^6 \right] - 1 \right\}.
$$
Using some Taylor expansions, the fact that $|\stackrel{\circ}{T} \var| =
O\left( \frac{1}{\l^2 \rho^6} \right)$ and \eqref{errorXY} we then
deduce
$$
  |\n_b \var| = O \left( \frac{1}{\l^2 \rho_0^5} \right) \qquad \quad
  \hbox{ in } \{ \rho \leq 2 \rho_0 \} \setminus \{ \rho \leq \rho_0 \} .
$$
From the above estimates on $\psi$ this implies
\begin{equation}\label{eq:estfipsi}
  |\n_b (\var \psi)| = O \left( \frac{1}{\l^2 \rho_0^5} \right)
  \qquad \quad \hbox{ in } \{ \rho \leq 2 \rho_0 \} \setminus \{ \rho \leq \rho_0 \} .
\end{equation}
It is also easy to see that
\begin{equation}\label{eq:estnG}
   |\n_b \tilde{G}_p| = O\left(\frac{1}{\rho_0^3}\right) \qquad
   \quad \hbox{ in } \{ \rho \leq 2 \rho_0 \} \setminus \{ \rho \leq \rho_0 \} .
\end{equation}
Next, from \eqref{eq:testfn}, \eqref{eq:e0}, the expression of $\tilde{G}_p$,
$\tilde{w} = O(\rho)$, \eqref{eq:bau} and $\l \rho_0 \gg 1$ we get
$$
  |u^2 - \e_0^2 \tilde{G}_p^2| \leq \e_0 \rho_0^{-2} \left( \e_0 \rho_0 + \frac{1}{\l^3
  \rho_0^4} \right) + \e_0^2 \rho_0^2 + \frac{1}{\l^6 \rho_0^8} \leq
  C \left( \frac{1}{\l^2 \rho_0} + \frac{1}{\l^4 \rho_0^6} \right).
$$
This and \eqref{eq:expebster} imply
$$
  \frac 14 \int_{\{ \rho \leq 2 \rho_0 \} \setminus \{ \rho \leq \rho_0 \} } R (u^2 - \e_0^2
  \tilde{G}_p^2) \th \wedge d \theta \leq C \left( \frac{\rho_0^5}{\l^2}
  + \frac{1}{\l^4} \right).
$$
Therefore, integrating and using the fact that the volume of $B_{2
\rho_0} \setminus \{ \rho \leq \rho_0 \} $ is of order $\rho_0^4$, from
\eqref{eq:intext} and the above estimates we obtain
\begin{eqnarray}\label{eq:intext2} \nonumber
  \int_{M \setminus \{ \rho \leq \rho_0 \} } \left( |\n_b u|^2 +
\frac 14 R u^2 \right) \th \wedge d \theta & = &
  \e_0^2 \int_{M \setminus \{ \rho \leq \rho_0 \} } \left( |\n_b \tilde{G}_p|^2
+ \frac 14 R \tilde{G}_p^2 \right) \th \wedge d
  \theta \\ & + & O \left( \frac{1}{\l^6 \rho_0^6} + \e_0^2 \rho_0 +
\frac{\e_0}{\l^3 \rho_0^4} \right).
\end{eqnarray}
The first term in this expression can be evaluated using an
integration by parts, yielding a boundary term: recall that the
Green's function satisfies $- \D_b \tilde{G}_p + \frac 14 R \tilde{G}_p = 0$ outside the
singularity. Therefore, as for \eqref{eq:intparts} one finds
$$
  \int_{M \setminus \{ \rho \leq \rho_0 \} } \left( |\n_b \tilde{G}_p|^2
+ \frac 14 R \tilde{G}_p^2 \right)
  \th \wedge d \theta = - i \oint_{S_{\rho_0}} \tilde{G}_p
(Z_{\ou} \tilde{G}_p) \th^{\ou} \wedge \th
  + i \oint_{S_{\rho_0}} \tilde{G}_p (Z_1 \tilde{G}_p) \th^1 \wedge \th,
$$
the minus sign coming from the fact that we are working in the exterior of
$\{ \rho \leq \rho_0 \}$.

By direct computations we have that
\begin{eqnarray*}
   & & - i \tilde{G}_p (Z_{\ou} \tilde{G}_p) \th^{\ou} \wedge \th + i \tilde{G}_p
(Z_1 \tilde{G}_p) \th^1 \wedge
  \th  \\
   & = & - i \tilde{G}_p \left[ \frac{1}{\sqrt{2}} \frac{\ov{z}(|z|^2 + i t)}{\rho^6}
  \sqrt{2} dz \wedge (dt + i z d \ov{z}) - \frac{1}{\sqrt{2}}
  \frac{z (|z|^2 - i t)}{\rho^6} \sqrt{2} d\ov{z} \wedge (dt - i \ov{z} dz)
  + O(\rho) d \var \wedge dt \right].
\end{eqnarray*}
Using some cancelations we get
\begin{equation}\label{eq:intG2}
  \int_{M \setminus \{ \rho \leq \rho_0 \} } \left(
  |\n_b \tilde{G}_p|^2 + \frac 14 R \tilde{G}_p^2 \right)
  \th \wedge d \theta = \oint_{S_{\rho_0} } \left( \frac{1}{\rho^2}
  + \tilde{A} + \tilde{w} \right) \left( \frac{2}{\rho^2}
  + O(\rho) \right) d \var \wedge dt,
\end{equation}
%It can be shown that, as $\rho \to 0$ (see Subsection \ref{ss:pf}
%below for the proof)
%\begin{equation}\label{eq:surprise}
%    \oint_{S_{\rho_0} } \frac{1}{4|z|^6+t^2}  d \hat{\s} =
%    O \left( \frac{1}{\rho_0^2} \right),
%\end{equation}
therefore after integration we obtain
\begin{eqnarray}\label{eq:intG3}
   \int_{M \setminus \{ \rho \leq \rho_0 \} } \left( |\n_b \tilde{G}_p|^2
   + \frac 14 R \tilde{G}_p^2 \right) \th \wedge d \theta
  = \oint_{S_{\rho_0} }  \frac{2 (1+\tilde{A}\rho^2)}{\rho^4}
  d \var \wedge dt + O(\rho_0^2).
\end{eqnarray}
Now, using the formula for $\e_0$ in \eqref{eq:e0} and some
straightforward computations one finds that the contribution of the
 boundary term in \eqref{eq:contrbub} together with the one in the last
formula becomes (there is cancelation of the two main terms after
subtraction)
$$
  \oint_{S_{\rho_0} } 2 \frac{3 \l^6 \rho^4 |z|^2
+ 1 + 4 \l^4 |z|^4 + 2 \l^4 \rho^4 + 4 \l^2 |z|^2 - \tilde{A} \l^8 \rho^{10}
- \tilde{A} \rho^6 \l^6 |z|^2}{\rho^4 \l^{10} (1+\tilde{A} \rho^2) \left( \rho^4 +
\frac{1}{\l^4} + 2 \frac{|z|^2}{\l^2} \right)^2} d \var \wedge dt.
$$
To get some asymptotics of the integrand (especially its sign) we can neglect
the term $\tilde{A} \rho^2$ in the denominator, while in the second bracket of
the denominator we can simply take $(\rho^4)^2$. For the numerator, we
recall that $\l \rho_0 \gg 1$. In this way we have to find the
asymptotics of
$$
  \oint_{S_{\rho_0} } 2
  \frac{3 |z|^2 - \tilde{A} \l^2 \rho^6}{\rho^8 \l^4} d \var \wedge dt.
$$
If we also choose $\l$ and $\rho_0$ so that $\l^2 \rho_0^4 \gg 1$
then the last term in the numerator dominates, and we are reduced to
$$
  - 2 \frac{\tilde{A}}{\rho_0^2 \l^2} \oint_{S_{\rho_0} }
  d \var \wedge dt = - 8 \pi \frac{\tilde{A}}{\l^2}.
$$

\subsubsection{Final estimates}

Collecting all the  terms in
\eqref{eq:webin}, \eqref{eq:contrbub}, \eqref{eq:intext2}, \eqref{eq:intG2} and
\eqref{eq:intG3}, and taking into account the previous computations,
with our choices of $\rho_0$ and $\l$ we  find
\begin{eqnarray}\label{eq:estfingrad} \nonumber
  \int_M \left( |\n_b u|^2 + \frac 14 R u^2 \right)
  \th \wedge d \th &=& \int_{ \{ \rho \leq\rho_0 \} }
  \o_\l^4 \stackrel{\circ}{\th} \wedge d \stackrel{\circ}{\th}
  - 8 \pi \frac{\tilde{A}}{\l^2} + O \left( \frac{1}{\l^6 \rho_0^6} \right)
  \\ & + & O \left( \frac{\rho_0}{\l^2} \right)
  + O \left( \frac{1}{\l^4 \rho_0^4} \right).
\end{eqnarray}
From $\l \rho_0 \gg 1$ we get
\begin{equation}\label{eq:estfingrad2}
    \int_M \left( |\n_b u|^2 + \frac 14 R u^2 \right) \th \wedge d \th =
  \int_{ \{ \rho \leq\rho_0 \} } \o_\l^4 \stackrel{\circ}{\th} \wedge
  d \stackrel{\circ}{\th} - (8 \pi + o(1)) \frac{\tilde{A}}{\l^2}.
\end{equation}

\

\noindent Finally, we estimate the denominator in \eqref{eq:sobq}. Clearly the
contribution inside $\{ \rho \leq \rho_0 \} $ is just the same as the one
of the bubble. In
$\{ \rho \leq 2 \rho_0 \} \setminus \{ \rho \leq \rho_0 \} $ we have that $|u| \leq
\frac{C}{\l \rho_0^2}$ (by the expression of $\e_0$ and he
asymptotics of $\tilde{G}_p$), so integrating we get
$$
  \int_{\{ \rho \leq 2 \rho_0 \} \setminus \{ \rho \leq \rho_0 \} } |u|^4 \th \wedge d \th
  = O \left( \frac{1}{\l^4 \rho_0^4} \right) \ll \frac{1}{\l^2}.
$$
Similarly we also find
$$
  \int_{M \setminus \{ \rho \leq 2 \rho_0 \}} |u|^4 \th \wedge d \th \leq C \frac{1}{\l^4}
  \int_{\rho_0}^\infty \frac{\rho^3 d \rho}{\rho^8} = O \left(
  \frac{1}{\l^4 \rho_0^4} \right) \ll \frac{1}{\l^2}.
$$
Moreover, from \eqref{eq:CRvol} and the expression of $\o_\l$ we
obtain
$$
   \int_{ \{ \rho \leq\rho_0 \} } |u|^4 \th \wedge d \th = \int_{ \{ \rho \leq\rho_0 \} } \o_\l^4
   \stackrel{\circ}{\th} \wedge d \stackrel{\circ}{\th} + O \left( \frac{1}{\l^4} \log (\l \rho_0)
  \right).
$$
In conclusion, using \eqref{eq:estfingrad2} and the last three
formulas, the Tanaka-Webster quotient \eqref{eq:sobq} for the test function
\eqref{eq:testfn} becomes
$$
  \frac{\int_M |\n_b u|^2 \th \wedge d \th}{\left( \int_{M} |u|^4 \th \wedge d
\th \right)^{\frac 12}} = \frac{\int_{ \{ \rho \leq\rho_0 \} }
\o_\l^4 \stackrel{\circ}{\th} \wedge d \stackrel{\circ}{\th} -
(C_1 + o(1)) \frac{\tilde{A}}{\l^2 \rho_0^2}}{\left(\int_{ \{
\rho \leq\rho_0 \} } \o_\l^4 \stackrel{\circ}{\th}
\wedge d \stackrel{\circ}{\th} \right)^{\frac 12} + O \left( \frac{1}{\l^4 \rho_0^4}
+ O \left(\frac{1}{\l^4} \log (\l \rho_0) \right)\right)}.
$$
Since $O \left( \frac{1}{\l^4 \rho_0^4} \right) + O \left(
\frac{1}{\l^4} \log (\l \rho_0) \right)\ll \frac{1}{\l^2}$,
from the last formula and \eqref{eq:best} we find
$$
  \frac{\int_M |\n_b u|^2 \th \wedge d \th}{\left( \int_{M} |u|^4 \th \wedge d
\th \right)^{\frac 12}}  <
\frac{\int_{\H^1} |\n_b \o_l|^2 \stackrel{\circ}{\th} \wedge d
    \stackrel{\circ}{\th}}{\left( \int_{\H^1} |\o_l|^4 \stackrel{\circ}{\th}
    \wedge d \stackrel{\circ}{\th} \right)^{\frac 12}} = \mathcal{Y}_0,
$$
which is the desired inequality.

\

\noindent The proof of Theorem \ref{t:y} then follows from \cite{JL1/2} in the
case when $M$ is CR equivalent to the 3-sphere, and from the variational
argument in \cite{JL0} in the complementary case.

\begin{rem}\label{r:nonmin} Since the minimization procedure for the Tanaka-Webster quotient
is related to the positivity of the mass, in view of Proposition \ref{p:example} there
might be examples of pseudohermitian manifolds of positive Tanaka-Webster class for which there
is no minimizer. In this case, the use of variational
or topological methods (as in \cite{GA}, \cite{GJ}) seems to be a necessary tool
in order to find conformal structures of constant Tanaka-Webster curvature.
\end{rem}

\section{Appendix: useful facts in pseudohermitian geometry}\label{s:app}

In this appendix we collect some useful facts in pseudohermitian geometry:
in particular we discuss the variations of some geometric quantities when
the contact form or the CR structure vary, and then the CR normal
coordinates introduced by Jerison and Lee in \cite{JL}.

\subsection{Variations of interesting quantities}\label{ss:variations}

Here we consider variations of some geometric quantities either under the
conformal change of contact form (pp. 421-422 in \cite{L1}), or under the
change of the pseudohermitian structure (pp. 231-232 in \cite{CL1}).

\subsubsection{Conformal changes of contact form}

\label{sss:conf}

Consider a new contact form $\hat{\th }=e^{2f}\th $, for a given smooth
function $f$. Then it turns out  that
\begin{equation}
\hat{\th }^{1}=e^{f}\left( \th ^{1}+2if^{1}\th \right) ;\qquad \quad \hat{h}%
_{1\overline{1}}=h_{1\overline{1}}\equiv 1,  \label{eq:th1f}
\end{equation}%
and that
\begin{equation}\label{eq:Z1conf}
\hat{Z}_{1}=e^{-f}Z_{1};\qquad \quad \hat{T}=e^{-2f}\left( T+2if^{\overline{1%
}}Z_{\overline{1}}-2if^{1}Z_{1}\right) .
\end{equation}%
Moreover we also have
\begin{equation}
\hat{\o }_{1}^{1}=\o _{1}^{1}+3\left( f_{1}\th ^{1}-f_{\overline{1}}\th ^{%
\overline{1}}\right) +i\left( f_{\overline{1}1}+f_{1\overline{1}}+8f_{1}f_{%
\overline{1}}\right) \th .  \label{eq:varoconf}
\end{equation}%
Furthermore, concerning the torsion and the Tanaka-Webster curvature one has
\begin{equation}
\hat{A}_{11}=e^{-2f}\left( A_{11}+2if_{11}-4if_{1}f_{1}\right) ;\qquad \quad
\hat{R}=e^{-2f}\left( R-4\D_{b}f-8f_{1}f_{\overline{1}}\right) .
\label{eq:varAconf}
\end{equation}

\subsubsection{Deformation of $J$}\label{sss:varJ}

Here we consider instead a variation of the CR structure $J(s)$ for which
\begin{equation*}
  \frac{d}{ds}_{|s=0} J(s) = \dot{J}=2E=2E_{11}\th ^{1}\otimes
Z_{\overline{1}}+2E_{\overline{1}\overline{%
1}}\th ^{\overline{1}}\otimes Z_{1}.
\end{equation*}%
This implies
\begin{equation}
\dot{Z}_{1}=-iE_{1}^{\overline{1}}Z_{\overline{1}};\qquad \quad \dot{\th }%
^{1}=-iE_{\overline{1}}^{1}\th ^{\overline{1}},  \label{eq:varZ}
\end{equation}%
and
\begin{equation}
\dot{\o }_{1}^{1}=i\left( A_{11}E_{\overline{1}\overline{1}}+A_{\overline{1}%
\overline{1}}E_{11}\right) \th -iE_{\overline{1},1}^{1}\th ^{\overline{1}%
}-iE_{1,\bar{1}}^{\bar{1}} \th ^{1};\qquad \quad {\dot{A}}_{%
\overline{1}}^{1}=-iE_{\overline{1},0}^{1}.  \label{eq:varoA}
\end{equation}%
We also have
\begin{equation}
\dot{R}=i\left( E_{11,\overline{1}\overline{1}}-E_{\overline{1}\overline{1}%
,11}\right) -\left( A_{11}E_{\overline{1}\overline{1}}+A_{\overline{1}%
\overline{1}}E_{11}\right) .  \label{eq:varW}
\end{equation}%
It is  useful to include the variation of $- \D_{b}$ as well: one has
\begin{equation*}
\D_{b}=\left( Z_{\overline{1}}Z_{1}-\o _{1}^{1}(Z_{\overline{1}%
})Z_{1}\right) +\hbox{conj.}.
\end{equation*}%
From this formula, differentiating with respect to $s$ one finds
\begin{eqnarray*}
- \dot{\D}_{b} &=&-iE_{\overline{1}\overline{1}}Z_{1}Z_{1}+iE_{11,\overline{1}%
}Z_{\overline{1}}+iE_{11}Z_{\overline{1}}Z_{\overline{1}} \\
&-&iE_{\overline{1}\overline{1},1}Z_{1}+iE_{\overline{1}\overline{1}}\o %
_{1}^{1}(Z_{1})Z_{1}-iE_{11}\o _{1}^{1}(Z_{\overline{1}})Z_{\overline{1}}+%
\hbox{conj.} \\
&=&2\left( -iE_{\overline{1}\overline{1}}Z_{1}Z_{1}+iE_{11}Z_{\overline{1}%
}Z_{\overline{1}}+iE_{11,\overline{1}}Z_{\overline{1}}-iE_{\overline{1}%
\overline{1},1}Z_{1}\right) \\
&+&iE_{\overline{1}\overline{1}}\o _{1}^{1}(Z_{1})Z_{1}-iE_{11}\o %
_{1}^{1}(Z_{\overline{1}})Z_{\overline{1}}-iE_{11}\o _{\overline{1}}^{%
\overline{1}}(Z_{\overline{1}})Z_{\overline{1}}+iE_{\overline{1}\overline{1}}%
\o _{\overline{1}}^{\overline{1}}(Z_{1})Z_{1}.
\end{eqnarray*}%
From the fact that $\o _{1}^{1}$ is purely imaginary we find
\begin{equation*}
- \dot{\D}_{b}=2i\left( E_{11}Z_{\overline{1}}Z_{\overline{1}}+E_{11,\overline{%
1}}Z_{\overline{1}}\right) +\hbox{conj.}.
\end{equation*}%
To derive these formulas we have used the above transformation laws for $J$,
$Z_{1}$ and $\o _{1}^{1}$.

If we differentiate once more we find that
\begin{equation*}
(E_{11,\overline{1}})^{\cdot }=iE_{\overline{1}\overline{1}}E_{11,1}+2iE_{%
\overline{1}\overline{1},1}E_{11}+(\dot{E}_{11})_{,\overline{1}},
\end{equation*}%
and hence
\begin{eqnarray}
- \ddot{\D}_{b} &=&2i\left( \dot{E}_{11}Z_{\overline{1}}Z_{\overline{1}}+E_{11}%
\dot{Z}_{\overline{1}}Z_{\overline{1}}+E_{11}Z_{\overline{1}}\dot{Z}_{%
\overline{1}}+(E_{11,\overline{1}})^{\cdot }Z_{\overline{1}}+E_{11,\overline{%
1}}\dot{Z}_{\overline{1}}\right) +\hbox{conj.}  \notag  \label{eq:ddotD} \\
&=&2i\dot{E}_{11}Z_{\overline{1}}Z_{\overline{1}}-2i\dot{E}_{\bar{1}\bar{1}%
}Z_{1}Z_{1}-4|E_{11}|^{2}\D_{b}-\left( 4E_{11}E_{\overline{1}\overline{1},%
\overline{1}}+6E_{\overline{1}\overline{1}}E_{11,\overline{1}}\right) Z_{1}
\\
&&-\left( 4E_{\overline{1}\overline{1}}E_{11,1}+6E_{11}E_{\overline{1}%
\overline{1},1}\right) Z_{\overline{1}}+2i(\dot{E}_{11})_{,\overline{1}}Z_{%
\overline{1}}-2i(\dot{E}_{\overline{1}\overline{1}})_{,1}Z_{1}.  \notag
\end{eqnarray}
Using similar formulas we can also deduce the expression of $\ddot{R}$.

\subsection{CR normal coordinates}
\label{ss:CRcoord}

We recall the following result in \cite{JL} on page 313, Proposition 2.5.

\begin{pro}\label{p:phcoord} Let ${\tilde{Z}}_1$ be a special frame dual to
${\tilde{\th}}^1$ (such that ${\tilde{h}}_{1 \ou}=2$). Let ${\th^1}={\sqrt{2}}
{\tilde{\th}}^1$ be a unitary coframe (i.e., $h_{1 \ou}=1$). Then in
pseudohermitian normal coordinates $(z, t)$ with respect to ${\tilde{Z}}_1$,
${\tilde{\th}}^1$, we have

\noindent $(a)$ $\th_{(2)} = \theta_0; \qquad \th_{(3)} = 0; \qquad
\th_{(m)} =
 \frac 1 m {\sqrt 2} \left( i  z \th^{\ou} - i \ov{z} \th^1 \right)_{(m)},
\quad m \geq 4$;

\noindent $(b)$ $\th^1_{(1)} = {\sqrt 2}d z; \quad \th^1_{(2)} = 0; \qquad
  \th^1_{(m)} = \frac 1 m \left( {\sqrt 2}z \o^1_1 + 2 t A_{\ou \ou} \th^{\ou}
  - {\sqrt 2}\ov{z} A_{\ou \ou} \th  \right)_{(m)}, \quad m \geq 3$

\noindent $(c)$  $(\o^1_1)_{(1)} = 0; \quad (\o^1_1)_{(m)} = \frac 1
m \left( {\sqrt 2}R ( z \th^{\ou} - \ov{z} \th^1 )  + A_{11,\ou} ( {\sqrt 2}z \th
- 2 t \th^1) - A_{\ou \ou, 1} ( {\sqrt 2}\ov{z} \th - 2 t \th^{\ou})
\right)_{(m)}$,

\noindent $m \geq 2$.
\end{pro}

\begin{df} Given a three dimensional pseudohermitian manifold
$(M,\th)$ we define a real symmetric tensor $Q$ as
$$
  Q = Q_{jk} \th^j \otimes \th^k, \qquad j,k \in \left\{ 0, 1, \ou \right\}
$$
whose components with respect to any admissible coframe are given by
$$
  Q_{11} = \overline{Q_{\ou \ou}} = 3 i A_{11}; \qquad \quad Q_{1 \ou} = Q_{\ou 1}
  = h_{1 \ou}R;
$$
$$
  Q_{01} = Q_{10} = \overline{Q_{0 \ou}} = \overline{Q_{\ou 0}} = 4 A_{11,}^{\;\;\;\;\;1}
  + i R_{,1}; \qquad \quad Q_{00} = 16 \hbox{Im } A_{11,}^{\;\; \;\;\; 11} - 2 \D_b R.
$$
\end{df}

\

\noindent We have then the following result, see page 315 in \cite{JL},
Theorem 3.1.

\begin{pro}\label{p:CRnormcoord} Suppose $M$ is a strictly pseudoconvex pseudohermitian manifold of
dimension 3 and let $q \in M$. Then for any integer $N \geq 2$ there exists a choice of
contact form $\th$ such that all symmetrized covariant derivatives
of $Q$ with total order less or equal than $N$ vanish at $q$, that
is
\begin{equation}\label{eq:sqv}
    Q_{\langle jk, L \rangle} = 0 \qquad \qquad \hbox{ if }
   \quad \mathbb{O}(jkL) \leq N.
\end{equation}
\end{pro}

\begin{rem}
(a) For a multi index $L = (l_1, \dots, l_s)$ we count its order as
$$
  \mathbb{O}(J) = \mathbb{O}(j_1) + \cdots + \mathbb{O}(j_s),
$$
where $\mathbb{O}(1) = \mathbb{O}(\ov{1}) = 1$ and where $\mathbb{O}(0) = 2$.

(b) The symmetrized covariant derivatives are defined by
$$
  Q_{\langle L \rangle} = \frac{1}{s!} \sum_{\s \in \mathbb{S}_s} Q_{\s L}; \qquad
  \s L = \left( l_{\s(1)}, \dots, l_{\s(s)} \right).
$$
\end{rem}

\

\noindent Let us apply Proposition \ref{p:CRnormcoord} for $N=4$
and derive some consequences. At order $2$, at $q$ we have
\begin{equation*}
0=Q_{11}=3iA_{11};\qquad \quad 0=Q_{1\overline{1}}=R.
\end{equation*}

At order $3$ we have
\begin{equation*}
0=Q_{01}=4A_{11,}^{\;\;\;\;\;1}+iR_{,1},
\end{equation*}%
and
\begin{equation*}
0=3!Q_{\langle 11,\overline{1}\rangle }=Q_{11,\overline{1}}+Q_{1\overline{1}%
,1}+Q_{\overline{1}1,1}=3iA_{11,\overline{1}}+R_{,1}+R_{,1}=3iA_{11,%
\overline{1}}+2R_{,1}.
\end{equation*}%
These two equations imply that
\begin{equation*}
A_{11,\overline{1}}=R_{,1}=0\qquad \quad \hbox{ at }q.
\end{equation*}%
We also have
\begin{equation*}
0=Q_{\langle 11,1\rangle }=Q_{11,1}=3iA_{11,1}.
\end{equation*}

At order $4$ we find
\begin{equation*}
0=Q_{00}=16\hbox{Im }A_{11,}^{\;\;\;\;\;11}-2\D_{b}R
\end{equation*}%
and
\begin{equation*}
0=Q_{1\overline{1},0}+Q_{01,\overline{1}}+Q_{0\overline{1},1}=R_{,0}+\left(
4A_{11,}^{\;\;\;\;\;1}+iR_{,1}\right) _{,\overline{1}}+\left( 4A_{\overline{1%
}\overline{1},}^{\;\;\;\;\;\overline{1}}-iR_{,\overline{1}}\right) _{,%
\overline{1}}.
\end{equation*}%
This quantity is equal to
\begin{eqnarray*}
0 &=&R_{,0}+4A_{11,\;\overline{1}}^{\;\;\;\;\;\;\;1}+iR_{,1\overline{1}}+4A_{%
\overline{1}\overline{1},\;1}^{\;\;\;\;\;\;\;\overline{1}}-iR_{,\overline{1}%
1} \\
&=&R_{,0}+8\hbox{Re }A_{11,}^{\;\;\;\;\;11}+i(iR_{,0})=8\hbox{Re }%
A_{11,}^{\;\;\;\;\;11}.
\end{eqnarray*}%
We also have that
\begin{equation}
0=Q_{11,0}+Q_{10,1}+Q_{01,1}=3iA_{11,0}+2\left(
4A_{11,}^{\;\;\;\;\;1}+iR_{,1}\right) _{,1}=3iA_{11,0}+8A_{11,\bar{1}%
}^{\;\;\;\;\;\;\;\bar{1}}+2iR_{,11};  \label{eq:d1}
\end{equation}%
\begin{equation*}
0=Q_{\langle 11,\overline{1}\overline{1}\rangle }=4Q_{11,\overline{1}%
\overline{1}}+4Q_{1\overline{1},1\overline{1}}+4Q_{1\overline{1},\overline{1}%
1}+4Q_{\overline{1}\overline{1},11}+4Q_{\overline{1}1,\overline{1}1}+4Q_{%
\overline{1}1,1\overline{1}},
\end{equation*}%
namely
\begin{equation*}
0=3iA_{11,\overline{1}\overline{1}}+R_{,1\overline{1}}+R_{,\overline{1}1}+%
\hbox{conj.}=-6\hbox{Im }A_{11,\overline{1}\overline{1}}-2\D_{b}R.
\end{equation*}%
From these relations we deduce that
\begin{equation*}
\hbox{Im }A_{11,}^{\;\;\;\;\;11}=\D_{b}R=0,
\end{equation*}%
and hence
\begin{equation*}
A_{11,}^{\;\;\;\;\;11}=0\quad \Rightarrow \quad R_{,0}=0\quad \Rightarrow
\quad \Box _{b}R=0.
\end{equation*}%
We also have
\begin{equation*}
24Q_{\langle 11,1\overline{1}\rangle }=6\left( Q_{11,1\overline{1}}+Q_{11,%
\overline{1}1}+Q_{1\overline{1},11}+Q_{\overline{1}1,11}\right) =6\left(
3i(A_{11,1\overline{1}}+A_{11,\overline{1}1})+2R_{,11}\right) ,
\end{equation*}%
and
\begin{equation}
0=3i\left( A_{11,\overline{1}1}+iA_{11,0}+2RA_{11}+A_{11,\overline{1}%
1}\right) +2R_{,11}=-3A_{11,0}+6iA_{11,\overline{1}1}+2R_{,11}.
\label{eq:d2}
\end{equation}%
From \eqref{eq:d1}, \eqref{eq:d2} and the definition of Cartan tensor
$\mathfrak{Q}$ (see page 227 in \cite{CL1}) we have that, at $q$
\begin{equation*}
\left\{
\begin{array}{ll}
3iA_{11,0}+8A_{11,\overline{1}1}+2iR_{,11}=0; &  \\
-3A_{11,0}+6iA_{11,\overline{1}1}+2R_{,11}=0; &  \\
-6A_{11,0}-4iA_{11,\overline{1}1}+R_{,11}=6\mathfrak{Q}_{11}, &
\end{array}%
\right.
\end{equation*}%
where
\begin{equation*}
\mathfrak{Q}_{11}=\frac{1}{6}R_{,11}+\frac{i}{2}RA_{11}-A_{11,0}-\frac{2}{3}%
iA_{11,\overline{1}1}.
\end{equation*}%
Therefore we deduce that, at $q$
\begin{equation*}
A_{11,0}=-\frac{4}{5}\mathfrak{Q}_{11};\qquad A_{11,\overline{1}1}=\frac{12}{%
35}i\mathfrak{Q}_{11};\qquad R_{,11}=-\frac{6}{35}\mathfrak{Q}_{11}.
\end{equation*}
For $N=4$, suppose we have chosen a contact form $\hat{\th }$ such that %
\eqref{eq:sqv} holds true. Then one can check  that
\begin{equation*}
\hat{\o }_{1\,(2)}^{1}=\hat{\o }_{1\,(3)}^{1}=0;\qquad \hat{\th }_{(3)}^{1}=%
\hat{\th }_{(4)}^{1}=0;\qquad \hat{\th }_{(4)}=\hat{\th }_{(5)}=0.
\end{equation*}%
For example we have
\begin{eqnarray*} \hat{\o }_{1\,(2)}^{1}
&=&\frac{1}{2}\left( \sqrt{2}R(z\hat{\th }^{\overline{1}}-\overline{z}\hat{%
\th }^{1})+A_{11,\overline{1}}(\sqrt{2}z\hat{\th }-2t\hat{\th }^{1})-A_{%
\overline{1}\overline{1},1}(\sqrt{2}\overline{z}\hat{\th }-2t\hat{\th }^{%
\overline{1}})\right) _{(2)} \\
&=&\frac{1}{2}\left( \sqrt{2}z(R\hat{\th }^{\overline{1}})_{(1)}-\sqrt{2}%
\overline{z}(R\hat{\th }^{1})_{(1)}+\sqrt{2}z(A_{11,\overline{1}}\hat{\th }%
)_{(1)}-2t(A_{11,\overline{1}}\hat{\th }^{1})_{(0)} \right.
\\ & - & \left. \sqrt{2}\overline{z}(A_{%
\overline{1}\overline{1},1}\hat{\th })_{(1)}+2t(A_{\overline{1}\overline{1}%
,1}\hat{\th }^{\overline{1}})_{(0)}\right) = 0,
\end{eqnarray*}%
and similarly $\hat{\o }_{1\,(3)}^{1}$ $=$ $0.$ Now we have that
\begin{eqnarray*}
\hat{\o }_{1\,(4)}^{1}
&=&\frac{1}{4}\left( \sqrt{2}R(z\hat{\th }^{\overline{1}}-\overline{z}\hat{%
\th }^{1})+A_{11,\overline{1}}(\sqrt{2}z\hat{\th }-2t\hat{\th }^{1})-A_{%
\overline{1}\overline{1},1}(\sqrt{2}\overline{z}\hat{\th }-2t\hat{\th }^{%
\overline{1}})\right) _{(4)} \\
&=&\frac{1}{4}\left( \sqrt{2}z(R\hat{\th }^{\overline{1}})_{(3)}-\sqrt{2}%
\overline{z}(R\hat{\th }^{1})_{(3)}+\sqrt{2}z(A_{11,\overline{1}}\hat{\th }%
)_{(3)}-2t(A_{11,\overline{1}}\hat{\th }^{1})_{(2)} \right. \\
& - & \left. \sqrt{2}\overline{z}(A_{%
\overline{1}\overline{1},1}\hat{\th })_{(3)}+2t(A_{\overline{1}\overline{1}%
,1}\hat{\th }^{\overline{1}})_{(2)}\right)  \\
&=&\frac{1}{4}\left( 2R_{(2)}zd\overline{z}-2R_{(2)}\overline{z}dz+\sqrt{2}%
(A_{11,\overline{1}})_{(1)}z\theta _{0}-2(A_{11,\overline{1}})_{(1)}t\sqrt{2}%
dz \right.
\\ & - & \left. \sqrt{2}(A_{\overline{1}\overline{1},1})_{(1)}\overline{z}\theta
_{0}+2(A_{\overline{1}\overline{1},1})_{(1)}t\sqrt{2}d\overline{z}\right)  \\
&=&\frac{1}{4}\left( -\sqrt{2}R_{(2)}\overline{z}-2(A_{11,\overline{1}%
})_{(1)}t\right) \sqrt{2}dz+\frac{1}{4}\left( \sqrt{2}R_{(2)}z+2(A_{%
\overline{1}\overline{1},1})_{(1)}t\right) \sqrt{2}d\overline{z}
\\ & + & \frac{\sqrt{%
2}}{4}\left( (A_{11,\overline{1}})_{(1)}z-(A_{\overline{1}\overline{1}%
,1})_{(1)}\overline{z}\right) \theta _{0},
\end{eqnarray*}%
so we deduce
\begin{eqnarray*}
\hat{\o }_{1\,(4)}^{1} &=&\left( -R_{,11}(q)z^{2}\overline{z}-R_{,\overline{1%
}\overline{1}}(q)\overline{z}^{3}-A_{11,\overline{1}1}(q)zt\right) dz+\left(
R_{,11}(q)z^{3}+R_{,\overline{1}\overline{1}}(q)\overline{z}^{2}z+A_{%
\overline{1}\overline{1},1\overline{1}}(q)\overline{z}t\right) d\overline{z}
\\
&+&\frac{1}{2}\left( A_{11,\overline{1}1}(q)z^{2}-A_{\overline{1}\overline{1}%
,1\overline{1}}(q)\overline{z}^{2}\right) \theta _{0}.
\end{eqnarray*}%
Hence we conclude
\begin{equation*}
\hat{\o }_{1}^{1}=O(\rho ^{3})dz+O(\rho ^{3})d\overline{z}+O(\rho
^{2})\theta _{0}.
\end{equation*}
We also have
\begin{eqnarray*}
\hat{\th }_{(5)}^{1} &=&\frac{1}{5}\left( \sqrt{2}z\hat{\o }_{1}^{1}+2tA_{%
\overline{1}\overline{1}}\hat{\th }^{\bar{1}}-\sqrt{2}\bar{z}A_{\overline{1}%
\overline{1}}\th \right) _{(5)}=\frac{1}{5}\left( \sqrt{2}z(\hat{\o }%
_{1}^{1})_{(4)}+2t(A_{\overline{1}\overline{1}}\hat{\th }^{\overline{1}%
})_{(3)}-\sqrt{2}\overline{z}(A_{\overline{1}\overline{1}}\th )_{(4)}\right)
\\
&=&\frac{1}{5}\left( \sqrt{2}z(\hat{\o }_{1}^{1})_{(4)}+2t(A_{\overline{1}%
\overline{1}})_{(2)}\sqrt{2}d\overline{z}-\sqrt{2}\overline{z}(A_{\overline{1%
}\overline{1}})_{(2)}\theta _{0}\right) ,
\end{eqnarray*}%
and therefore we get
\begin{eqnarray*}
 & & \hat{\th }_{(5)}^{1} = \frac{1}{20}\left( -4R_{,11}(q)z^{3}\overline{z}%
-4R_{,\overline{1}\overline{1}}(q)\overline{z}^{3}z-4A_{11,\overline{1}%
1}(q)z^{2}t\right) \sqrt{2}dz \\
&+&\left[ \frac{1}{5}\left( R_{,11}(q)z^{4}+R_{,\overline{1}\overline{1}%
}(q)z^{2}\overline{z}^{2}+A_{\overline{1}\overline{1},1\overline{1}}(q)z%
\overline{z}t\right) +\frac{2}{5}A_{\overline{1}\overline{1},0}(q)t^{2}+%
\frac{4}{5}(2A_{\overline{1}\overline{1},1\overline{1}}(q)-iA_{\overline{1}%
\overline{1},0}(q))z\overline{z}t\right] \sqrt{2}d\overline{z} \\
&+&\left[ \frac{\sqrt{2}}{10}\left( A_{11,\overline{1}1}(q)z^{3}-A_{%
\overline{1}\overline{1},1\overline{1}}(q)z\overline{z}^{2}\right) -\frac{%
\sqrt{2}}{5}A_{\overline{1}\overline{1},0}(q)\overline{z}t-\frac{2\sqrt{2}}{5}%
\left( 2A_{\overline{1}\overline{1},1\overline{1}}(q)-iA_{\overline{1}%
\overline{1},0}(q)\right) z\overline{z}^{2}\right] \theta _{0},
\end{eqnarray*}%
which implies
\begin{equation*}
\hat{\th }^{1}=\left( 1+O(\rho ^{4})\right) \sqrt{2}dz+O(\rho ^{4})d%
\overline{z}+O(\rho ^{3})\theta _{0}.
\end{equation*}%
We also have
\begin{equation*}
\hat{\th }_{(6)}=\frac{1}{6}\left( \sqrt{2}iz\hat{\th }^{\overline{1}}-\sqrt{%
2}i\overline{z}\hat{\th }^{1}\right) _{(6)}=\frac{\sqrt{2}}{6}\left( iz\hat{%
\th }_{(5)}^{\overline{1}}-i\overline{z}\hat{\th }_{(5)}^{1}\right) ,
\end{equation*}%
and hence
\begin{equation*}
\hat{\th }=\left( 1+O(\rho ^{4})\right) \theta _{0}+O(\rho ^{5})dz+O(\rho
^{5})d\overline{z}.
\end{equation*}

Let us now try to understand the behavior of the dual vectors. Let
\begin{equation*}
\hat{W}_{1}=(1+a)\mathbb{Z}_{1}+b\mathbb{Z}_{\overline{1}}+c\frac{\partial }{%
\partial t},
\end{equation*}%
where $\mathbb{Z}_{1}=\frac{1}{\sqrt{2}}(\frac{\partial }{\partial z}+i%
\overline{z}\frac{\partial }{\partial t})$. From these formulas, if we set
\begin{equation*}
\hat{\th }=\theta _{0}+\sum_{m\geq 6}\hat{\th }_{(m)};\qquad \hat{\th }^{1}=%
\sqrt{2}dz+\sum_{n\geq 5}\hat{\th }_{(n)}^{1};\qquad \hat{\th }^{\overline{1}%
}=\sqrt{2}d\overline{z}+\sum_{n\geq 5}\hat{\th }_{(n)}^{\overline{1}},
\end{equation*}%
we have that
\begin{eqnarray*}
0 &=&\hat{\th }(\hat{W}_{1})=\left( \theta _{0}+\sum_{m\geq 6}\hat{\th }%
_{(m)}\right) \left( (1+a)\mathbb{Z}_{1}+b\mathbb{Z}_{\overline{1}}+c\frac{%
\partial }{\partial t}\right)  \\
&=&(1+a)\theta _{0}(\mathbb{Z}_{1})+b\theta _{0}(\mathbb{Z}_{\overline{1}%
})+c\theta _{0}(\frac{\partial }{\partial t})+(1+a)\sum_{m\geq 6}\hat{\th }%
_{(m)}(\mathbb{Z}_{1})+b\sum_{m\geq 6}\hat{\th }_{(m)}(\mathbb{Z}_{\overline{%
1}})+c\sum_{m\geq 6}\hat{\th }_{(m)}(\frac{\partial }{\partial t}) \\
&=&c+O(\rho ^{5})+aO(\rho ^{5})+bO(\rho ^{5})+cO(\rho ^{4})
\end{eqnarray*}%
and hence we find that $c=O(\rho ^{5})$.

Similarly
\begin{equation*}
1=\hat{\th }^{1}(\hat{W}_{1})\quad \Rightarrow \quad a=O(\rho ^{4});\qquad
\qquad 0=\hat{\th }^{\overline{1}}(\hat{W}_{1})\quad \Rightarrow \quad
b=O(\rho ^{4}).
\end{equation*}%
Also, if we set
\begin{equation*}
\hat{T}=(1+a)\frac{\partial }{\partial t}+b\mathbb{Z}_{1}+c\mathbb{Z}_{%
\overline{1}},
\end{equation*}%
then from
\begin{equation*}
1=\hat{\th }(\hat{T});\qquad 0=
\hat{\th }^{1}(\hat{T});\qquad 0=\hat{\th }^{%
\overline{1}}(\hat{T})
\end{equation*}%
we deduce that $a=O(\rho ^{4}),b=O(\rho ^{3}),c=O(\rho ^{3})$, and hence we
obtain
\begin{equation*}
\hat{T}=\left( 1+O(\rho ^{4})\right) \frac{\partial }{\partial t}+O(\rho
^{3})\mathbb{Z}_{1}+O(\rho ^{3})\mathbb{Z}_{\overline{1}}.
\end{equation*}%
In conclusion, we arrive to the following result.

\begin{pro}\label{p:CRcoord} In CR normal coordinates, for $N = 4$ we can choose a
contact form $\hat{\th}$ such that in pseudohermitian normal
coordinates $(z, t)$ we have that
$$
  \hat{\th} = \left( 1 + O(\rho^4) \right) \theta_0 + O(\rho^5) dz
  + O(\rho^5) d\ov{z}; \qquad \hat{\th}^1 = \left( 1 + O(\rho^4)
  \right) {\sqrt 2}dz + O(\rho^4) d\ov{z} + O(\rho^3) \theta_0;
$$
$$
  \hat{\o}^1_1 = O(\rho^3) d z + O(\rho^3) d\ov{z} + O(\rho^2) \theta_0;
$$
$$
  \hat{W}_1 = \left( 1 + O(\rho^4) \right) \mathbb{Z}_1 + O(\rho^4)
  \mathbb{Z}_{\ou} + O(\rho^5) \frac{\pa}{\pa t}; \qquad \quad \hat{T}
  = \left( 1 + O(\rho^4) \right) \frac{\pa}{\pa t} +
  O(\rho^3) \mathbb{Z}_1 + O(\rho^3) \mathbb{Z}_{\ou},
$$
where
$$
 \theta_0 = dt + i z d\ov{z} - i \ov{z} dz; \qquad \mathbb{Z}_1 =
 \frac{1}{{\sqrt 2}} \left( \frac{\pa}{\pa z} + i \ov{z} \frac{\pa}{\pa t}
 \right); \qquad \rho^4 = t^2 + |z|^4.
$$
\end{pro}

\noindent The above computation concluding the proof of  Proposition 6.5 was first
done by H.-L. Chiu in the summer of 2007 when he visited J.-H. Cheng.

\end{document}